\documentclass[12pt]{amsart}
\begin{document}

\input psfig
\let\Bbb\mathbb
\let\goth\mathfrak
\let\cal\mathcal
\makeatletter
\def\part{\@startsection{part}{0}%
  \z@{\linespacing\@plus\linespacing}{.5\linespacing}%
  {\large\bfseries\centering}}
\def\subsection{\@startsection{subsection}{2}%
  \z@{.5\linespacing\@plus.7\linespacing}{8pt plus 4pt}%
  {\normalfont\itshape}}
\def\subsubsection{\@startsection{subsubsection}{3}%
  \z@{.5\linespacing\@plus.7\linespacing}{8pt plus 4pt}%
  {\normalfont\itshape}}
\makeatother

\let\appsection\section 
\let\bppsection\section 

% Copied from plain tex
\def\eqalign#1{\null\,\vcenter{\openup\jot
  \ialign{\strut\hfil$\displaystyle{##}$&$\displaystyle{{}##}$\hfil
      \crcr#1\crcr}}\,}
\def\cases#1{\left\{\,\vcenter{\normalbaselines
    \ialign{$##\hfil$&\quad##\hfil\crcr#1\crcr}}\right.}
\def\matrix#1{\null\,\vcenter{\normalbaselines
    \ialign{\hfil$##$\hfil&&\quad\hfil$##$\hfil\crcr
      \mathstrut\crcr\noalign{\kern-\baselineskip}
      #1\crcr\mathstrut\crcr\noalign{\kern-\baselineskip}}}\,}
\def\pmatrix#1{\left(\matrix{#1}\right)}

% Cancel out pedantry
\makeatletter
\let\atopwithdelims=\@@atopwithdelims    
\let\overwithdelims=\@@overwithdelims    
\let\abovewithdelims=\@@abovewithdelims    
\let\atop=\@@atop
\let\above=\@@above
\let\over\@@over
\makeatother
% ********************* DEFINITIONS*************************

\def\smallskip{\vskip\smallskipamount\noindent}
\def\medskip{\vskip\medskipamount\noindent}
\def\bigskip{\vskip\bigskipamount\noindent}

\let\bigbreak\bigskip
\let\medbreak\medskip
\let\smallbreak\smallskip

\def\C{\ensuremath{\Bbb C}}
\def\F{\ensuremath{\Bbb F}}
\def\N{\ensuremath{\Bbb N}}
\def\Z{\ensuremath{\Bbb Z}}
\def\R{\ensuremath{\Bbb R}}
\def\CG{\mathrel{\C\mkern-.5mu G}}
\def\CH{\mathrel{\C\mkern-.5mu H}}

\def\Card{\mathrm{Card}}
\def\chr{\mathrm{char}}
\def\End{\mathrm{End}}
\def\Hom{\mathrm{Hom}}  %without mathop will place subscripts correctly
\def\id{\mathrm{id}}
\def\im{\mathrm{im}}
\def\ind{\mathrm{ind}}
\def\inv{\mathrm{inv}}
\def\rad{\mathrm{rad}}
\def\span{\mathrm{span}}
\def\sgn{\mathrm{sgn}}
\def\supp{\mathrm{supp}}
\def\trace{\mathrm{trace}}
\def\Tr{\mathrm{Tr}}
\def\wt{\mathrm{wt}}

\def\coeff{{\big\vert}}
\def\One{{\bf 1}}

\def\longrightmapsto{\mapstochar\longrightarrow}

%******************* MATHEMATICAL LABELS **************************

\newtheorem{Prop}{Proposition}[section]
\newtheorem{Lemma}[Prop]{Lemma}
\newtheorem{Cor}[Prop]{Corollary}
\newtheorem{Thm}[Prop]{Theorem}

\theoremstyle{definition}
\newtheorem{Remark}[Prop]{Remark}

\def\prop{\begin{Prop}}
\def\lemma{\begin{Lemma}}
\def\remark{\begin{Remark}}
\def\cor{\begin{Cor}}
\def\appcor{\begin{Cor}}
\def\thm{\begin{Thm}}
\def\appthm{\begin{Thm}}
\def\bppthm{\begin{Thm}}

\def\endthm{\end{Thm}}
\def\thmend{\end{Thm}}
\def\endlemma{\end{Lemma}}
\def\endremark{\end{Remark}}
\def\endprop{\end{Prop}}
\def\endcor{\end{Cor}}

\def\pf{\begin{proof}}
\def\endpf{\end{proof}}
\def\pfend{\end{proof}}
\def\note{\smallbreak\noindent{Note: \enspace}}

\def\formula{\refstepcounter{Prop}\eqno(\theProp)}
\let\appformula=\formula
\def\formulano{\refstepcounter{Prop}(\theProp)}
\def\lformula{\refstepcounter{Prop}\leqno(\theProp)}

%************Commutative diagrams**********************

\def\mapright#1{\smash{\mathop
	{\longrightarrow}\limits^{#1}}}

\def\mapleftright#1{\smash{\mathop
	{\longleftrightarrow}\limits^{#1}}}

\def\mapsrightto#1{\smash{\mathop
	{\longmapsto}\limits^{#1}}}

\def\mapleft#1{\smash{
   \mathop{\longleftarrow}\limits^{#1}}}

\def\mapdown#1{\Big\downarrow
   \rlap{$\vcenter{\hbox{$\scriptstyle#1$}}$}}

\def\lmapdown#1{{\hbox{$\scriptstyle#1$}}
   \llap {$\vcenter{\hbox{\Big\downarrow}}$} }

\def\mapup#1{\Big\uparrow
   \rlap{$\vcenter{\hbox{$\scriptstyle#1$}}$}}

%********************* Post Script **********************

%\def\tableau{ {\psfig{figure=tableau.ps,height=1.5in}} }
%\def\stdtab{ {\psfig{figure=stdtab.ps,height=.75in}} }
%\def\ynglat{ {\psfig{figure=ynglat.ps,height=4.5in}} }

\def\shape{{\psfig{figure=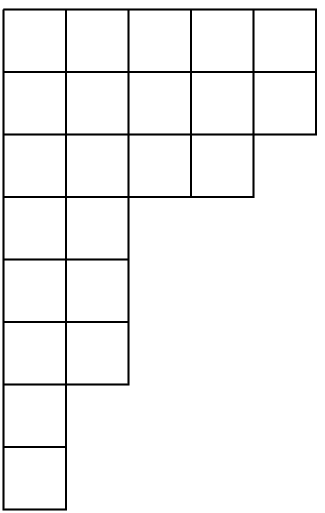,height=1.5in}}}
\def\contents{{\psfig{figure=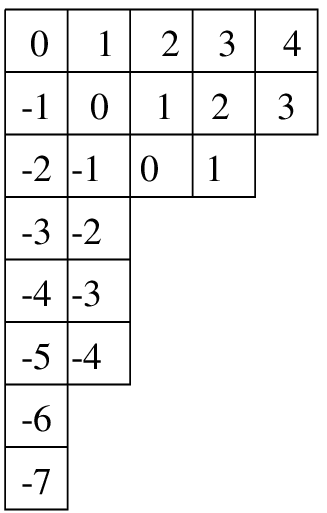,height=1.5in}}}
\def\hooks{{\psfig{figure=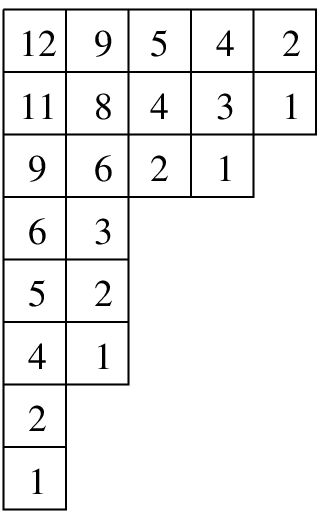,height=1.5in}}}
\def\skewshape{{\psfig{figure=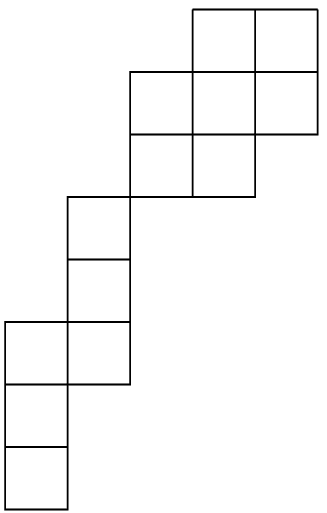,height=1.5in}}}
\def\notlp{{\psfig{figure=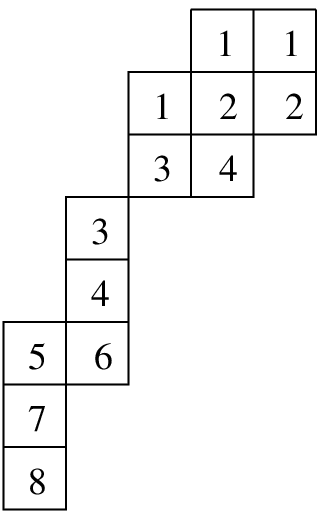,height=1.5in}}}
\def\latperm{{\psfig{figure=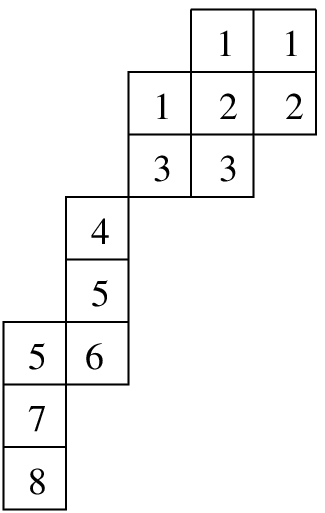,height=1.5in}}}
\def\stdtab{{\psfig{figure=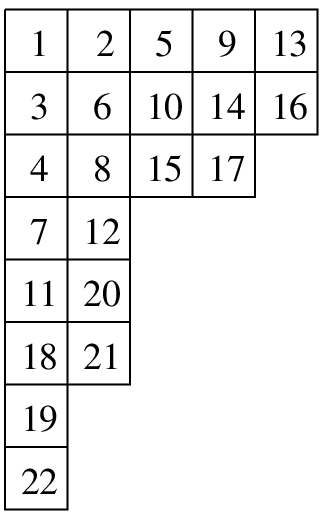,height=1.5in}}}
\def\colsttab{{\psfig{figure=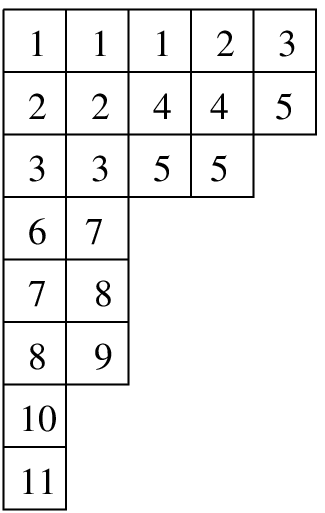,height=1.5in}}}
\def\border{{\psfig{figure=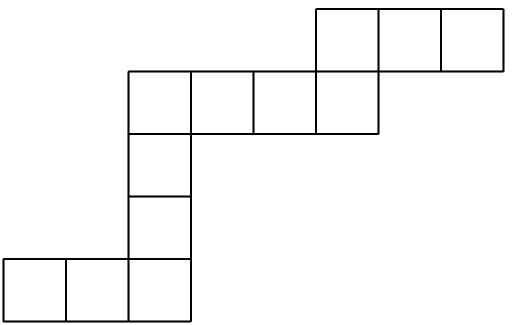,height=.75in}}}

\def\intweight{{\psfig{figure=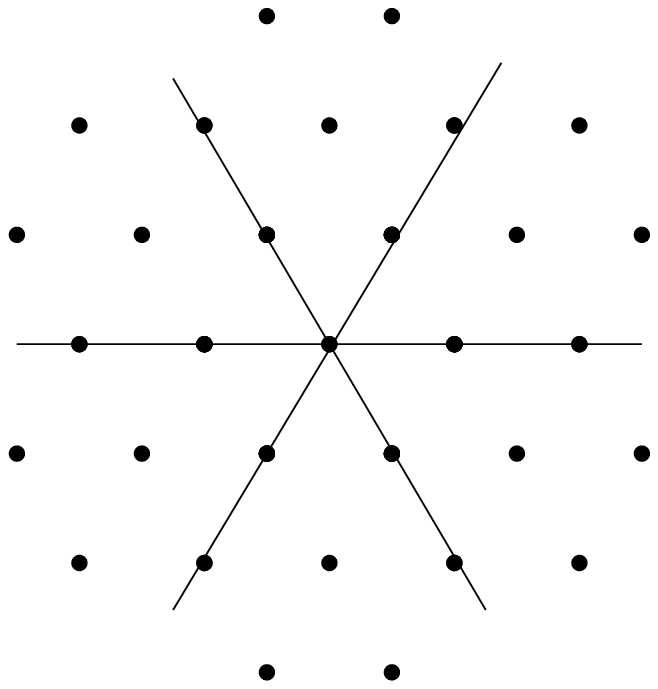,height=2in}}}
\def\cone{{\psfig{figure=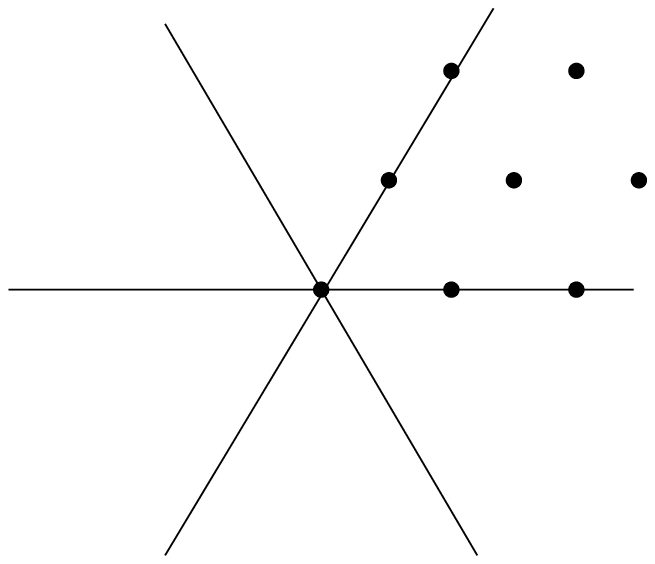,height=2in}}}
\def\piecewiselin{{\psfig{figure=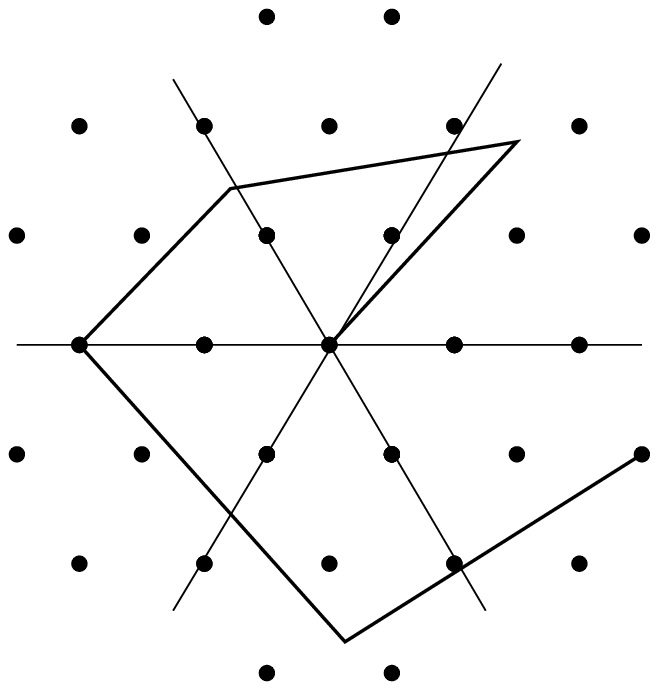,height=2in}}}
\def\root{{\psfig{figure=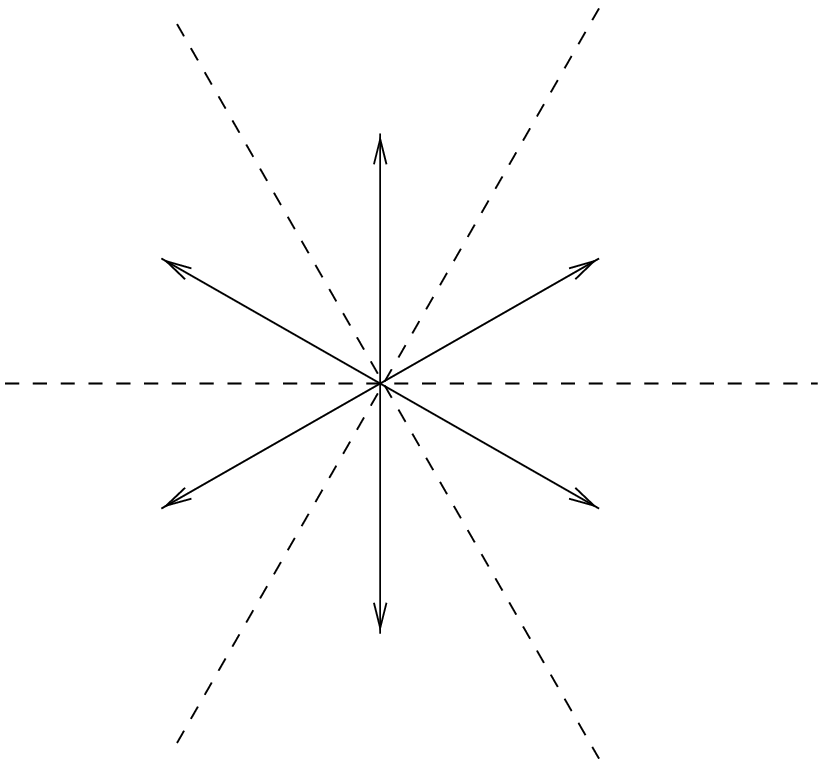,height=2in}}}
\def\chamber{{\psfig{figure=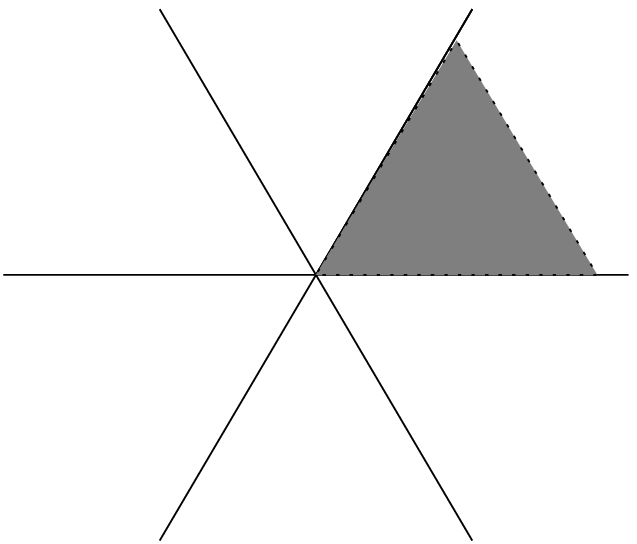,height=2in}}}
\def\hyparr{{\psfig{figure=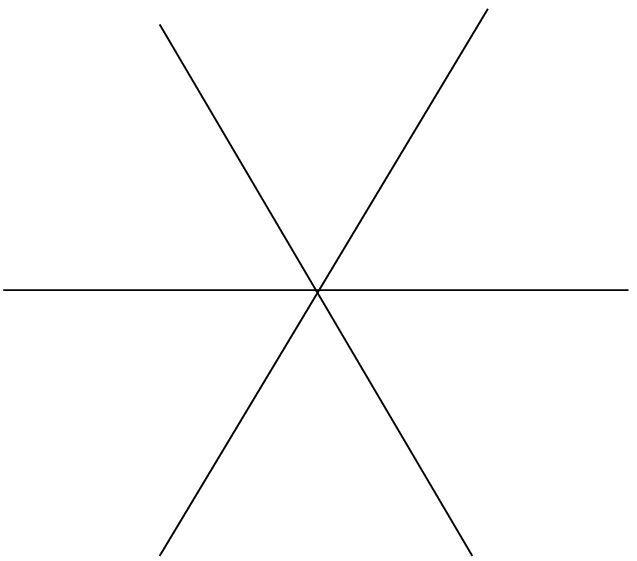,height=2in}}}
\def\path{{\psfig{figure=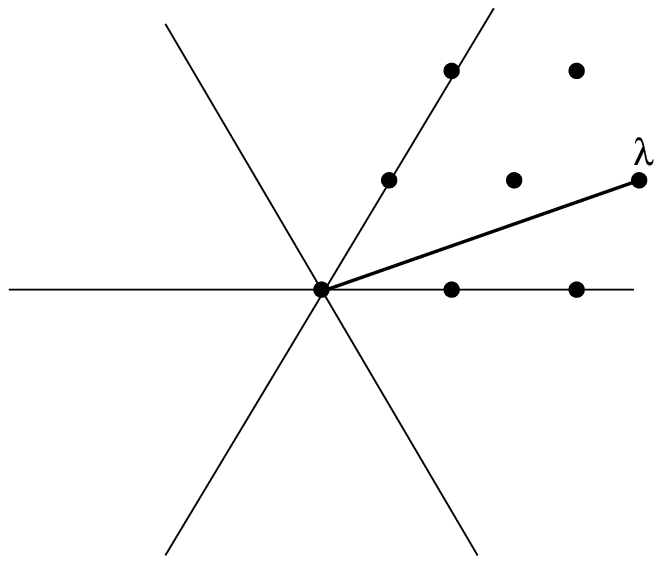,height=2in}}}

\def\TLkdiagram{{\psfig{figure=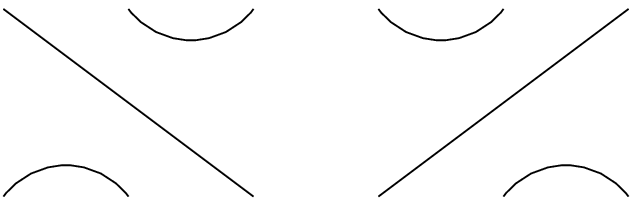,height=.4in}}}
\def\idkdhAa{\lower.2in\hbox{\psfig{figure=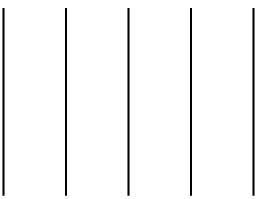,height=.4in}} }
\def\idkdhAb{\lower.2in\hbox{\psfig{figure=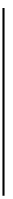,height=.4in}} }
\def\idkdhBa{\lower.2in\hbox{\psfig{figure=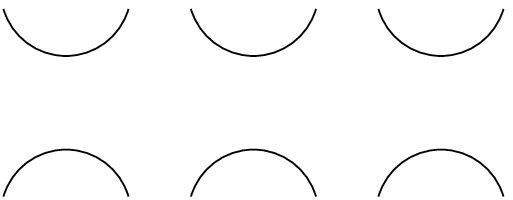,height=.4in}} }
\def\idkdhBb{\lower.2in\hbox{\psfig{figure=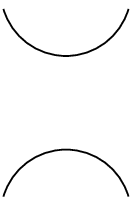,height=.4in}} }
\def\gammuua{\lower.2in\hbox{\psfig{figure=gammu1a.ps,height=.4in}} }
\def\gammuub{\lower.2in\hbox{\psfig{figure=gammu1b.ps,height=.4in}} }
\def\gammud{\lower.2in\hbox{\psfig{figure=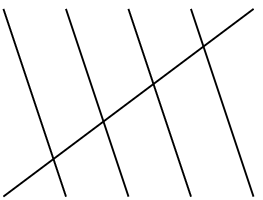,height=.4in}} }
\def\gammut{\lower.2in\hbox{\psfig{figure=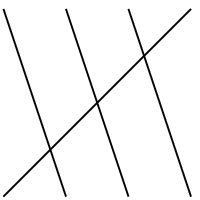,height=.4in}} }
\def\gammuq{\lower.2in\hbox{\psfig{figure=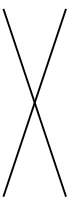,height=.4in}} }

\def\perm{\lower.3in\hbox{\psfig{figure=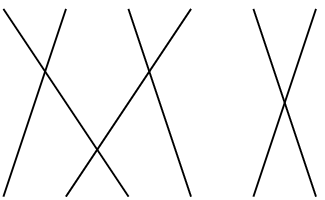,height=.6in}} }
\def\Tw{\lower.3in\hbox{\psfig{figure=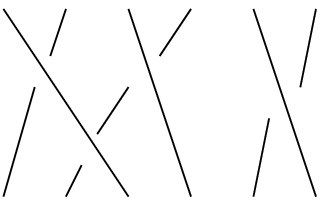,height=.6in}}}
\def\poscross{\lower.2in\hbox{\psfig{figure=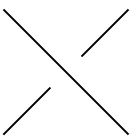,height=.4in}} }
\def\negcross{\lower.2in\hbox{\psfig{figure=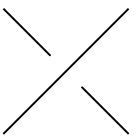,height=.4in}} }
\def\idcross{\lower.2in\hbox{\psfig{figure=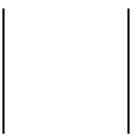,height=.4in}} }

\def\stanglesa{\lower.3in\hbox{\psfig{figure=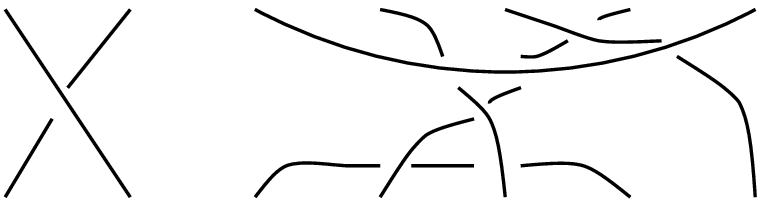,height=.6in}} }
\def\stanglesb{\lower.3in\hbox{\psfig{figure=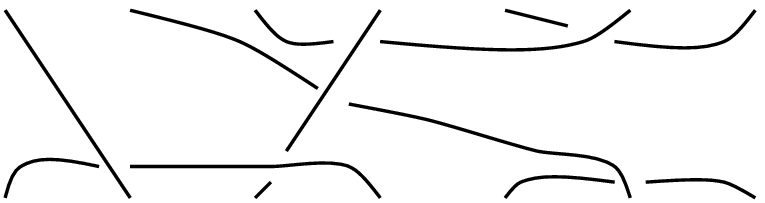,height=.6in}} }
\def\brauertang{\lower.3in\hbox{\psfig{figure=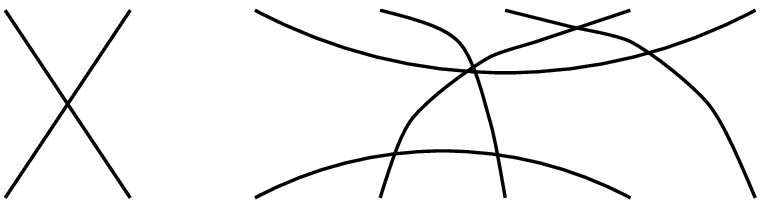,height=.6in}} }
\def\tangleida{\lower.2in\hbox{\psfig{figure=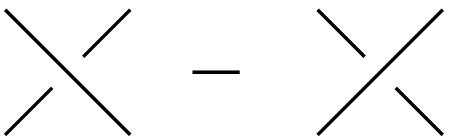,height=.4in}} }
\def\tangleidb{\lower.2in\hbox{\psfig{figure=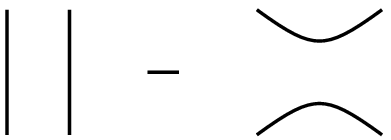,height=.4in}} }
\def\Tda{\lower.2in\hbox{\psfig{figure=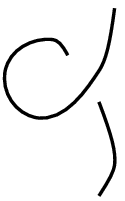,height=.4in}} }
\def\Tdb{\lower.2in\hbox{\psfig{figure=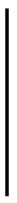,height=.4in}} }
\def\Tdaa{\lower.2in\hbox{\psfig{figure=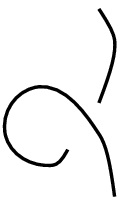,height=.4in}} }
\def\Tt{\lower.1in\hbox{\psfig{figure=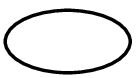,height=.2in}} }
\def\BMWa{\lower.2in\hbox{\psfig{figure=BMWa.ps,width=.375in}} }
\def\BMWb{\lower.2in\hbox{\psfig{figure=BMWb.ps,width=1in}} }
\def\BMWc{\lower.2in\hbox{\psfig{figure=BMWc.ps,width=.70in}} }

\def\bdone{{\lower.25in\hbox{\psfig{figure=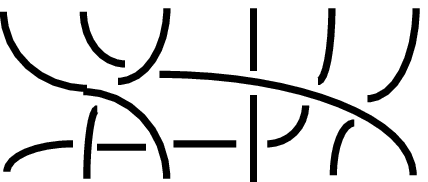,height=.5in}}}}
\def\bdtwo{{\lower.25in\hbox{\psfig{figure=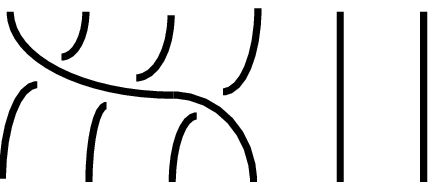,height=.5in}}}}
\def\prodbd{{\lower.5in\hbox{\psfig{figure=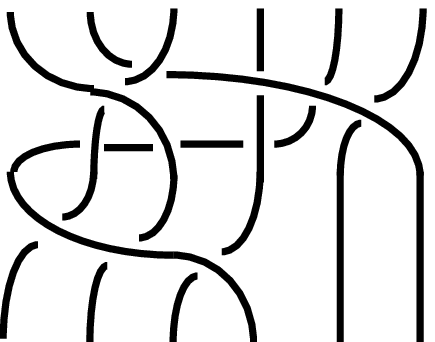,height=1in}}}}

\def\dcomp{\lower.375in\hbox{\psfig{figure=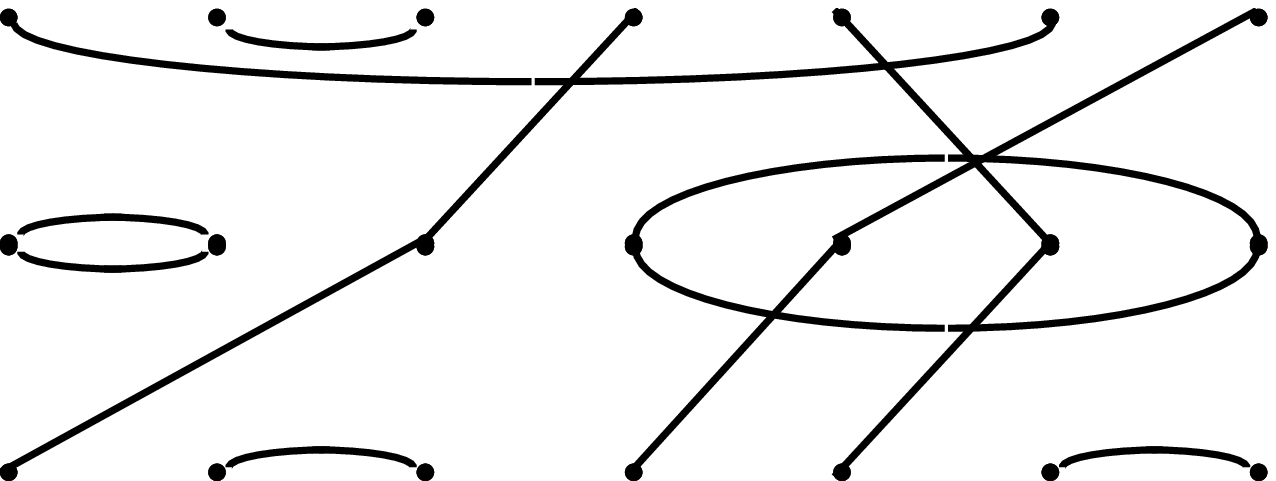,height=.75in}} }
\def\dprod{\lower.2in\hbox{\psfig{figure=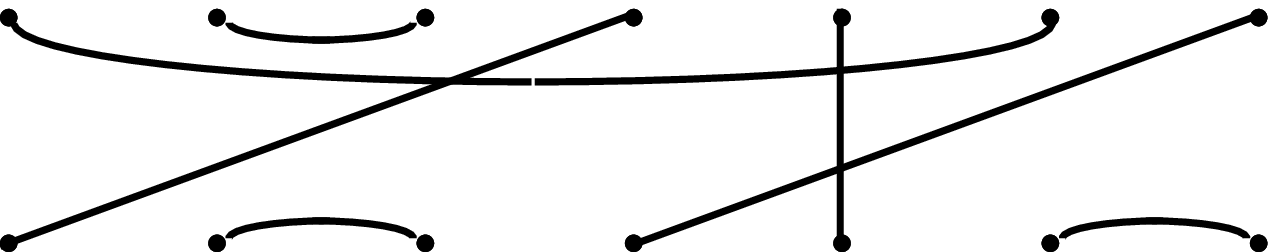,height=.4in}} }
\def\donediag{\lower.2in\hbox{\psfig{figure=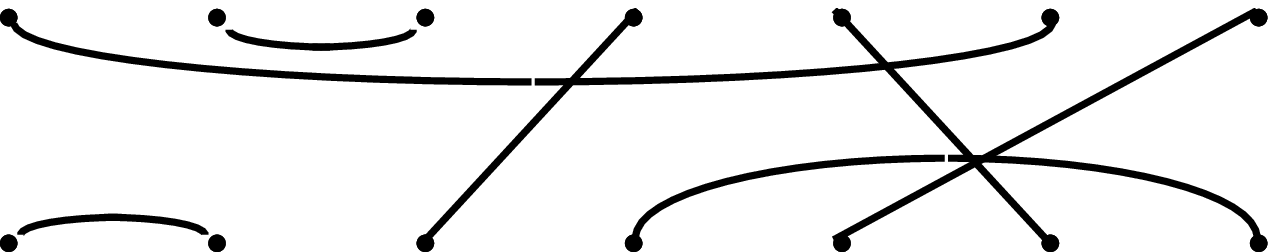,height=.4in,width=1.7in}} }
\def\dtwodiag{\lower.2in\hbox{\psfig{figure=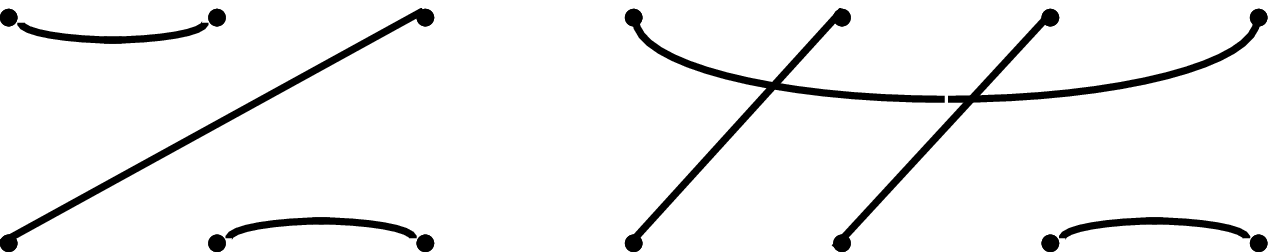,height=.4in,width=1.7in}} }
\def\iddiag{\lower.25in\hbox{\psfig{figure=iddiag.eps,height=.6in}} }
\def\sidiag{\lower.2in\hbox{\psfig{figure=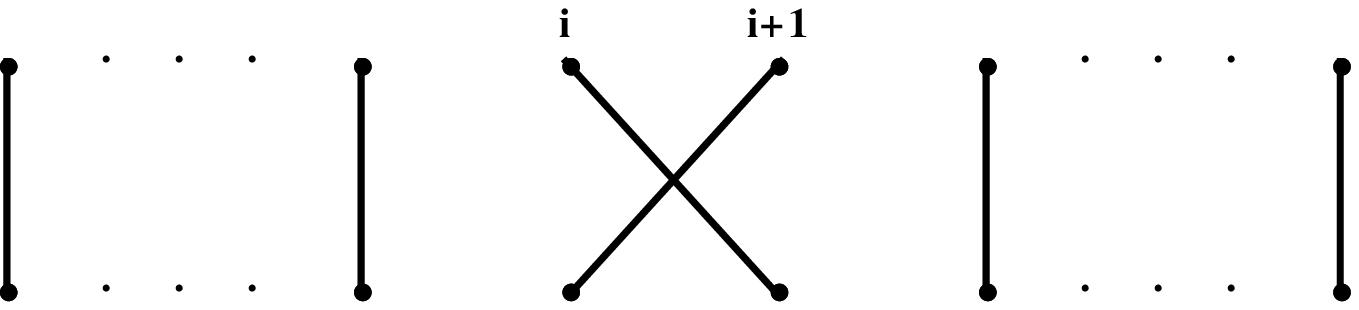,height=.5in}} }
\def\eidiag{\lower.2in\hbox{\psfig{figure=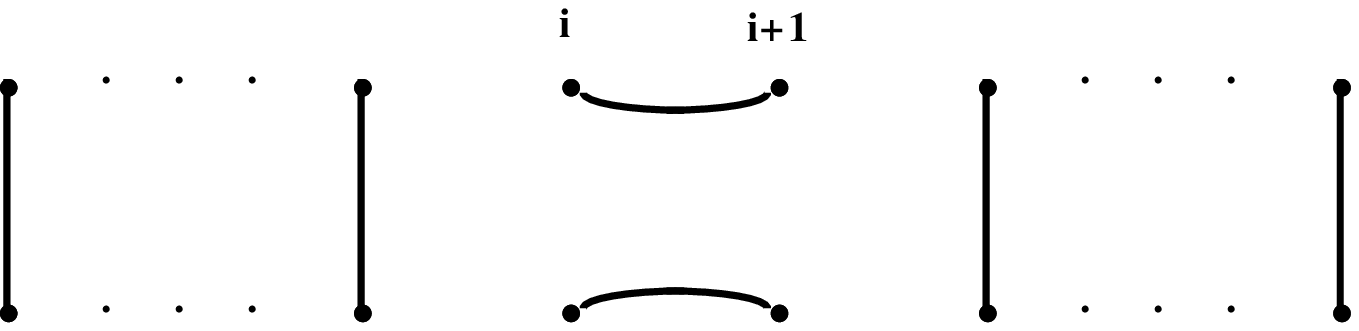,height=.5in}} }
\def\diage{\lower.2in\hbox{\psfig{figure=diage.eps,height=.5in}} }
\def\gammadiag{\lower.2in\hbox{\psfig{figure=gamma.eps,height=.5in}} }

\def\ronea{{\lower.1in\hbox{\psfig{figure=ronea.eps,width=.75in}}}}
\def\roneb{{\psfig{figure=roneb.eps,width=.75in}}}
\def\ronec{{\lower.1in\hbox{\psfig{figure=ronec.eps,width=.75in}}}}
\def\rtwoa{{\lower.15in\hbox{\psfig{figure=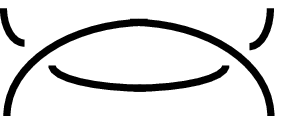,width=.75in}}}}
\def\rtwob{{\lower.15in\hbox{\psfig{figure=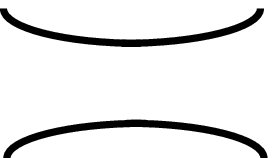,width=.75in}}}}
\def\rtwoc{{\lower.15in\hbox{\psfig{figure=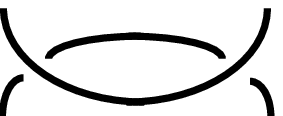,width=.75in}}}}
\def\rthreea{{\lower.25in\hbox{\psfig{figure=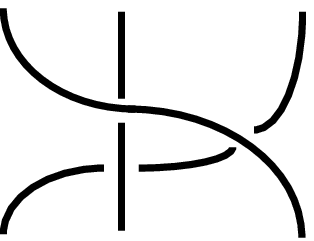,height=.5in}}}}
\def\rthreeb{{\lower.25in\hbox{\psfig{figure=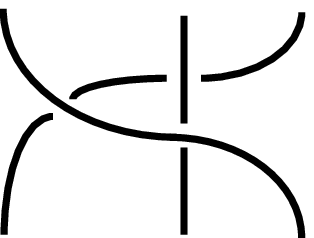,height=.5in}}}}

%**************** TITLE *************************
\vphantom{$ $}  %My kludge to get the first page to move down a bit
\vskip.5truein

\title{Combinatorial Representation Theory}
%\bigbreak
%\centerline{Draft: \today}
%\vskip.4truein
\author{H\'el\`ene Barcelo}
\thanks{\hskip-\parindent
Barcelo was supported in part by National
Science Foundation grant DMS-9510655.
\smallbreak
Ram was supported in part by National
Science Foundation grant DMS-9622985.
\smallbreak
This paper was written while both authors were in residence at MSRI.
We are grateful for the hospitality and financial support of MSRI.}
\address{\hskip-\parindent
H\'el\`ene Barcelo\\
Department of Mathematics\\
Arizona State University\\
Tempe, AZ 85287-1804}
\author{Arun Ram}
\address{\hskip-\parindent
Department of Mathematics\\
Princeton University\\
Princeton, NJ 08544}


\begin{abstract}
We attempt to survey the field of combinatorial
representation theory, describe the main results 
and main questions and give an update of its current status.
We give a personal viewpoint on the field, while remaining aware 
that there is much important and beautiful work 
that we have not been able to mention.
\end{abstract}

\maketitle

\tableofcontents

\section*{Introduction}

In January 1997, during the special year in combinatorics at MSRI,
at a dessert party at H\'el\`ene's house,
Gil Kalai, in his usual fashion, began asking very pointed questions
about exactly what all the combinatorial representation
theorists were doing their research on.
After several unsuccessful attempts at giving answers that Gil would
find satisfactory,
it was decided that some talks should be given
in order to explain to other combinatorialists what the specialty
is about and what its main questions are.

In the end, Arun gave two talks at MSRI in which he tried 
to clear up the situation.  
After the talks several people suggested
that it would be helpful if someone would write a survey article
containing what had been covered in the two talks and including
further interesting details.  After some arm twisting
it was agreed that Arun and 
H\'el\`ene would write such a paper on
combinatorial representation theory.
What follows is our attempt to define the field of combinatorial
representation theory, describe the main results 
and main questions and give an update of its current status.

Of course this is wholly impossible.  Everybody in the field
has their own point of view and their own preferences of
questions and answers.  Furthermore, there is much too
much material in the field to possibly collect it all
in a single article (even conceptually).
We therefore feel that we must
stress the obvious; in this article we give
a personal viewpoint on the field while remaining aware 
that there is much important and beautiful work 
that we have not been able to mention.

On the other hand, we have tried very hard to give a focused
approach and to make something that will be useful 
to both specialists and non specialists, 
for understanding what we do, for learning the concepts
of the field, and for tracking down history and references.
We have chosen to write in an informal style in the hope that
this way we can better convey the conceptual aspects of the field.  
Readers should keep this in mind
and refer to the notes and references and the
appendices when there are questions about the precision in 
definitions and statements of results.
We have included a table of contents at the end of the paper 
which should help with navigation.
Having made these points, and put in a lot of work,
we leave it to you the reader, with
the earnest hope that you find it useful.

We would like to thank the many people in residence at the
special year 1996-97 in combinatorics at MSRI, for their interest,
their suggestions, and for continually encouraging
us to explain and write about the things that we enjoy doing.
We both are extremely indebted to our graduate
advisors, A. Garsia and H. Wenzl, who (already many years ago)
introduced us to and taught us this wonderful field.

\part*{Part I}

% Section 1
\section{What is Combinatorial Representation Theory?}

What do we mean by ``combinatorial representation
theory''?
First and foremost, combinatorial representation theory is 
representation theory.
The adjective ``combinatorial'' will refer to the way in which
we answer representation theoretic questions.
We will discuss this more fully
later.  For the moment let us begin with, what is representation theory?

\section*{\qquad What is representation theory?\qquad\null}

If representation theory is a black box, or a machine,
then the input is an algebra $A$. 
The output of the machine is information about the modules 
for $A$.

\smallskip\noindent
{\bf 
An {\it algebra} is a vector space $A$ over $\C $ with a 
multiplication.
}
\smallskip\noindent
An important example is the case of the group algebra.
\smallskip\noindent
{\bf 
Define the {\it group algebra} of a group  $G$ to be 
$$A = \C G =  {\hbox{$\C$-span}~\{ g\in G\}},$$
so that the elements of $G$ form a basis of $A$.  The multiplication
in the group algebra is inherited from the multiplication in the group.
}
\smallskip\noindent
We want to study the algebra $A$ via its actions on vector spaces.
\smallskip\noindent
{\bf 
An {\it $A$-module} is a finite dimensional vector space $M$ over $\C$
with an $A$ action. 
}
\smallskip\noindent
See Appendix A1 for a complete definition.
We shall use the words module and {\it representation} interchangeably.
Representation theorists are always trying to break up modules into
pieces.
\smallskip\noindent
{\bf 
An $A$-module $M$ is {\it indecomposable} if $M\not\cong M_1\oplus M_2$
where $M_1$ and $M_2$ are nonzero $A$-modules.
\smallskip\noindent
An $A$-module $M$ is {\it irreducible} or {\it simple} if it has no
submodules.
}
\smallskip\noindent
In reference to modules the words ``irreducible'' and ``simple''
are used completely interchangeably.
\smallskip\noindent
{\bf 
The algebra $A$ is {\it semisimple} if 
$$\hbox{indecomposable $=$ irreducible}$$
for $A$-modules.}

\smallskip\noindent
The ``non-semisimple'' case, i.e. where indecomposable is not the same 
as irreducible, is called modular representation theory.  
We will not consider this case much in these notes.  
However, before we banish it completely let us
describe the flavor of modular representation theory.

\medskip\noindent
{\bf A {\it composition series} for $M$ 
$$M=M_0\supseteq M_1\supseteq \cdots \supseteq M_k=0$$
is a sequence of submodules of $M$ such that each $M_i/M_{i+1}$ is
simple.}

\smallskip\noindent
The Jordan-H\"older theorem says that two different composition series of
$M$ will always produce the same multiset $\{M_i/M_{i+1}\}$ of
simple modules.  Modular representation theorists are always trying
to determine this multiset of {\it composition factors of $M$}.

\bigskip\noindent
{\bf Remarks.}
\smallskip\noindent
\begin{enumerate}
\item[(1)] We shall not make life difficult in this article but 
one should note that it is common to work over general fields 
rather than just using the field $\C$.
\item[(2)] If one is bold
one can relax things in the definition of module and
let $M$ be infinite dimensional. 
\item[(3)]
Of course the definition of irreducible modules is not correct 
since $0$ and $M$ are always submodules of $M$.  So we are ignoring
these two submodules in this definition.  But conceptually the definition
is the right one, we want a simple module to be something
that has no submodules.
\item[(4)]  The definition of semisimple above is not 
a technically correct definition.  Look in Appendix A1
for the proper definition.  However, the power of semisimplicity
is exactly that it makes all indecomposable modules irreducible.
So ``indecomposable $=$ irreducible'' is really the right way to think
of semisimplicity.
\item[(5)] A good reference for the basics of representation
theory is [CR1].  The book [Bou2]
contains a completely general and comprehensive treatment of the theory 
of semisimple algebras.
The appendix, \S A1, to this article also contains a brief
(and technically correct) introduction with more specific references.
\end{enumerate}

\section*{\qquad Main questions in representation theory\qquad\null }

\bigskip\noindent
{\bf 
I.  What are the irreducible $A$-modules?
}
\medskip\noindent
What do we mean by this question?  We would like to be able
to give some kind of answers to the following more specific questions.
\smallskip\noindent
\begin{enumerate}
{\bf
     \item[(a)] How do we index/count  them?
\smallskip
     \item[(b)] What are their dimensions?  }
\smallskip
\begin{enumerate}
     \item[]  The dimension of a module is just its dimension as a vector space.
\end{enumerate}
{\bf
     \item[(c)] What are their characters?  }
\smallskip
\begin{enumerate}
     \item[]  The character of a module $M$ is the function 
$\chi_M\colon A\to \C$, 
where $\chi_M(a)$ is the trace of the linear transformation 
determined by the action
of $a$ on $M$.  More precisely,
$$\chi_M(a) = \sum_{b_i\in B} ab_i\big|_{b_i},$$
where the sum is over a basis $B$ of the module $M$ and 
$ab_i\big|_{b_i}$ denotes
the coefficient of $b_i$ in $ab_i$ when we expand in terms of the basis $B$.
\end{enumerate}
\end{enumerate}

\bigskip\noindent
{\bf 
C.  How do we construct the irreducible modules?
}

\bigskip\noindent
{\bf 
S.  Special/Interesting representations $M$}
\medskip\noindent
\begin{enumerate}
{\bf 
     \item[(a)] How does $M$ decompose into irreducibles?  }
\smallskip
\begin{enumerate}
     \item[] If we are in the semisimple case then $M$ will always 
be a direct sum of irreducible
modules.  If we group the irreducibles of the same type together we can write
$$M\cong \bigoplus_\lambda (V^\lambda)^{\oplus c_\lambda},$$
where the modules $V^\lambda$ are the irreducible $A$-modules and $c_\lambda$ is
the number of times an irreducible of type $V^\lambda$ appears as a summand in $M$.
It is common to abuse notation and write
$$M=\sum_\lambda c_\lambda V^\lambda.$$
\end{enumerate}
{\bf 
     \item[(b)] What is the character of $M$?}
\smallskip
\begin{enumerate}
     \item[]  Special modules often have particularly nice formulas
describing their characters.  It is important to note that having
a nice character formula for $M$ does not necessarily mean that it
is easy to see how $M$ decomposes into irreducibles.
Thus this question really is different from the previous one.
\end{enumerate}
\smallskip
{\bf 
     \item[(c)] How do we find interesting representations?}
\smallskip
\begin{enumerate}
     \item[] 
Sometimes special representations turn up by themselves 
and other times one has to work hard to construct the right representation
with the right properties.  Often very interesting representations come from
other fields.
\end{enumerate}
\smallskip
{\bf 
     \item[(d)] Are they useful?}
\smallskip
\begin{enumerate}
     \item[]  A representation may be particularly interesting just because 
of its structure while other times it is a special representation
that helps to prove some particularly elusive theorem.
Sometimes these representations lead to a
completely new understanding of previously known facts.  
A famous example (which unfortunately we won't have space to discuss,
see [Hu1]) is the Verma module, which was discovered in the mid 1960s and
completely changed representation theory.
\end{enumerate}
\end{enumerate}

\bigskip\noindent
{\bf 
M. The modular case}
\medskip\noindent
In the modular case we have the following important question in addition
to those above.
\smallskip
\begin{enumerate}
{\bf 
     \item[(a)] What are the indecomposable representations?
\smallskip
     \item[(b)] What are the structures of their composition series?
}
\smallskip
\end{enumerate}
For each indecomposable module $M$ there is 
a multiset of irreducibles $\{M_i/M_{i+1}\}$
determined by a composition series of $M$.
One would like to determine this multiset.

\smallskip
Even better (especially for combinatorialists),
the submodules of $M$ form a lattice under inclusion
of submodules and one would like to understand this lattice.
This lattice is always a modular lattice
and we may imagine that each edge of the Hasse
diagram is labeled by the simple module
$N_1/N_2$ where $N_1$ and $N_2 $ are the modules on the
ends of the edge. With this point of view the 
various compositions series of $M$ are 
the maximal chains in this lattice of modules.
The Jordan-H\"older
theorem says that every maximal chain in the
lattice of submodules of $M$ has the same multiset
of labels on its edges.  What modular representation
theorists try to do is determine the set of
labels on a maximal chain.

\bigskip\noindent
{\bf Remark.}
The abuse of notation which allows us to write
$M=\sum_\lambda c_\lambda V^\lambda$ has been given a formal
setting which is called the {\it Grothendieck ring}.  In other words,
the formal object which allows 
us to write such identities has been defined carefully.
See [Se1] for precise definitions of the Grothendieck ring.

\section*{\qquad Answers should be of the form \dots \qquad\null}

Now we come to the adjective ``Combinatorial.''
It refers to the way in which we give the
answers to the main questions of representation theory.

\bigskip\noindent
{\bf 
I.  What are the irreducible $A$-modules?}
\medskip\noindent
\begin{enumerate}
{\bf 
     \item[(a)] How do we index/count  them?}
\smallskip
\begin{enumerate}
     \item[] We want to answer with {\bf a {\it bijection}
between
$$\hbox{Nice combinatorial objects $\lambda$}
\qquad
\mapleftright{1-1}
\qquad
\hbox{Irreducible representations $V^\lambda$.}$$
}
\end{enumerate}

{\bf 
     \item[(b)] What are their dimensions?  }
\smallskip\noindent
\begin{enumerate}
     \item[]We should answer with {\bf a formula of 
the form
$$\dim(V^\lambda) = 
\hbox{\# of nice combinatorial objects.}
$$
}
\end{enumerate}

{\bf
     \item[(c)] What are their characters?  }
\smallskip\noindent
\begin{enumerate}
     \item[]  We want
{\bf a character formula of the type
$$\chi(a) = \sum_{T} {\bf wt}^a(T),$$
}
where the sum runs over all $T$ in a set of
nice combinatorial objects and $\wt^a$ is
a weight on these objects which depends on the
element $a\in A$ where we are evaluating the character.
\end{enumerate}
\end{enumerate}

\bigskip\noindent
{\bf C.  How do we construct the irreducible modules? }
\medskip\noindent
\begin{enumerate}
     \item[]  We want to give {\bf  constructions that
have a very explicit and very combinatorial flavor.}
What we mean by this will be more clear from the examples, see {\bf (C1-2)}
of Section 2.
\end{enumerate}

\bigskip\noindent
{\bf 
S.  Special/Interesting representations $M$}
\medskip\noindent
\begin{enumerate}
{\bf 
     \item[(a)] How does $M$ decompose into irreducibles?  }
\smallskip\noindent
\begin{enumerate}
     \item[]  If $M$ is an interesting representation
we want to {\bf determine the positive integers $c_\lambda$ }
in the decomposition
$$M\cong \bigoplus_{\lambda} (V^\lambda)^{\oplus c_\lambda}$$
{\bf in the form
$$c_\lambda = \hbox{\# of nice combinatorial objects.}$$
}
In the formula for the decomposition of $M$
the sum is over all $\lambda$ which are objects
indexing the irreducible representations of $A$.
\end{enumerate}

{\bf
     \item[(b)] What is the character of $M$? }
\smallskip\noindent
\begin{enumerate}
     \item[] As in the case I(c) we want
{\bf a character formula of the type
$$\chi_M(a) = \sum_{T} {\bf wt}^a(T)$$
}
where the sum runs over all $T$ in some set of
nice combinatorial objects and $\wt^a$ is
a weight on these objects which depends on the
element $a\in A$ where we are evaluating the character.
\end{enumerate}

\smallskip\noindent
{\bf
     \item[(c)] How do we find interesting representations? }
\smallskip\noindent
\begin{enumerate}
     \item[]
It is particularly pleasing when interesting representations
arise in {\bf  other parts of combinatorics!}  One such example is a
representation on the homology of the partition lattice which
also, miraculously, appears as a representation on the free Lie algebra.
We won't have space to discuss this here, see the original
references [Hn],[Jy], [Kl], [Sta2], the article 
[Gar] for some further basics, and 
[Ba2] for a study of how it can be that this representation appears in two 
completely different places.
\end{enumerate}

\smallskip\noindent
{\bf
     \item[(d)] Are they useful? }
\smallskip\noindent
\begin{enumerate}
     \item[] {\bf How about for solving combinatorial problems?  
Or making new combinatorial problems? }
Sometimes a representation is exactly what is most helpful for solving
a combinatorial problem.  One example of this is in the recent
solution of the the last few plane partition conjectures. 
See [Sta1], [Mac] I \S 5 Ex.\ 13-18,
for the statement of the problem and [Ku1-3] and [Ste5-7] for the solutions.
These solutions were motivated by the method of Proctor [Prc].
\end{enumerate}
\end{enumerate}

\medskip\noindent
The main point of all this is that a combinatorialist thinks in
a special way (nice objects, bijections, weighted objects, etc.)
and this method of thinking should be an integral part
of the form of the solution to the problem.

%Section 2
\section{Answers for $S_n$, the symmetric group}

Most people in the field of combinatorial representation
theory agree that the field begins with the fundamental results 
for the symmetric group $S_n$.  
Let us give the answers to the main questions for the case of $S_n$.
The precise definitions of all the objects used below can be found
in Appendix A2.  As always,
by a representation of the symmetric group we mean
a representation of its group algebra $A=\C S_n$.

\bigskip\noindent
{\bf I.  What are the irreducible $S_n$-modules?}
\medskip\noindent
\begin{enumerate}
\item[(a)] How do we index/count  them?
\smallskip
\begin{enumerate}
{\bf 
     \item[]
There is a bijection
$$
\hbox{Partitions $\lambda$ of $n$}
\qquad\mapleftright{1-1}\qquad
\hbox{Irreducible representations $S^\lambda$. }
$$
}
\end{enumerate}

\item[(b)] What are their dimensions?  
\smallskip\noindent
\begin{enumerate}
     \item[]
The dimension of the irreducible representation
$S^\lambda$ is given by
$$
\eqalign{
{\bf dim}(S^\lambda)
&= \hbox{{\bf \# of standard tableaux of shape $\lambda$}}\cr
&= {n!\over \prod_{x\in \lambda} h_x },\cr}
$$
where $h_x$ is the hook length at the box $x$ in $\lambda$,
see Appendix A2.
\end{enumerate}

\item[(c)] What are their characters?  
\smallskip\noindent
\begin{enumerate}
     \item[]  Let $\chi^\lambda(\mu)$ be the character of the irreducible
representation $S^\lambda$ evaluated at a permutation of cycle type 
$\mu=(\mu_1,\mu_2,\ldots,\mu_\ell)$.
Then {\bf the character $\chi^\lambda(\mu)$ is given by 
$$
\chi^\lambda(\mu) = \sum_{T} {\bf wt}^\mu(T), $$
where the sum is over all standard tableaux $T$ of
shape $\lambda$ and 
$$
{\bf wt}^\mu(T) = \prod_{i=1}^n f(i,T),$$
where
$$
f(i,T) = \cases{
-1, &if $i\not\in B(\mu)$ and $i+1$ is sw of $i$,\cr
0, &if $i,i+1\not\in B(\mu)$, 
$i+1$ is ne of $i$, and 
$i+2$ is sw of $i+1$,\cr
1, &otherwise,\cr}
$$
and
$B(\mu)= \{ \mu_1+\mu_2+\cdots+\mu_k | 1\le k\le \ell \}
$.
}
In the formula for $f(i,T)$, {\bf sw} means strictly south and weakly west
and {\bf ne} means strictly north and weakly east.
\end{enumerate}
\end{enumerate}

\bigskip\noindent
{\bf C.  How do we construct the irreducible modules? }
\smallskip\noindent
There are several interesting constructions
of the irreducible $S^\lambda$.

\medskip\noindent
\begin{enumerate}
\item[(i)] via {\bf Young symmetrizers}.
\smallskip\noindent
Let $T$ be a tableau.  Let
$$\eqalign{
R(T) &= \hbox{ permutations which fix the rows of $T$, as sets;} \cr
C(T) &= \hbox{ permutations which fix the columns of $T$, as sets;} \cr
}$$
$$
P(T) = \sum_{w\in R(T)} w,\quad\hbox{and}\quad
N(T) = \sum_{w\in C(T)} \varepsilon(w)w,
$$
where $\varepsilon(w)$ is the sign of the permutation $w$. Then
$$S^\lambda \cong \C S_n P(T)N(T),$$
where the action of the symmetric group is by left multiplication.

\item[(ii)] {\bf Young's seminormal construction}.

\smallskip\noindent
Let 
{\bf 
$$S^\lambda = 
\hbox{\C -span-}
\{ v_T \ |\ \hbox{$T$ are standard tableaux of shape $\lambda$} \}
$$
}
so that the vectors $v_T$ are a basis of $S^\lambda$.
{\bf The action of $S_n$ on $S^\lambda$ is given by
$$
s_i v_T = (s_i)_{TT}v_T + (1+(s_i)_{TT})v_{s_iT},
\qquad\hbox{where}\qquad
(s_i)_{TT} = { 1\over c(T(i+1))-c(T(i)) }.$$
and $s_i = (i,i+1)$. }
In this formula
\smallskip
\begin{enumerate}
{\bf 
     \item[]
$T(i)$ denotes the box containing $i$ in $T$; 
     \item[]
$c(b)=j-i$ is the {\it content} of the box $b$, where $(i,j)$ is 
the position of $b$ in $\lambda$;
     \item[]
$s_iT$ is the same as $T$ except that the entries
$i$ and $i+1$ are switched;
     \item[]
$v_{s_iT} = 0$ if $s_iT$ is not a standard tableau.
}
\end{enumerate}

\medskip\noindent
There are other important constructions of the irreducible
representations $S^\lambda$.  We do not have room to discuss 
these constructions
here, see the Notes and References (10-12) below and A3 
in the appendix.
The main ones are:

     \item[(iii)]  {\bf Young's orthonormal construction},

     \item[(iv)] {\bf The Kazhdan-Lusztig construction},

     \item[(v)] {\bf The Springer construction}.
\end{enumerate}

\bigskip\noindent
{\bf S.  Particularly interesting representations}

\medskip\noindent
\begin{enumerate}
\item[\bf (S1)] Let $k+\ell = n$.
The module
$S^\lambda \big\downarrow^{S_n}_{S_k\times S_\ell}
$
is the same as $S^\lambda$ except that we only look at the
action of the subgroup $S_k\times S_\ell$.
Then
$$S^\lambda \big\downarrow^{S_n}_{S_k\times S_\ell}
= \bigoplus_{\mu\vdash k, \nu\vdash \ell}
(S^\mu\otimes S^\nu)^{\oplus c^\lambda_{\mu\nu} }
= \sum_{\mu,\nu} c_{\mu\nu}^\lambda(S^\mu\otimes S^\nu),$$
where {\bf $c_{\mu\nu}^\lambda$ is the number of 
column strict fillings of $\lambda/\mu$ of content $\nu$
such that the word of the filling is a lattice permutation.}
The positive integers
$c^\lambda_{\mu\nu}$ are the {\it Littlewood-Richardson
coefficients}.
See Appendix A2.

\item[\bf (S2)] Let $\mu=(\mu_1,\ldots,\mu_\ell)$ be a partition of $n$.
Let $S_\mu= S_{\mu_1}\times\cdots\times S_{\mu_\ell}$.  
%Let $\One$ be the trivial representation of $S_\mu$, i.e.
%the one dimensional module $\One = \C v$ where
%$wv = v$ for $w\in S_\mu$.  
The module $\One \big\uparrow_{S_\mu}^{S_n}$ is 
the vector space 
$$\One \big\uparrow_{S_\mu}^{S_n} = \C(S_n/S_\mu) = 
\hbox{$\C$-span}\{ wS_\mu\ |\ w\in S_n\}$$
where the action of $S_n$ on the cosets is by left multiplication.
Then
$$
\One \big\uparrow^{S_n}_{S_\mu}
= \sum_\lambda K_{\lambda\mu} S^\lambda,$$
where 
{\bf
$$K_{\lambda\mu} = \hbox{\# of column strict
tableaux of shape $\lambda$ and weight $\mu$. }$$
}
This representation also occurs in 
the following context: 
$$
\One \big\uparrow^{S_n}_{S_\mu}
\cong
H^*({\cal B}_u),$$
where $u$ is a unipotent element of $GL(n,\C)$
with Jordan decomposition $\mu$ and ${\cal B}_u$ is the variety 
of Borel sugroups in $GL(n,\C)$ containing $u$. This representation
is related to the Springer construction mentioned in {\bf C(v)} above.
See Appendix A3 for further details.

\item[\bf (S3)]
 If $\mu,\nu\vdash n$ then the tensor product 
$S_n$-module $S^\mu\otimes S^\nu$ is 
defined by 
$w(m\otimes n) = wm\otimes wn,$
for all $w\in S_n$, $m\in S^\mu$ and $n\in S^\nu$.
{\bf There are positive integers $\gamma_{mu\nu\lambda}$ such that
$$
S^\mu\otimes S^\nu = \sum_{\lambda\vdash n}
\gamma_{\mu\nu\lambda} S^\lambda.$$
}
Except for a few special cases the positive integers
$\gamma_{\mu\nu\lambda}$ 
are still unknown. See [Rem] for a combinatorial
description of the cases for which the coefficients 
$\gamma_{\mu\nu\lambda}$ are known.
\end{enumerate}

\bigskip\noindent
{\bf Notes and references}
\smallskip\noindent
\begin{enumerate}
\item[(1)]  The bijection in {\bf (Ia)}, between irreducible representations
and partitions, is due to Frobenius [Fr].  Frobenius is 
the founder of representation theory and the symmetric group
was one of the first examples that he worked out.
\item[(2)]  The formula in {\bf (Ib)} for the dimension of $S^\lambda$
as the number of standard tableaux is immediate from
the work of Frobenius, but it really came into the
fore from the work of Young [Y].
The ``hook formula'' for $\dim(S^\lambda)$ is
due to Frame-Robinson-Thrall [FRT].
\item[(3)] The formula for the characters of the symmetric group
which is given in {\bf (Ic)} is due to
Fomin and Greene [FG].  For them, this formula
arose by application of their theory of noncommutative
symmetric functions.
Roichman [Ro] discovered this formula independently
in the  more general case of the Iwahori-Hecke algebra.
The formula for the Iwahori-Hecke algebra is exactly the same as the
formula for the $S_n$ case except that 
the $1$ appearing in case 3 of the definition of
$f(i,T)$ should be changed to a $q$.
\item[(4)]
There is a different and more classical formula
for the characters than the formula given in {\bf (Ic)}
which is called the Murnaghan-Nakayama rule [Mur] [Nak].
We have described the Murnaghan-Nakayama rule 
in the appendix, Theorem A2.2.
Once the formula in {\bf (Ic)} is given it is not hard to show 
combinatorially that
it is equivalent to the Murnaghan-Nakayama rule
but if one does not know the formula
it is nontrivial to guess it from the Murnaghan-Nakayama
rule.
\item[(5)]  We do not know if anyone has compared the
algorithmic complexity of the formula given in {\bf (Ic)} with the 
algorithmic complexity of the Murnaghan-Nakayama rule.
One would expect that they have the same complexity:
the formula above is a sum over more objects
than the sum in the Murnaghan-Nakayama rule but
these objects are easier to create and many of them have
zero weight.  
\item[(6)]  One of the beautiful things about the formula for the 
character of $S^\lambda$ which is given in {\bf (Ic)} is that it is a sum
over the same set that we have used to describe the dimension
of $S^\lambda$.
\item[(7)]  The construction of $S^\lambda$ by Young symmetrizers
is due to Young [Y1-2] from 1900.  
It is used so often and has so many applications
that it is considered classical.  A review and generalization of this
construction to skew shapes appears in [GW].
\item[(8)]  The seminormal form construction of $S^\lambda$ 
is also due to Young [Y4-5] 
although it was discovered some thirty years after the 
Young symmetrizer construction. 
\item[(9)]
Young's orthonormal construction differs
from the seminormal construction only by
multiplication of the basis vectors by certain
constants.  A comprehensive treatment of all three constructions of Young
is given in the book by Rutherford [Ru].
\item[(10)]  The Kazhdan-Lusztig construction uses the
Iwahori-Hecke algebra in a crucial way.  It is combinatorial
but relies crucially on certain polynomials which seem to
be impossible to compute in practice except for very small $n$,
see [Bre] for further information.  
This construction has important
connections to geometry and other parts of representation
theory.
The paper [GM] and the book [Hu2] give elementary treatments of the
Kazhdan-Lusztig construction.
\item[(11)]
Springer's construction is a geometric construction.  In this
construction the irreducible module $S^\lambda$ is realized
as the top cohomology group of a certain variety, see [Spr], [CG],
and Appendix A3.
\item[(12)]  There are many ways of constructing new representations
from old ones.  Among the common techniques are 
{\it restriction}, {\it induction}, and {\it tensoring}.  
The special representations
{\bf (S1)}, {\bf (S2)}, and {\bf (S3)} given above are particularly nice 
examples of these constructions.  One should note that tensoring
of representations works for group algebras (and Hopf algebras)
but not for general algebras.
\end{enumerate}

%Section 3
\section{Answers for $GL(n,\C)$, the general linear group}

The results for the general linear group are just as beautiful
and just as fundamental as those
for the symmetric group.  The results are surprisingly similar and yet
different in many crucial ways.  We shall see that the results for 
$GL(n,\C)$ have been generalized to a very wide class of groups
whereas the results for $S_n$ have only been generalized successfully
to groups that look very similar to symmetric groups.
The representation theory of $GL(n,\C)$ was put on a very firm footing
from the fundamental work of Schur [Sc1-2] in 1901 and 1927.

\bigskip\noindent
{\bf I. What are the irreducible $GL(n,\C)$-modules?}
\medskip\noindent
\begin{enumerate}
\item[(a)]  How do we index/count them?
\smallskip\noindent
\begin{enumerate}
\item[]  {\bf There is a bijection
$$
\hbox{Partitions $\lambda$ with at most $n$ rows}
\quad
\mapleftright{1-1}
\quad
\hbox{Irreducible polynomial representations $V^\lambda$.}
$$
}
See Appendix A4 for a definition and discussion of what it
means to be a {\it polynomial} representation.
\end{enumerate}

\item[(b)] What are their dimensions?
\smallskip\noindent
\begin{enumerate}
\item[]  The {\bf dimension of the irreducible representation $V^\lambda$
is given by
$$\eqalign{
{\bf dim}(V^\lambda) &= \hbox{\# of column strict tableaux
of shape $\lambda$} \cr
&\phantom{= \qquad} \hbox{filled with entries from $\{1,2,\ldots,n\}$}\cr
&= \prod_{x\in \lambda} {n+c(x)\over h_x},\cr}$$
}
where $c(x)$ is the content of the box $x$ and
$h_x$ is the hook length at the box $x$.
\end{enumerate}

\item[(c)] What are their characters?
\smallskip\noindent
\begin{enumerate}
\item[] Let $\chi^\lambda(g)$ be the character of the 
irreducible representation $V^\lambda$ evaluated at an element
$g\in GL(n,\C)$.  The {\bf character $\chi^\lambda(g)$ is
given by
$$
\eqalign{
\chi^\lambda(g) &= \sum_T x^T \cr
&= {\sum_{w\in S_n} \varepsilon(w)wx^{\lambda+\delta}
\over \sum_{w\in S_n} \varepsilon(w) wx^\delta}
= {\det(x_i^{\lambda_j+n-j}) \over \det(x_i^{n-j}) }, \cr
}$$
where the sum is over all column strict tableaux $T$ of
shape $\lambda$ filled with entries from $\{1,2,\ldots,n\}$
and 
$$x^T = x^{\mu_1}x_2^{\mu_2}\cdots x_n^{\mu_n},
\quad \hbox{ where $\mu_i=$ \# of $i$'s in $T$}
$$
and $x_1,x_2,\ldots, x_n$ are the eigenvalues of the matrix $g$.
}
Let us not worry about the first expression in the second
line at the moment. Let us only say that it is routine
to rewrite it as the second expression in that 
line which is one of the standard expressions for the 
Schur function, see [Mac] I \S 3.
\end{enumerate}
\end{enumerate}

\bigskip\noindent
{\bf 
C. How do we construct the irreducible modules?}
\medskip\noindent
There are several interesting constructions of the irreducible
$V^\lambda$.

\begin{enumerate}
\item[\bf (C1)]  
via {\bf Young symmetrizers}.
\smallskip\noindent
\begin{enumerate}
\item[] Recall that the irreducible $S^\lambda$ of the
symmetric group $S_k$ was constructed via Young symmetrizers
in the form 
$$S^\lambda\cong \C S_nP(T)N(T).$$
We can construct the irreducible $GL(n,\C)$-module
in a similar form.  
{\bf If $\lambda$ is a partition of $k$ then 
$$V^\lambda \cong V^{\otimes k}P(T)N(T).$$
}
This important construction is detailed in Appendix A5.
\end{enumerate}

\item[\bf (C2)]  {\bf Gelfan'd-Tsetlin bases}
\smallskip\noindent
\begin{enumerate}
\item[] This construction of the irreducible $GL(n,\C)$ representations
$V^\lambda$ is analogous to 
the Young's seminormal construction of the irreducible representations
$S^\lambda$ of the symmetric group.  Let
{\bf\small
$$V^\lambda =
\hbox{span-}
\{ v_T \ |\ \hbox{$T$ are column strict tableaux of shape $\lambda$
filled with elements of $\{1,2,\ldots,n\}$ } 
\}
$$
}
so that the vectors $v_T$ are a basis of $V^\lambda$.
Define an action of symbols $E_{k-1,k}$, $2\le k\le n$,
on the basis vectors $v_T$ by
$$
E_{k-1,k}v_T = \sum_{T^-} a_{T^-T}(k)v_{T^-},
$$
where the sum is over all column strict tableaux
$T^-$ which are obtained from $T$ by
changing a $k$ to a $k-1$ and the coefficients $a_{T^-T}(k)$ are given by
$$a_{T^-T}(k)
= - {\displaystyle{\prod_{i=1}^k (T_{ik}-T_{j,k-1}+j-k)}
\over
\displaystyle{\prod_{i=1\atop i\ne j}^{k-1} (T_{i,k-1}-T_{j,k-1}+j-k)} },
$$
where $j$ is the row number of the entry where $T^-$ and $T$
differ and $T_{ik}$ is the position of the rightmost entry $\le k$
in row $i$ of $T$.
Similarly, define an action of symbols $E_{k,k-1}$, $2\le k\le n$,
on the basis vectors $v_T$ by
$$
E_{k,k-1}v_T = \sum_{T^+} b_{T^+T}(k)v_{T^+},
$$
where the sum is over all column strict tableaux
$T^+$ which are obtained from $T$ by
changing a $k-1$ to a $k$ and the coefficients $b_{T^+T}(k)$ are given by
$$b_{T^+T}(k)
=  {\displaystyle{\prod_{i=1}^{k-2} (T_{i,k-2}-T_{j,k-1}+j-k)}
\over
\displaystyle{\prod_{i=1\atop i\ne j}^{k-1} (T_{i,k-1}-T_{j,k-1}+j-k) } },
$$
where $j$ is the row number of the entry where $T^+$ and $T$
differ and $T_{ik}$ is the position of the rightmost entry $\le k$
in row $i$.

Since 
$$\matrix{
&g_i(z) = 
\pmatrix{ 1      &0      &  &\cdots  &   &         &0       \cr
          0      &\ddots &  &        &   &         &        \cr
                 &       &1 &        &   &         &\vdots  \cr
          \vdots &       &  &z       &   &         &        \cr
                 &       &  &        &1  &         &        \cr
                 &       &  &        &   &\ddots   &0       \cr
          0      &       &  &\cdots  &   &0        &1       \cr},
\quad &z\in \C^*,\cr
\cr
\cr
&g_{i-1,i}(z) = 
\pmatrix{ 1      &0      &\cdots  &  &         &0       \cr
          0      &\ddots &        &  &         &        \cr
                 &       &1       &z &         &\vdots  \cr
          \vdots &       &0       &1 &         &        \cr
                 &       &        &  &\ddots   &0       \cr
          0      &       &\cdots  &  &0        &1       \cr},
\quad &z\in \C,\cr
\cr
\cr
&g_{i,i-1}(z) = 
\pmatrix{ 1      &0      &\cdots  &  &         &0       \cr
          0      &\ddots &        &  &         &        \cr
                 &       &1       &  &         &\vdots  \cr
          \vdots &       &z       &1 &         &        \cr
                 &       &        &  &\ddots   &0       \cr
          0      &       &\cdots  &  &0        &1       \cr},
\quad &z\in \C,\cr
}
$$
generate $GL(n,\C)$, the action of
these matrices on the basis vectors $v_T$ will determine the
action of all of $GL(n,\C)$ on the space $V^\lambda$.  
The action of these generators is given by:
$$\eqalign{
g_i(z) v_T &= z^{\hbox{(\# of $i$'s in $T$)}} v_T,\cr
g_{i-1,i}(z) v_T &= e^{zE_{i-1,i}} v_T
= (1 + zE_{i-1,i} + {1\over 2!}z^2E_{i-1,i}^2 + \ldots) v_T,\cr
g_{i,i-1}(z) v_T &= e^{zE_{i,i-1}} v_T
= (1 + zE_{i,i-1} + {1\over 2!}z^2E_{i,i-1}^2 + \ldots) v_T.\cr
}
$$
\end{enumerate}

%$$a_{ST} = 
%{\prod_{i=1}^k  p(\le k,i)-i+r(b(S,T))-p(b(S,T))+1
%\over 
%\prod_{i=1\atop i\ne r(b(S,T))}^{k-1} 
%p(\le k-1,i)-i+r(b(S,T))-p(b(S,T))+1}.$$
%where 
%\smallskip
%\begin{enumerate}
     \item[] $b(S,T)$ is the box where $S$ and $T$ differ,
%\smallskip
%\begin{enumerate}
     \item[] $r(b(S,T))$ is the row number of the box $b(S,T)$,
%\smallskip
%\begin{enumerate}
     \item[] $p(b(S,T))$ is the position of the box $b(S,T)$ in its row,
%and
%\smallskip
%\begin{enumerate}
     \item[] $p(\le k,i)$ is the position of the rightmost entry .
%$\le k$ in row $i$ of $T$.

\item[\bf (C3)]  {\bf The Borel-Weil-Bott construction}
\smallskip\noindent
\begin{enumerate}
\item[]
Let $\lambda$ be a partition.  Then $\lambda$ defines a character 
(one-dimensional representation) of the group $T_n$ of diagonal
matrices in $G=GL(n,\C)$.  This character can be extended to
the group $B=B_n$ of upper triangular matrices in $G=GL(n,\C)$ by
letting it act trivially on $U_n$ the group of upper unitriangular
matrices in $G=GL(n,\C)$.  Then the fiber product
$${\cal L}_\lambda = G\times_B \lambda$$
is a line bundle on $G/B$.  Finally,
$$V^\lambda \cong H^0(G/B, {\cal L}_\lambda),$$
{\bf where $H^0(G/B, {\cal L}_\lambda)$ is the space of global sections
of the line bundle ${\cal L}_\lambda$.}
More details on the construction of the character $\lambda$
and the line bundle ${\cal L}_\lambda$ are given in Appendix A6.
\end{enumerate}
\end{enumerate}

\bigskip\noindent
{\bf S. Special/Interesting representations}

\medskip\noindent
\begin{enumerate}
\item[\bf (S1)] Let
$$GL(k)\times GL(\ell) = 
\let\quad=\enspace
\pmatrix{
\pmatrix{ &&\cr &GL(k,\C) &\cr && \cr} 
&\matrix{ &&\cr &0 &\cr && \cr} \cr 
\matrix{ &&\cr &0 &\cr && \cr} 
&\pmatrix{ &&\cr &GL(\ell,\C) &\cr && \cr} 
\cr }
\subseteq GL(n),
\qquad\hbox{where $k+\ell=n$.}$$
Then
$$V^\lambda\big\downarrow^{GL(n)}_{GL(k)\times GL(\ell)}
%= \bigoplus_{\mu,\nu} (V^\mu\otimes V^\nu)^{\oplus c_{\mu\nu}^\lambda}
=\sum_{\mu,\nu} c_{\mu\nu}^{\lambda}(V^\mu\otimes V^\nu),$$
where {\bf $c_{\mu\nu}^\lambda$ is the number of column strict
fillings of $\lambda/\mu$ with content $\nu$ such that
the word of the filling is a lattice permutation.}
The positive integers
$c^\lambda_{\mu\nu}$ are the {\it Littlewood-Richardson
coefficients} that appeared earlier in the decomposition
of $S^\lambda \big\downarrow^{S_n}_{S_k\times S_\ell}$
in terms of $S^\mu\otimes S^\nu$.

\begin{enumerate}
\item[]
We may write this expansion in the form
$$V^\lambda\big\downarrow^{GL(n)}_{GL(k)\times GL(\ell)}
=\sum_{F\hbox{ fillings}} V^{\mu(F)}\otimes V^{\nu(F)}.$$
\begin{enumerate}
\item[] We could do this precisely if we wanted. We won't do it now,
but the point is that it may be nice to write this
expansion as a sum over combinatorial objects.  This
will be the form in which this will be generalized later.
\end{enumerate}

\item[\bf (S2)] Let $V^\mu$ and $V^\nu$ be irreducible 
polynomial representations of $GL(n)$.  Then
$$V^\mu\otimes V^\nu = \sum_\lambda c_{\mu\nu}^\lambda V^\lambda,$$
where $GL(n)$ acts on $V^\mu\otimes V^\nu$ by 
$g(m\otimes n) = gm\otimes gn$,
for $g\in GL(n,\C)$, $m\in V^\mu$ and $n\in V^\nu$.
Amazingly, {\bf the coefficients $c_{\mu\nu}^\lambda$ are the 
Littlewood-Richardson coefficients again.}  These are the same coefficients
that appeared in the {\bf (S1)} case above and in the {\bf (S1)} case for
the symmetric group.
\end{enumerate}

\bigskip\noindent
{\bf Remarks}
\smallskip\noindent
\begin{enumerate}
\item[(1)] There is a strong similarity between the results
for the symmetric group and the results for
$GL(n,\C)$.  One might wonder whether there is any 
connection between these two pictures.

\smallskip\noindent
{\bf 
There are TWO DISTINCT ways of making concrete connections
between the representation theories of $GL(n,\C)$ and the
symmetric group.}  In fact these two are so different that
{\bf DIFFERENT SYMMETRIC GROUPS} are involved.
\smallskip\noindent
\begin{enumerate}
     \item[(a)]  {\it If $\lambda$ is a partition of $n$}
then the ``zero weight space'', or $(1,1,\ldots,1)$ weight space,
of the irreducible $GL(n,\C)$-module $V^\lambda$ is
isomorphic to the irreducible module $S^\lambda$
for the group $S_n$, where the
action of $S_n$ is determined by the fact that 
$S_n$ is the subgroup of permutation matrices in $GL(n,\C)$.
This relationship is reflected in the combinatorics:
the standard tableaux of shape $\lambda$ are exactly the column
strict tableaux of shape $\lambda$ which are of weight 
$\nu=(1,1,\ldots,1)$.
     \item[(b)] Schur-Weyl duality, see \S A5 in the appendix,
says that the action of the symmetric group $S_k$ on $V^{\otimes k}$ 
by permutation of the tensor factors generates the full
centralizer of the $GL(n,\C)$-action on $V^{\otimes k}$
where $V$ is the standard $n$-dimensional representation
of $GL(n,\C)$.  
By double centralizer theory, this duality induces a 
correspondence between the irreducible representations
of $GL(n,\C)$ which appear in $V^{\otimes k}$ and
the irreducible representations of $S_k$ which appear
in $V^{\otimes k}$.  These representations are indexed by 
{\it partitions $\lambda$ of $k$}.
\smallskip\noindent
\end{enumerate}
\item[(2)]  It is important to note that the word {\it character}
has two different and commonly used meanings and the use of
the word character in {\bf (C3)} is different than 
in Section 1.  In {\bf (C3)} above the word
character means  {\it one dimensional representation}.  This terminology
is used particularly (but not exclusively) in reference 
to representations of abelian groups (like the group $T_n$ 
in {\bf (C3)}).  In general
one has to infer from the context which meaning is intended.
\smallskip\noindent
\end{enumerate}
\item[(3)]  The indexing and the formula for the characters of 
the irreducible representations is due to Schur [Sc1].
\item[(4)]  The formula for the dimensions of the irreducibles as the
number of column strict tableaux follows from the work of 
Kostka [Kk] and Schur [Sc1].  The ``hook-content'' formula 
appears in [Mac] I \S 3 Ex. 4, where the book of Littlewood [Lw]
is quoted.
\item[(5)]  The construction of the irreducibles by Young symmetrizers
appeared in 1939 in the influential book [Wy1] of H. Weyl.
It was generalized to the symplectic and orthogonal groups
by H. Weyl in the same book.  Further important information
about this construction in the symplectic and orthogonal
cases is found in [Be2] and  [KWe].  It is not known how to
generalize this construction to arbitrary complex semisimple Lie groups.
\item[(6)] 
The Gelfan'd-Tsetlin basis construction originates from 1950 [GT1].
A similar construction was given for the orthogonal group at the
same time [GT2] and was generalized to the 
symplectic group by Zhelobenko, see [Zh1-2].  
This construction does not generalize
well to other complex semisimple groups since it depends crucially
on a tower $G \supseteq G_1\supseteq \cdots \supseteq G_k \supseteq \{1\}$
of ``nice'' Lie groups such that all the combinatorics is controllable.
\item[(7)]  The Borel-Weil-Bott construction is not a combinatorial
construction of the irreducible module $V^\lambda$.  
It is very important because it is a construction that generalizes well
to all other compact connected real Lie groups.
\item[(8)]  The facts about the special representations which we have given
above are found in Littlewood's book [Lw].
\end{enumerate}

%Section 4
\section{Answers for finite dimensional complex 
semisimple Lie algebras ${\goth g}$}

Although the foundations for generalizing the $GL(n,\C)$ results
to all complex semisimple Lie groups and Lie algebras were laid in 
the fundamental work of Weyl [Wy2] in 1925, it is only recently
that a complete generalization of the tableaux results for 
$GL(n,\C)$ has been obtained by Littelmann [Li2].  
The results which we state below are generalizations of 
those given for $GL(n,\C)$ in the last section; partitions get
replaced by points in a lattice called $P^+$, and column 
strict tableaux get replaced by paths.  See the Appendix A7 for
some basics on complex semisimple Lie algebras.

\bigskip\noindent
{\bf I. What are the irreducible ${\goth g}$-modules?}
\medskip\noindent
\begin{enumerate}
\item[(a)]  How do we index/count them?
\smallskip\noindent
\begin{enumerate}
\item[]  {\bf There is a bijection
$$\lambda\in P^+
\qquad
\mapleftright{1-1}
\qquad
\hbox{irreducible representations $V^\lambda$},
$$
where $P^+$ is the cone of dominant integral weights for ${\goth g}$.}
The set $P^+$ is described in Appendix A8.
\end{enumerate}
\item[\bf (Ib)] What are their dimensions?
\smallskip\noindent
\begin{enumerate}
\item[]  {\bf The dimension of the irreducible representation $V^\lambda$
is given by
$$\eqalign{
{\bf dim}(V^\lambda) &=
\hbox{\# of paths in ${\cal P}\pi_\lambda$ } \cr
&=
\prod_{\alpha>0} {
\langle \lambda+\rho,\alpha\rangle
\over
\langle \rho,\alpha\rangle }, \cr
}$$
where
\smallskip
\begin{enumerate}
     \item[] {\bf ${\displaystyle{
\rho=\hbox{$1\over 2$}\sum_{\alpha>0} \alpha}}$,\enspace 
is the half sum of the positive roots,}
     \item[] {\bf $\pi_\lambda$ is the straight line path from 
$0$ to $\lambda$, and}
     \item[]  {\bf ${\cal P}\pi_\lambda = 
\{ f_{i_1}\cdots f_{i_k}\pi_\lambda | 1\le i_1,\ldots, i_k\le n \}$, where }
     \item[] {\bf $f_1,\ldots,f_n$ are the path operators introduced in [Li2].}
\end{enumerate}
}
\item[]
We shall not define the operators $f_i$ here (or in the appendix, see 
[Li2]), let us just say that 
they act on paths and they are partial permutations in the sense
that if $f_i$ acts on a path $\pi$ then the result is either
$0$ or another path.  See Appendix A8 for a few more details.

\item[(c)] What are their characters?
\smallskip\noindent
\begin{enumerate}
\item[] {\bf   The character of the irreducible module $V^\lambda$ is 
given by
$$\eqalign
{\chr(V^\lambda) 
&= \sum_{\eta\in {\cal P}\pi_\lambda} e^{\eta(1)} \cr
&= 
{\sum_{w\in W} \varepsilon(w) e^{w(\lambda+\rho)}
\over
\sum_{w\in W} \varepsilon(w) e^{w\rho} }, \cr
}$$
where $\eta(1)$ is the endpoint of the path $\eta$.}
These expressions live in the group algebra of the weight lattice $P$,
$\C[P] = \span\{ e^\mu \ |\ \mu\in P\},$
where $e^\mu$ is a formal variable indexed by $\mu$ and
the multiplication is given by
$e^\mu e^\nu = e^{\mu+\nu},$ for $\mu,\nu\in P$.
See Appendix A7 for more details.
\end{enumerate}
\end{enumerate}

\bigskip\noindent
{\bf S. Special/Interesting representations}

\medskip\noindent
{\bf (S1)}  Let ${\goth l}\subseteq {\goth g}$ be a Levi
subalgebra of ${\goth g}$ (this is a Lie algebra corresponding
to a subgraph of the Dynkin diagram which corresponds to 
${\goth g}$).  The subalgebra ${\goth l}$ corresponds to a subset
$J$
of the set $\{\alpha_1,\ldots,\alpha_n\}$ of simple roots.
The {\bf restriction rule from ${\goth g}$ to ${\goth l}$ is
$$V^\lambda\big\downarrow^{\goth g}_{\goth l}
= \sum_\eta V^{\eta(1)}, $$
where 
\begin{enumerate}
     \item[] the sum is over all paths $\eta\in {\cal P}\pi_\lambda$
such that $\eta\in \overline{C}_{\goth l}$, 
     \item[] $\eta\in \overline{C}_{\goth l}$ means that 
$\langle \eta(t),\alpha_i\rangle \ge 0$,
for all $t\in [0,1]$ and all $\alpha_i\in J$. 
\end{enumerate}
}
\item[\bf (S2)]  {\bf The tensor product of two irreducible 
modules is given by
$$V^\mu\otimes V^\nu
= \sum_{\eta}
V^{\mu+\eta(1)}, $$
where the sum is over all paths $\eta\in {\cal P}\pi_\nu$ such that
$\pi_\mu*\eta\in \overline{C}$,  
\smallskip
\begin{enumerate}
     \item[] $\pi_\mu$ and $\pi_\nu$ are straight line paths 
from $0$ to $\mu$ and $0$ to $\nu$, respectively,
     \item[] ${\cal P}\pi_\nu$ is as in (Ib),
     \item[] $\pi_\mu*\eta$ is the path obtained by attaching $\eta$ to
the end of $\pi_\mu$, and 
     \item[] $(\pi_\mu*\eta)\in \overline{C}$ means that
$\langle (\pi_\mu*\eta)(t),\alpha_i\rangle \ge 0$,
for all $t\in [0,1]$ and all simple roots $\alpha_i$.
\end{enumerate}
}
\end{enumerate}

\bigskip\noindent
{\bf Notes and references}
\smallskip\noindent
\begin{enumerate}
\item[(1)]  The indexing of irreducible representations
given in {\bf (Ia)} is due to Cartan and Killing,
the founders of the theory, from around the turn of the century.  
Introductory treatments of this
result can be found in [FH] and [Hu1].
\item[(2)] The first equality in {\bf (Ib)} is due to Littelmann [Li1], but
his later article [Li2] has some improvements and
can be read independently, so we recommend the later article.
This formula for the dimension of the irreducible representation,
the number of paths in a certain set,
is exactly analogous to the formula in the $GL(n,\C)$ case,
the number of tableaux which satisfy a certain condition.
The second equality is the Weyl dimension formula which
was originally proved in [Wy2].  It can be proved easily
from the Weyl character formula given in {\bf (Ic)}, see
[Hu1] and [Ste4] Lemma 2.5.
This product formula is an analogue of the ``hook-content'' formula
given in the $GL(n,\C)$ case.
\item[(3)]  A priori, it might be possible that the
set ${\cal P}\pi_\lambda$ is an infinite set, at least the
way that we have defined it. In fact, this set is always finite
and there is a description of the paths that are contained in it.
The paths in this set are called {\it Lakshmibai-Seshadri paths},
see [Li2].
The explicit description of these paths is a generalization of the
types of indexings that were used in the
``standard monomial theory'' of Lakshmibai and Seshadri [LS].
\item[(4)] The first equality in {\bf (Ic)} is due to Littelmann [Li2].
This formula, a weighted sum over paths, is an analogue of the formula
for the irreducible character of $GL(n,\C)$ as a weighted sum of
column strict tableaux. 
The second equality in {\bf (Ic)} is the celebrated Weyl character formula which
was originally proved in [Wy2].  A modern treatment of this
formula can be found in [BtD], [Hu1], and [Va].
\item[(5)] The general restriction formula {\bf (S1)} is due to Littelmann [Li2].
This is an analogue of the rule given in {\bf (S1)} of the $GL(n,\C)$ results.
In this case the formula is as a sum over paths which satisfy certain
conditions whereas in the $GL(n,\C)$ case the formula is a
sum over column strict fillings which satisfy a certain condition.
\item[(6)] The general tensor product formula in {\bf (S2)}
is due to Littelmann [Li2].  This formula is an analogue of the
formula given in {\bf (S2)} of the $GL(n,\C)$ results.
\item[(7)]  The results of Littelmann given above are some of the most
exciting results of combinatorial representation theory in recent years.  
They were very much inspired by some very explicit conjectures
of Lakshmibai, see [LS], which arose out of the ``standard
monomial theory'' developed by Lakshmibai and Seshadri.
Although Littelmann's theory is actually much more general than we have
stated above, the special set of paths ${\cal P}\pi_\lambda$ used in
{\bf (Ib-c)} is a modified description of the same set which appeared in 
Lakshmibai's conjecture.
Another important influence on Littelmann in his work was
Kashiwara's work on crystal bases [Ksh].  
\end{enumerate}

\part*{Part II}

%Section 5
\section{Generalizing the $S_n$ results}

Having the above results for the symmetric group in hand we
would like to try to 
generalize as many of the $S_n$ results to other 
similar groups and algebras as we can.
Work along this line began almost immediately
after the discovery of the $S_n$ results and it continues today.
In the current state of results this has been largely
\smallskip\noindent
\begin{enumerate}
\item[\bf (1)]  {\bf successful} for the complex reflection
groups $G(r,p,n)$ and their ``Hecke algebras,''
\item[\bf (2)]  {\bf successful} for tensor power centralizer
algebras and their $q$-analogues, and
\item[\bf (3)] {\bf unsuccessful} for general
Weyl groups and finite Coxeter groups.
\end{enumerate}

\medskip\noindent
Let us give a brief description of what the objects are in {\bf (1)}, 
{\bf (2)}, and {\bf (3)}.
For more precise definitions and discussion of everything below
see Appendix B.

\section*{\qquad Definitions\qquad\null }

\medskip\noindent
\leftline{\bf (1) Complex reflection groups $G(r,p,n)$ and their
Hecke algebras}
\medskip\noindent
\leftline{\it The complex reflection groups $G(r,p,n)$}
\medskip
{\bf A finite Coxeter group is a finite group which is 
generated by reflections in $\R^n$.}  
In other words, take a bunch of
linear transformations of $\R^n$ which are reflections
(in the sense of reflections and rotations in
the orthogonal group)  and see what group they generate.  
If the group is finite then it is a finite Coxeter group.
Actually, this definition of finite Coxeter group is
not the usual one (for that see Appendix B1), but 
since we have the following theorem we are not too far astray.

\thm A group is a finite group generated by reflections if and 
only if it is a finite Coxeter group.
\endthm

\noindent
The finite Coxeter groups have been classified completely
and there is one group of each of the following ``types''
$$A_n,\enspace B_n,\enspace D_n,\enspace E_6,\enspace E_7, \enspace E_8,
\enspace F_4, \enspace H_3, \enspace H_4, \enspace\hbox{or}\enspace I_2(m).$$
The finite {\it crystallographic} reflection groups are called {\it Weyl
groups} because of their connection with Lie theory.  These are the finite
Coxeter groups of types
$$A_n,\enspace B_n,\enspace D_n,\enspace E_6,\enspace E_7, \enspace E_8,
\enspace F_4, \enspace\hbox{and}\enspace G_2=I_2(6).$$

\medskip
{\bf A complex reflection group is a group generated by 
complex reflections, 
i.e. invertible linear transformations of $\C^n$ which have finite
order and which have exactly one eigenvalue that is not $1$.}
Every finite Coxeter group is also a finite complex reflection group.
The finite complex reflection groups have been classified
by Shephard and Todd [ST] and each such
group is one of the groups
\begin{enumerate}
     \item[(a)] $G(r,p,n)$, where $r,p,n$ are positive integers such that
$p$ divides $r$, or
     \item[(b)] one of 34 ``exceptional'' finite complex reflection groups.
\end{enumerate}
The groups $G(r,p,n)$ are very similar to the symmetric group $S_n$ 
in many ways and this is probably why generalizing 
the $S_n$ theory has been so successful for these groups.
The symmetric groups and the
finite Coxeter groups of types $B_n$, and $D_n$ are all special
cases of the groups $G(r,p,n)$.

\bigskip\medskip\noindent
\leftline{\it The ``Hecke algebras'' of reflection groups}
\bigskip
{\bf The Iwahori-Hecke algebra of a finite Coxeter group $W$ 
is an algebra which is a $q$-analogue,
or $q$-deformation, of the group algebra $W$.}  
See Appendix B3 for a proper definition of
this algebra.  It has a basis $T_w$, 
$w\in W$, (so it is the same dimension as the group algebra of $W$) but
the multiplication in this algebra depends on a particular number 
$q\in \C$, which can be chosen arbitrarily.
These algebras are true Hecke algebras only when $W$ is a
finite Weyl group.

\medskip
{\bf The ``Hecke algebras'' of the groups $G(r,p,n)$ 
are $q$-analogues of the group algebras of the groups $G(r,p,n)$.}
It is only recently (1990-1994) that they have been defined. 
It is important to note that these algebras are not true Hecke algebras.  
In group theory, a Hecke algebra is a very specific kind
of double coset algebra and the ``Hecke algebras'' of the groups $G(r,p,n)$
do not fit this mold.  See Appendix B3-4 for the proper definition of a
Hecke algebra and some discussion of how the ``Hecke algebras'' 
for the groups $G(r,p,n)$ are defined.

\bigskip\bigskip\noindent
\leftline{\bf (2) Tensor power centralizer algebras }
\bigskip
{\bf A tensor power centralizer algebra is an algebra
which is isomorphic to $\End_G(V^{\otimes k})$
for some group (or Hopf algebra) $G$ and some representation $V$ of $G$.}  
In this definition 
$$\End_G(V^{\otimes k})=\{T\in \End(V^{\otimes k})\ |\ 
Tgv = gTv, \hbox{ for all $g\in G$ and all $v\in V^{\otimes k}$}\}.$$ 
There are some examples of tensor power centralizer algebras
that have been particularly important:
\smallskip\noindent
\begin{enumerate}
     \item[(a)] The group algebras, $\C S_k$, of the symmetric groups $S_k$, 
     \item[(b)] The Iwahori-Hecke algebras, $H_k(q)$, of type $A_{k-1}$,
     \item[(c)] The Temperley-Lieb algebras, $TL_k(x)$,
     \item[(d)] The Brauer algebras, $B_k(x)$,
     \item[(e)] The Birman-Murakami-Wenzl algebras, $BMW_k(r,q)$.
     \item[(f)] The spider algebras,
     \item[(f)] The rook monoid algebras,
     \item[(g)] The Solomon-Iwahori algebras,
     \item[(h)] The wall algebras, 
     \item[(i)]  The $q$-wall algebras,
     \item[(j)]  The partition algebras.
\end{enumerate}
We certainly do not have space to discuss all of these objects
in this paper, and thus we will limit ourselves to the 
cases (a)-(e) in our discussion below and in the Appendix, \S B5-8.
References for the remaining cases are as follows.
\smallskip\noindent
{\it The spider algebras.}
These algebras were written down combinatorially and studied by
G. Kuperberg [Ku4-5].
\smallskip\noindent
{\it The rook monoid algebras.}
L. Solomon (unpublished work) recognized that this very 
natural monoid (the combinatorics of which
has also been studied in [GR]) appears as a 
tensor power centralizer.
\smallskip\noindent
{\it The Solomon-Iwahori algebras.}
This algebra was introduced in [So1].  The fact that it is a 
tensor power centralizer algebra is an unpublished result
of L. Solomon, see [So2].
\smallskip\noindent
{\it The wall algebras.} These algebras were introduced in a nice combinatorial
form in [BC] and in other forms in [Ko] and Procesi [Pr] and
other older invariant theory works [Wy].  
All of these works were related to tensor power centralizers and/or
fundamental theorems of invariant theory.
\smallskip\noindent
{\it The q-wall algebras.}  These algebras were introduced by Kosuda
and Murakami [KM1-2] and studied subsequently in [Le] and [Ha1-2].
\smallskip\noindent
{\it The partition algebras.}
These algebras were introduced by V. Jones in [Jo1] and have been
studied subsequently by P. Martin [Ma].

\section*{\qquad Notes and references for answers to the main questions\qquad\null} 

Some partial results giving answers to the main questions for 
the complex reflection groups $G(r,p,n)$, their ``Hecke algebras'',
the Temperley-Lieb algebras, the Brauer algebras, 
and the Birman-Murakami-Wenzl
algebras can be found in Appendix B.
The appropriate references are as follows.

\bigskip\bigskip\noindent
\leftline{\bf (1) Complex reflection groups $G(r,p,n)$ and their
Hecke algebras}
\bigskip\medskip\noindent
\leftline{\it The complex reflection groups $G(r,p,n)$}
\bigskip

\noindent
{\bf I. What are the irreducible modules?}
\medskip
The indexing, dimension formulas and character formulas
for the representations
of the groups $G(r,p,n)$ are originally due to 
\smallskip
\halign{\indent\indent#\hfil\qquad&#\hfil\cr
Young [Y1] &for finite Coxeter groups of types $B_n$ and $D_n$, and\cr
Specht [Spc] &for the group $G(r,1,n)$.\cr}
\smallskip
\noindent
We do not know who first did the general $G(r,p,n)$ case but it is 
easy to generalize Young and Specht's results to this case.  
See [Ari] and [HR2] for recent accounts.

Essentially what one does to determine the indexing, 
dimensions and the characters of the irreducible
modules is to use Clifford theory
to reduce the $G(r,1,n)$ case to the case of the symmetric
group $S_n$.  Then one can use Clifford theory again
to reduce the $G(r,p,n)$ case to the $G(r,1,n)$ case.
The original reference for Clifford theory is [Cl] 
and the book by Curtis and Reiner [CR2] has a modern treatment.
The articles [Ste3] and [HR2] explain how the reduction
from $G(r,p,n)$ to $G(r,1,n)$ is done.
The dimension and character theory for the case $G(r,1,n)$
has an excellent modern treatment in [Mac], Appendix B to Chapter I.

\bigskip\noindent
{\bf C. How do we construct the irreducible modules?}  
\medskip
The construction
of the irreducible representations by Young symmetrizers was extended
to the finite Coxeter groups of types $B_n$ and $D_n$ by Young
himself in his paper [Y3].  
The authors don't know when the general case was first
treated in the literature, but it is not difficult to extend Young's results
to the general case $G(r,p,n)$.  
The $G(r,p,n)$ case does appear periodically
in the literature, see for example [Al].

Young's seminormal construction
was generalized to the ``Hecke algebras'' of $G(r,p,n)$ 
in the work of Ariki and Koike [AK] and Ariki [Ari].
One can easily set $q=1$ in the constructions of Ariki and Koike 
and obtain the appropriate analogues for the groups $G(r,p,n)$.
We do not know if the analogue of Young's seminormal construction
for the groups $G(r,p,n)$ appeared in the literature previous to the
work of Ariki and Koike on the ``Hecke algebra'' case.

\bigskip\noindent
{\bf S. Special/Interesting representations.}
\medskip
The authors do not know if 
the analogues of the $S_n$ results, {\bf (S1-3)} of Section 2,
have explicitly appeared 
in the literature.  It is easy to use symmetric functions and the character
formulas of Specht, see [Mac] Chpt.\ I, App.\ B, to derive 
formulas for the $G(r,1,n)$ case in terms of the symmetric group results.
Then one proceeds as described above to compute the
necessary formulas for $G(r,p,n)$ in terms of the $G(r,1,n)$ results.
See [Ste3] for how this is done.

\bigskip\medskip\noindent
\leftline{\it The ``Hecke algebras'' of reflection groups}
\bigskip\noindent
{\bf The definition.}
\medskip
The ``Hecke algebras'' corresponding to the groups $G(r,p,n)$
were defined by
\smallskip
\halign{\indent\indent#\hfil\qquad&#\hfil\cr
Ariki and Koike [AK], &for the case $G(r,1,n)$, and \cr
Brou\'e and Malle [BM] and Ariki [Ari], &for the general case $G(r,p,n)$. \cr}
\smallskip
See Appendix B4 for a definition of these algebras
and some partial answers to the main questions.
Let us give references to the literature for the answers to
the main questions for these algebras.

\bigskip\noindent
{\bf I. What are the irreducibles?}
\medskip
The results of Ariki-Koike [AK] and Ariki [Ari] say that the
``Hecke algebras'' of $G(r,p,n)$
are $q$-deformations of the group algebras of the
groups $G(r,p,n)$. Thus, it follows from the Tits deformation theorem
(see [Ca] Chapt 10, 11.2 and [CR2] \S 68.17) that
the indexings and dimension formulas for 
the irreducible representations of these algebras
must be the same as the indexings and dimension formulas
for the groups $G(r,p,n)$.
Finding analogues of the character formulas requires
a bit more work and a Murnaghan-Nakayama type rule for the
``Hecke algebras'' of $G(r,p,n)$ was given by Halverson and Ram [HR2].
As far as we know, the formula for the irreducible
characters of $S_n$ as a weighted sum of standard tableaux which we gave
in the symmetric group section has not yet been
generalized to the case of $G(r,p,n)$ and its ``Hecke
algebras''.

\bigskip\noindent
{\bf C. How do we construct the irreducible modules? }
\medskip
Analogues of Young's seminormal representations have been 
given by
\begin{enumerate}
     \item[] Hoefsmit [Hfs] and Wenzl [Wz1], independently, 
for Iwahori-Hecke algebras of type $A_{n-1}$,
     \item[] Hoefsmit [Hfs], for Iwahori-Hecke algebras of types $B_n$ and $D_n$,
     \item[] Ariki and Koike [AK] for the ``Hecke algebras'' of $G(r,1,n)$, and
     \item[] Ariki [Ari] for the general ``Hecke algebras'' of $G(r,p,n)$.
\end{enumerate}
There seems to be more than one appropriate choice
for the analogue of Young symmetrizers for Hecke algebras.
The definitions in the literature are due to
\smallskip
{\small
\halign{\indent\indent#\hfil\qquad&#\hfil\cr
Gyoja [Gy], for the Iwahori-Hecke algebras of type $A_{n-1}$, \cr
Dipper and James [DJ1] and Murphy [M1-2] for the
Iwahori-Hecke algebras of type $A_{n-1}$, \cr
King and Wybourne [KW] and Duchamp, et al [DK] 
for the Iwahori-Hecke algebras of type $A_{n-1}$, \cr
Dipper, James, and Murphy [DJ2], [DJM] for the Iwahori-Hecke algebras
of type $B_n$, \cr
Pallikaros [Pa], for the Iwahori-Hecke algebras of type $D_n$.\cr
Mathas [Mth] and Murphy [M3], for the ``Hecke algebras'' of $G(r,p,n)$. \cr}}
\smallskip
\noindent
The paper [GL] also contains important ideas in this direction.

\bigskip\noindent
{\bf S. Special/Interesting representations.}
\medskip
It follows from the Tits deformation theorem (or rather,
an extension of it) that the results for the ``Hecke algebras''
of $G(r,p,n)$ must be the same as for the case of the groups $G(r,p,n)$.

\bigskip\bigskip\noindent
\leftline {\bf (2) Tensor power centralizer algebras}
\bigskip

\noindent
{\bf The definitions.}
\medskip

The references for
the combinatorial definitions of the various centralizer algebras
are as follows. 
\medskip\noindent
{\it Temperley-Lieb algebras.}
These algebras are due, independently, to many different
people.  Some of the discoverers were
Rumer-Teller-Weyl [RTW], Penrose [P1-2], Temperley-Lieb [TL],
Kaufmann [Ka] and Jones[Jo2].  
The work of V. Jones was crucial in
making them so important in combinatorial representation theory today.

\medskip\noindent
{\it The Iwahori-Hecke algebras of type $A_{n-1}$.}
Iwahori [Iw] introduced these algebras in 1964 in connection with 
$GL(n,\F_q)$. Jimbo [Ji] realized that they
arise as tensor power centralizer algebras for quantum groups.
\medskip\noindent
{\it Brauer algebras.}
These algebras were defined by Brauer in 1937 [Br]. Brauer
also proved that they are tensor power centralizers.
\medskip\noindent
{\it Birman-Murakami-Wenzl algebras.} These algebras are due
to Birman and Wenzl [BW] and Murakami [Mu1].  It was realized early
[Re] [Wz3] that these arise as tensor power centralizers but there was no
proof in the literature for some time. 
See the references in [CP] \S 10.2.

\bigskip\noindent
{\bf I. What are the irreducibles?}
\medskip

Indexing of the representations of tensor power centralizer
algebras follows from double centralizer theory (see Weyl [Wy1])
and a good understanding of the indexings and tensor product
rules for the group or algebra which it is centralizing
(i.e. $GL(n,\C)$, $O(n,\C)$, $U_q{\goth{sl}}(n)$, etc.).
The references for 
resulting indexings and dimension formulas for the irreducible
representations are as follows:
\medskip\noindent
{\it Temperley-Lieb algebras.} These results are classical and can 
be found in the book by Goodman, de la Harpe, and Jones [GHJ].
\medskip\noindent
{\it Brauer algebras.} These results were known to Brauer [Br] 
and Weyl [Wy]. An important combinatorial point of view was
given by Berele [Be1-2] and further developed by Sundaram [Su1-3].
\medskip\noindent
{\it Iwahori-Hecke algebras of type $A_{n-1}$.}
These results follow
from the Tits deformation theorem and the
corresponding results for the symmetric group. 
\medskip\noindent
{\it Birman-Murakami-Wenzl algebras.}
These results follow
from the Tits deformation theorem and the
corresponding results for the Brauer algebra.
\medskip\noindent
The indexings and dimension formulas for the
Temperley-Lieb and Brauer algebras also follow
easily by using the techniques of the Jones basic construction,
see [Wz2] and [HR1].

\bigskip
The references for the irreducible characters of the various
tensor power centralizer algebras are as follows:
\medskip\noindent
{\it Temperley-Lieb algebras}. Character formulas can be derived 
easily by using 
Jones Basic Construction techniques [HR1].
\medskip\noindent
{\it Iwahori-Hecke algebras of type $A_{n-1}$}. 
The analogue of the formula for the
irreducible characters of $S_n$ as a weighted sum of standard 
tableaux was found by Roichman [Ro].  
Murnaghan-Nakayama type formulas were found by several authors
[KW], [vdJ], [VK], [SU], and [Ra2].  
\medskip\noindent
{\it Brauer algebras and Birman-Murakami-Wenzl algebras}. 
Murnaghan-Nakayama type formulas were derived in
[Ra1] and [HR1], respectively.

\medskip\noindent
Brauer algebra and Birman-Murakami-Wenzl algebra
analogues of the formula for
the irreducible characters of the symmetric groups as a 
weighted sum of standard tableaux have not appeared in
the literature.

\bigskip\noindent
{\bf C. How do we construct the irreducibles?}

\medskip\noindent
{\it Temperley-Lieb algebras.}
An application of the Jones Basic Construction (see [Wz2] and [HR1])
gives a construction of the irreducible representations of the 
Temperley-Lieb algebras.  This construction is classical and
has been rediscovered by many people.  
In this case the construction is an 
analogue of the Young symmetrizer construction. 
The analogue of the seminormal construction
appears in [GHJ].

\medskip\noindent
{\it Iwahori-Hecke algebras of type $A_{n-1}$.}
The analogue of Young's seminormal construction for this case
is due, independently, to Hoefsmit [Hfs] and Wenzl [Wz1].
Analogues of Young symmetrizers (different analogues) have been given by
Gyoja [Gy] and Dipper, James, and Murphy [DJ1], [M1-2], 
King and Wybourne [KW] and Duchamp, et al. [DK].

\medskip\noindent
{\it Brauer algebras.}
Analogues of Young's seminormal representations have been given,
independently, by Nazarov [Nz] and Leduc and Ram [LR].  An analogue of
the Young symmetrizer construction can be obtained by applying the 
Jones Basic Construction to the classical Young symmetrizer construction
and this is the one that has been used by many authors [BBL], [HW], 
[Ke], [GL].
The actual element of the algebra which is the analogue of 
the Young symmetrizer
involves a central idempotent for which there is no known explicit formula
and this is the reason that most authors work with a quotient
formulation of the appropriate module.

\medskip\noindent
{\it Birman-Murakami-Wenzl algebras.}
Analogues of Young's seminormal representations have been given
by Murakami [Mu2]  and Leduc and Ram [LR].  The methods in the
two works are different, the work of Murakami
uses the physical theory of Boltzmann weights and the work of
Leduc and Ram uses the theory of ribbon Hopf algebras and quantum groups.
Exactly in the same way as for the Brauer algebra, an analogue of
the Young symmetrizer construction can be obtained by applying the 
Jones Basic Construction to the Young symmetrizer constructions
for the Iwahori-Hecke algebra of type $A_{n-1}$.
As in the Brauer algebra case one should work with a quotient 
formulation of the module
to avoid using a central idempotent for which there is no known 
explicit formula.

\bigskip\bigskip\noindent
\leftline{\bf (3) Reflection groups 
of exceptional type.}
\bigskip

Generalizing the $S_n$ theory to
finite Coxeter groups of 
exceptional type, finite complex reflection groups
of exceptional type and the corresponding Iwahori-Hecke
algebras, has been largely unsuccessful.  This is
not to say that there haven't been some very nice
partial results only that at the moment nobody
has any understanding of how to make a
good combinatorial theory to encompass all the classical and exceptional
types at once.
Two amazing partial results along these lines are
$$\hbox{the Springer construction\qquad and \qquad
the Kazhdan-Lusztig construction.}$$ 

\medskip
The Springer construction is a construction
of the irreducible representations
of the crystallographic reflection groups
on cohomology of unipotent varieties [Spr].
It is a geometric construction and not a combinatorial construction.
See Appendix A3 for more information in the symmetric group case.
It is possible that this construction may be combinatorialized in the
future, but to date no one has done this.

\medskip
The Kazhdan-Lusztig construction [KL1] is a construction of certain
representations called {\it cell representations} and it works for all
finite Coxeter groups.  The cell representations are almost
irreducible but unfortunately not irreducible in general, and
nobody understands how to break them up into irreducibles,
except in a case by case fashion.
%\footnote{*}{We have recently heard that
%G. Lusztig has solved this problem.}
The other problem with these
representations is that they depend crucially on certain polynomials,
the Kazhdan-Lusztig polynomials, which seem to be impossible
to compute or understand well except in very small cases,
see [Bre] for more information.
See [Ca] for a summary and tables of the known facts about representation
of finite Coxeter groups of exceptional type.

\bigskip\noindent
{\bf Remarks}  
\smallskip\noindent
\begin{enumerate}
\item[(1)] A Hecke algebra is a specific ``double coset
algebra'' which depends on a group $G$ and a subgroup $B$.
Iwahori [Iw] studied these algebras in the case that $G$ is a finite Chevalley
group and $B$ is a Borel subgroup of $G$ and defined
what are now called {\it Iwahori-Hecke algebras}.  These are
$q$-analogues of the group algebras of finite Weyl groups.
The work of Iwahori yields a presentation for these algebras which
can easily be extended to define Iwahori-Hecke algebras for all Coxeter
groups but, except for the original Weyl group case, these have never been
realized as true Hecke algebras, i.e. double coset algebras corresponding
to an appropriate $G$ and $B$.
The ``Hecke algebras'' corresponding to the groups
$G(r,p,n)$ are $q$-analogues of the group algebras
of $G(r,p,n)$.  Although these algebras are not true Hecke
algebras either,
Brou\'e and Malle [BM] have shown that
many of these algebras arise in connection with
non-defining characteristic representations of
finite Chevalley groups and Deligne-Lusztig
varieties.
\item[(2)]
There is much current research on generalizing
symmetric group results 
to affine Coxeter groups and affine Hecke algebras.
The case of affine Coxeter groups was done by Kato [Kat] 
using Clifford theory ideas.
The case of affine Hecke algebras has been intensely
studied by Lusztig [Lu1-7], Kazhdan-Lusztig [KL2],
and Ginzburg [G], [CG], but most of this work is very geometric
and relies on intersection cohomology/K-theory methods.
Hopefully some of their work will be made combinatorial
in the near future.
\item[(3)]  Wouldn't it be great if we had a nice
combinatorial representation theory for finite simple
groups!!
\end{enumerate}

\bigskip

%Section 6
\section{Generalizations of $GL(n,\C)$ results}

\bigskip
There have been successful generalizations of 
the $GL(n,\C)$ results for the questions 
{\bf (Ia-c), (S1), (S2)} to the following classes of groups
and algebras.

\medskip\noindent
\begin{enumerate}
\item[\bf (1)] {\bf Connected complex semisimple Lie groups}.
\smallskip\noindent
Examples: $SL(n,\C)$, \enspace $SO(n,\C)$, \enspace $Sp(2n,\C)$,
\enspace $PGL(n,\C)$, \enspace $PSO(2n,\C)$, \enspace $PSp(2n,\C)$.
\item[\bf (2)]  {\bf Compact connected real Lie groups.}
\smallskip\noindent
Examples: $SU(n,\C)$, \enspace $SO(n,\R)$, \enspace $Sp(n)$,
where $Sp(n)= Sp(2n,\C)\cap U(2n,\C)$.  
\item[\bf (3)] {\bf Finite dimensional complex semisimple Lie algebras.}
\smallskip\noindent
Examples: ${\goth{sl}}(n,\C)$, \enspace ${\goth{so}}(n,\C)$, \enspace
${\goth{sp}}(2n,\C).$
See Appendix A7 for the complete list of the finite dimensional
complex semisimple Lie algebras.
\item[\bf (4)] {\bf Quantum groups corresponding to complex semisimple
Lie algebras.}
\end{enumerate}

\bigskip
The method of generalizing the $GL(n,\C)$ results
to the objects in {\bf (1-4)}  is to reduce them all to
case {\bf (3)} and then solve case {\bf (3)}.  
The results for case {\bf (3)} are given in Section 4.
The reduction of cases {\bf (1)} and {\bf (2)} to case 
{\bf (3)} are outlined in [Se2],
and given in more detail in [Va] and [BtD].
The reduction of {\bf (4)} to {\bf (3)} is given in [CP] and in [Ja].

\section*{\qquad Partial results for further generalizations\qquad\null }

\medskip
Some partial results along the lines of the results
{\bf (Ia-c)} and {\bf (S1-2)} for $GL(n,\C)$ and complex
semisimple Lie algebras
have been obtained for the following groups and algebras.
{\bf \smallskip\noindent
\begin{enumerate}
     \item[(1)] Kac-Moody Lie algebras and groups
     \item[(2)] Yangians
     \item[(3)] Simple Lie superalgebras
\end{enumerate}
}

\medskip\noindent
Other groups and algebras, for which the combinatorial 
representation theory is not understood very well, are
{\bf 
\smallskip\noindent
\begin{enumerate}
     \item[(4)]  Finite Chevalley groups
     \item[(5)] ${\goth p}$-adic Chevalley groups
     \item[(6)]  Real reductive Lie groups
     \item[(7)] The Virasoro algebra
\end{enumerate}
}
\smallskip\noindent
There are many many possible ways that we could extend this list
but probably these four cases are the most fundamental cases where
the combinatorial representation theory has not been formulated.
There has been intense work on all of these cases, but hardly any by
combinatorialists. Thus there are many beautiful results known
but very few of them have been stated or interpreted through a
combinatorialists eyes.  The world is a gold mine, yet to be mined!

\bigskip\noindent
{\bf Notes and references}
\smallskip\noindent
\begin{enumerate}
\item[(1)]  An introductory reference to Kac-Moody Lie algebras is
[Kc] .  This book contains a good description of the basic representation
theory of these algebras.  We don't know of a good introductory reference
for the Kac-Moody groups case, we would suggest beginning with the
paper [KK] and following the references there.
\item[(2)]  The basic introductory
reference for Yangians and their basic representation
theory is [CP], Chapter 12.  See also the references given there.
\item[(3)]  The best introductory reference for Lie superalgebras is
Scheunert's book [Sch].  For an update on the combinatorial representation
theory of these cases see the papers [Srg], [BR], [BRS], and [Snv].
\item[(4)]  Finding a general combinatorial representation theory for
finite Chevalley groups has been elusive for many
years.  After the fundamental work of J.A. Green [Gr] in 1955
which established a 
combinatorial representation theory for $GL(n,\F_q)$ there has been a concerted
effort to extend these results to other finite Chevalley groups.  
G. Lusztig [Lu8-11] has made important contributions
to this field; in particular, the results of 
Deligne-Lusztig [DL] are fundamental.  
However, this is a geometric approach
rather than a combinatorial one and there is much work 
to be done for combinatorialists,
even in interpreting the known results from combinatorial viewpoint.  
A good introductory treatment of this 
theory is the book by Digne and Michel [DM].  The original work of Green
is treated in [Mac] Chapt. IV.
\item[(5)]  The representation theory of ${\goth p}$-adic Lie groups has
been studied intensely by representation theorists but essentially 
not at all by
combinatorialists.  It is clear that there is a beautiful (although possibly
very difficult) combinatorial representation theory lurking here.
The best introductory reference to this work is the paper of R. Howe [Ho]
on ${\goth p}$-adic $GL(n)$.  Recent results of G. Lusztig [Lu7]
are a very important step in providing a general combinatorial
representation theory for ${\goth p}$-adic groups.
\item[(6)]  The best place to read about
the representation theory of real reductive groups
is in the books of D. Vogan and N. Wallach 
[AV], [Vg1] , [Vg2], [Wa].
\item[(7)]  The Virasoro algebra is a Lie algebra that seems to turn up
in every back alley of representation theory.  One can only surmise
that it must have a beautiful combinatorial representation theory that is
waiting to be clarified.  A good place to read about the Virasoro algebra
is in [FF].
\end{enumerate}

\clearpage
\setcounter{section}{0}
\def\thesection{A\arabic{section}}

\part*{Appendix A}

\appsection{Basic Representation Theory}

\medskip
An {\it algebra} $A$ is a vector space over $\C$ with
a multiplication that is associative, distributive,
has an identity and satisfies the following equation
$$
(ca_1)a_2=a_1(ca_2)=c(a_1a_2),
\qquad
\hbox{for all $a_1,a_2\in A$ and $c\in \C$.}$$
An {\it $A$-module} is a vector space $M$ over $\C$ 
with an $A$-action
$$\matrix{
A\times M &\longrightarrow &M \cr
(a,m) &\longrightmapsto &am, \cr}$$
which satisfies
$$\eqalign{
1m&=m,\cr
a_1(a_2m)&=(a_1a_2)m,\cr
(a_1+a_2)m&=a_1m+a_2m,\cr
a(c_1m_1+c_2m_2)&=c_1(am_1)+c_2(am_2). \cr}
$$
for all $a, a_1,a_2\in A$, $m,m_1,m_2\in M$ and $c_1,c_2\in \C$.
We shall use the words module and {\it representation} interchangeably.

A module $M$ is {\it indecomposable} if there do not exist
non zero $A$-modules $M_1$ and $M_2$ such that
$$M\cong M_1\oplus M_2.$$
A module $M$ is {\it irreducible} or {\it simple} 
if the only submodules
of $M$ are the zero module $0$ and $M$ itself.
A module $M$ is {\it semisimple} if it
is the direct sum of simple submodules. 
%(Chapter 8, sec 3, $n^o 3$, [Bour]). 

An algebra is {\it simple} 
if the only ideals of $A$ are the zero ideal
$0$ and $A$ itself.
The {\it radical} $\rad(A)$ of an algebra $A$ is the intersection
of all the maximal left ideals of $A$. 
%(Chapter 8, sec. 6 $n^o 3$ of [Bour]). 
An algebra $A$ is 
{\it semisimple} if all its modules are semisimple. 
%(Chapter 8, sec 5, $n^o 1$,  [Bour]). 
An algebra $A$ is {\it Artinian} if every decreasing sequence of left ideals
of $A$ stabilizes, that is for every chain 
$$A_1\supseteq A_2\supseteq A_3 \supseteq \cdots$$ 
of left ideals of $A$ there is an integer $m$ such that
$A_i=A_m$ for all $i\geq m.$
%Definition on p. 157 of CR and p.24 of Bourbaki.

The following statements follow directly from the definitions.

\medskip\noindent
{\sl
Let $A$ be an algebra.
\smallskip
\begin{enumerate}
\item[(a)]  Every irreducible $A$-module is indecomposable.
\smallskip
\item[(b)]  The algebra $A$ is semisimple if and only if
every indecomposable $A$-module is irreducible.
\end{enumerate}
}

\medskip\noindent
The proofs of the following statements are more involved and can be found in
[Bou2] Chpt. VIII, \S 6, $n^o 4$ and \S 5, $n^o 3$. 

\appthm 
\smallskip\noindent
\begin{enumerate}
\item[(a)] 
If $A$ is an Artinian algebra then the radical of $A$ is the largest nilpotent
ideal of $A$. 
%(Chapter 8, sec 6 $n^o4$ [Bour]). 
\item[(b)]
%(Chapter 8, sec 6, $n^o 4$, [Bour]),
An algebra $A$ is semisimple if and only if $A$ is Artinian and
$\rad(A) = 0.$
\item[(c)] 
%([Bour], Chapter 8, sec 5, $n^o 3$.) 
Every semisimple algebra is a direct sum of simple algebras.
\end{enumerate}
\endthm

The case when $A$ is not necessarily semisimple
is often called {\it modular representation theory}.
Let $M$ be an $A$-module.  
A {\it composition series} of $M$ is a chain 
$$M=M_k\supseteq M_{k-1}\supseteq \cdots \supseteq M_{1}
\supseteq M_0=0,$$
such that, for each $1\le i\le k$, the modules
$M_i/M_{i-1}$ are irreducible.  The irreducible
modules $M_i/M_{i-1}$ are the {\it factors} of the
composition series.
The following theorem is proved in [CR1] (13.7).

\appthm  (Jordan-H\"older)
If there exists a composition series for $M$
then any two composition series must have the same
multiset of factors (up to module isomorphism).
\endthm

\noindent
An important combinatorial point of view
is as follows:
The analogue of the subgroup lattice of a group
can be studied for any $A$-module $M$.  
More precisely, 
the {\it submodule lattice} $L(M)$ of $M$ is the lattice
defined by the submodules of $M$ with the
order relations given by inclusions of submodules.
The composition series are maximal chains in this
lattice.   

\bigskip\noindent
{\bf References}  
\smallskip\noindent
All of the above results can be found in
[Bou2] Chapt. VIII and [CR1].

\appsection{Partitions and tableaux}

\subsubsection*{Partitions}

\medskip
A {\it partition} is a sequence $\lambda = (\lambda_1,\ldots,\lambda_n)$
of integers such that $\lambda_1\ge \cdots \ge\lambda_n\ge 0$.
It is conventional to identify a partition with its Ferrers diagram
which has $\lambda_i$ boxes in the $i$th row.
For example the partition $\lambda = (55422211)$ has Ferrers diagram
\smallskip
$$
\displaylines{
\shape \cr
\lambda = (55422211) \cr
}
$$
We number the rows and columns of the Ferrers diagram
as is conventionally done for matrices.
If $x$ is a box in $\lambda$ then
the {\it content} and the {\it hook length} of $x$ are respectively given by
$$\matrix{
c(x) = j-i, 
&\quad 
&\hbox{if $x$ is in position $(i,j)\in \lambda$,\enspace and }\cr
h_x = \lambda_i-i+\lambda_j'-j+1,
&&\hbox{
where  $\lambda_j'$ is the length of the $j$th column of
$\lambda$.}
\cr
}$$
\medskip
$$
\matrix{
\contents &\qquad\qquad\qquad &\hooks \cr 
\hbox{Contents of the boxes} & 
&\hbox{Hook lengths of the boxes} \cr }
$$

If $\mu$ and $\lambda$ are partitions such that the Ferrers
diagram of $\mu$ is contained in the Ferrers diagram of $\lambda$ 
then we write $\mu\subseteq\lambda$ and we denote the difference
of the Ferrers diagrams by $\lambda/\mu$.  We refer to $\lambda/\mu$
as a {\it shape} or, more specifically, a {\it skew shape}.
$$
\displaylines{
\skewshape \cr
\lambda/\mu = (55422211)/(32211) \cr
}
$$

\subsubsection*{Tableaux}

\medskip
Suppose that $\lambda$ has $k$ boxes.
A {\it standard tableau} of shape $\lambda$ is a filling 
of the Ferrers diagram of $\lambda$ with $1,2,\ldots,k$
such that the rows and columns are increasing from left to right
and from top to bottom respectively.
\bigskip
$$
\stdtab
$$
%\centerline{Standard tableau}
\bigskip

\noindent
Let $\lambda/\mu$ be a shape.  A {\it column strict tableau}
of shape $\lambda/\mu$ filled with $1,2,\ldots n$ is a filling of the
Ferrers diagram of $\lambda/\mu$ with elements of the set 
$\{1,2,\ldots,n\}$ such that the rows are weakly increasing
from left to right and the columns are strictly increasing from
top to bottom.
The {\it weight} of a column strict tableau $T$ is the
sequence of positive integers $\nu = (\nu_1,\ldots,\nu_n)$,
where $\nu_i$ is the number of $i$'s in $T$.
\bigskip
$$
\colsttab
$$
%\centerline{Column strict tableau}
\smallskip
\centerline{Shape $\lambda = (55422211)$}
\centerline{Weight $\nu = (33323122111)$} 

\bigskip\noindent
The {\it word} of a column strict tableau $T$ is the sequence
$w=w_1w_2\cdots w_p$
obtained by reading the entries of $T$ from right to left
in successive rows, starting with the top row.
A word $w=w_1\cdots w_p$ is a {\it lattice permutation} if for each
$1\le r\le p$ and each $1\le i\le n-1$ the number of occurrences of the symbol
$i$ in $w_1\cdots w_r$ is not less than the number of occurences 
of $i+1$ in $w_1\cdots w_r$.
\bigskip
$$
\matrix{
\notlp &\qquad\qquad\qquad &\latperm \cr 
w=1122143346578 && w=1122133456578 \cr
\hbox{Not a lattice permutation} 
&&\hbox{Lattice permutation} \cr 
}$$

\bigskip
A {\it border strip} is a skew shape $\lambda/\mu$ which
is
\smallskip\noindent
\begin{enumerate}
\item[(a)] connected (two boxes are connected if they share an edge), and 
\item[(b)] does not contain a $2\times2$ block of boxes.
\smallskip\noindent
The weight of a border strip $\lambda/\mu$ is given by
$$\wt(\lambda/\mu) = (-1)^{r(\lambda/\mu)-1},$$
where $r(\lambda/\mu)$ is the number of rows in $\lambda/\mu$.
\bigskip
$$
\border  
$$
$$
\eqalign{
\lambda/\mu &= (86333)/(5222) \cr
\wt(\lambda/\mu) &= (-1)^{5-1} \cr }
$$
\end{enumerate}

\bigskip\noindent
Let $\lambda$ and $\mu=(\mu_1,\ldots,\mu_{\ell})$ be 
partitions of $n$.
A {\it $\mu$-border strip tableau} of shape $\lambda$ is
a sequence of partitions
$$T = (\emptyset = \lambda^{(0)}
\subseteq \lambda^{(1)} \subseteq \cdots\subseteq
\lambda^{(\ell-1)}\subseteq \lambda^{(\ell)} = \lambda)$$
such that, for each $1\le i\le \ell$, 
\smallskip\noindent
\begin{enumerate}
\item[(a)] $\lambda^{(i)}/\lambda^{(i-1)}$ is a border strip, and
\item[(b)] $|\lambda^{(i)}/\lambda^{(i-1)}| = \mu_i$.
\smallskip\noindent
The weight of a $\mu$-border strip tableau $T$ of shape $\lambda$ is
$$\wt(T) = \prod_{i=1}^{\ell-1} \wt(\lambda^{(i)}/\lambda^{(i-1)}).
\appformula$$
\end{enumerate}

\appthm  (Murnaghan-Nakayama rule)  
Let $\lambda$ and $\mu$ be partitions of $n$ and let
$\chi^\lambda(\mu)$ denote the irreducible character
of the symmetric group $S_n$ indexed by $\lambda$ evaluated 
at a permutation of cycle type $\mu$.
Then
$$\chi^\lambda(\mu) = \sum_T \wt(T),$$
where the sum is over all $\mu$-border strip tableaux $T$ of
shape $\lambda$ and $\wt(T)$ is as given in (A2.1).
\endthm

\bigskip\noindent
{\bf References}
\smallskip\noindent
All of the above facts can be found in [Mac] Chapt. I.
The proof of theorem (A2.2) is given in [Mac] Ch. I \S 7, Ex. 5.

\appsection{The flag variety, unipotent varieties, and Springer
theory for $GL(n,\C)$}

\subsubsection*{Borel subgroups, Cartan subgroups, and unipotent elements}

\medskip
The groups 
$$
B_n = \left\{ 
\pmatrix{ *      &*      &\cdots &*      \cr
          0      &*      &       &\vdots \cr
          \vdots &       &\ddots &*      \cr
          0      &\cdots &0      &*      \cr}
\right\},
$$
$$
T_n = \left\{ 
\pmatrix{ *      &0      &\cdots &0      \cr
          0      &*      &       &\vdots \cr
          \vdots &       &\ddots &0      \cr
          0      &\cdots &0      &*      \cr}
\right\},
$$
$$
U_n = \left\{ 
\pmatrix{ 1      &*      &\cdots &*      \cr
          0      &1      &       &\vdots \cr
          \vdots &       &\ddots &*      \cr
          0      &\cdots &0      &1      \cr}
\right\}, 
$$
are the subgroups of $GL(n,\C)$ consisting of upper triangular, diagonal,
and upper unitriangular matrices, respectively.
\smallskip\noindent
A {\it Borel subgroup} of $GL(n,\C)$ is a subgroup which is conjugate to
$B_n$.
\smallskip\noindent
A {\it Cartan subgroup} of $GL(n,\C)$ is a subgroup which is conjugate to
$T_n$.
\smallskip\noindent
A matrix $u\in GL(n,\C)$ is {\it unipotent} 
if it is conjugate to an upper unitriangular matrix.

\subsubsection*{The flag variety}

\noindent
There is a one-to-one correspondence between each of the following 
sets:
\smallskip\noindent
\begin{enumerate}
     \item[(1)] ${\cal B}= \{\hbox{Borel subgroups of $GL(n,\C)$} \}$,
     \item[(2)] $G/B$, where $G=GL(n,\C)$ and $B=B_n$,
     \item[(3)] $\{ \hbox{flags $0\subseteq V_1\subseteq V_2
\subseteq \cdots \subseteq V_n = \C^n$ such that $\dim(V_i)=i$} \}$.
\end{enumerate}
Each of these sets naturally has the structure of a complex algebraic variety,  
which is called the {\it flag variety}.

\subsubsection*{The unipotent varieties}

Given a unipotent element $u\in GL(n,\C)$ with Jordan blocks given by
the partition $\mu=(\mu_1,\ldots,\mu_{\ell})$ of $n$, 
define an algebraic variety
$${\cal B}_\mu = 
{\cal B}_u = \{\hbox{Borel subgroups of $GL(n,\C)$ which contain $u$} \}.$$
By conjugation, the structure of the subvariety ${\cal B}_u$
of the flag variety depends only on the partition $\mu$.
Thus ${\cal B}_\mu$ is well defined, as an algebraic variety.

\subsubsection*{Springer theory}

It is a deep theorem of Springer [Spr] (which holds in the generality
of semisimple algebraic groups and their corresponding Weyl groups)
that there is an action of
the symmetric group $S_n$ on the cohomology $H^*({\cal B}_u)$ 
of the variety ${\cal B}_u$.  This action can be interpreted nicely as
follows.  The imbedding 
$${\cal B}_u \subseteq {\cal B}
\quad\hbox{induces a surjective map}\quad 
H^*({\cal B}) \longrightarrow H^*({\cal B}_u).$$
It is a famous theorem of Borel that
there is a ring isomorphism
$$H^*({\cal B}) ~\cong~ \C[x_1,\ldots,x_n]/ I^+,\appformula$$
where $I^+$ is the ideal generated by symmetric functions
without constant term.  It follows that
$H^*({\cal B}_u)$ is also a quotient of $\C[x_1,\ldots,x_n]$.
From the work of Kraft [Kr], DeConcini and Procesi [DP]
and Tanisaki [Ta], one has that the ideal ${\cal T}_u$
which it is necessary to quotient by in order to obtain an isomorphism
$$H^*({\cal B}_u) ~\cong~ \C[x_1,\ldots,x_n]/{\cal T}_u,$$
can be described explicitly.

The symmetric group $S_n$ acts on the polynomial ring
$\C[x_1,\ldots, x_n]$ by permuting the variables. It turns out
that the ideal ${\cal T}_u$ remains invariant under this action, 
thus yielding a well defined action of $S_n$ on
$\C[x_1,\ldots,x_n]/{\cal T}_u.$
This action coincides with the Springer action on
$H^*({\cal B}_u)$.  Hotta and Springer [HS] have established that, 
if $u$ is a unipotent element of shape $\mu$ then,
for every permutation $w\in S_n$,
$$\sum_i q^i \varepsilon(w)\trace(w^{-1},H^{2i}({\cal B}_u))
=\sum_{\lambda\vdash n} \tilde K_{\lambda\mu}(q)\chi^\lambda(w),$$
where 
\smallskip
\begin{enumerate}
     \item[] $\varepsilon(w)$ is the sign of the permutation $w$, 
     \item[] $\trace(w^{-1},H^{2i}({\cal B}_u))$ is the trace 
of the action of $w^{-1}$ on $H^{2i}({\cal B}_u)$, 
     \item[]  $\chi^\lambda(w)$ is the irreducible character 
of the symmetric group evaluated at $w$, and 
     \item[]  $\tilde K_{\lambda\mu}(q)$ is a variant of
the Kostka-Foulkes polynomial,
see [Mac] III \S 7 Ex. 8, and \S 6.
\end{enumerate}

\medskip\noindent
It follows from this discussion and some basic facts about the
polynomials $\tilde K_{\lambda\mu}(q)$ that the top degree
cohomology group in $H^*({\cal B}_\mu)$ is a realization of the
irreducible representation of $S_n$ indexed by $\mu$,
$$S^\mu ~\cong~ H^{\rm top}({\cal B}_\mu).$$
This construction of the irreducible modules of $S_n$ is the
{\it Springer construction}.

\bigskip\noindent
{\bf References}
\smallskip\noindent
See [Mac] II \S 3 Ex. 1 for a description of the variety
${\cal B}_u$ and its structure.
The theorem of Borel stated in (A3.1) is given in [Bo]
and [BGG].  The references quoted in the text above
will provide a good introduction to the Springer theory.
The beautiful combinatorics of Springer theory has
been studied by Barcelo [Ba], Garsia-Procesi [GP], 
Lascoux [L], Lusztig [Lu12], Shoji [Shj], Spaltenstein [Sp], 
Weyman [Wm], and others. 

\appsection{Polynomial and rational representations of $GL(n,\C)$}

\medskip
If $V$ is a $GL(n,\C)$-module of dimension $d$ then,
by choosing a basis of $V$, we can define a map
$$\matrix{
\rho_V\colon GL(n,\C) &\longrightarrow &GL(d,\C) \cr
g &\longrightmapsto &\rho(g), \cr}$$
where $\rho(g)$ is the transformation of $V$ that is
induced by the action of $g$ on $V$.  
Let 
\begin{enumerate}
     \item[] $g_{ij}$ denote the $(i,j)$ entry of the matrix $g$,
and
     \item[] $\rho(g)_{kl}$ denote the $(k,l)$ entry of the matrix $\rho(g)$.
\end{enumerate}
The map $\rho$ depends on the choice of the basis of $V$, but
the following definitions do not.

\medskip\noindent
The module $V$ is a {\it polynomial representation} if
there are polynomials $p_{kl}(x_{ij})$, $1\le k,l\le d$,
such that 
$$\rho(g)_{kl}=p_{kl}(g_{ij}),\qquad\hbox{for all $1\le k,l\le d$.}$$
In other words $\rho(g)_{jk}$ is the same as the polynomial
$p_{kl}$ evaluated at the entries $g_{ij}$ of the matrix $g$.
\medskip\noindent
The module $V$ is a {\it rational representation} if
there are rational functions (quotients of two polynomials)
$p_{kl}(x_{ij})/q_{kl}(x_{ij})$, $1\le k,l\le d$,
such that 
$$\rho(g)_{kl}=p_{kl}(g_{ij})/q_{kl}(g_{ij}),\qquad
\hbox{for all $1\le k,l\le n$.}
$$

\noindent
Clearly, every polynomial representation is a rational one.

\medskip
The theory of rational representations of $GL(n,\C)$ can be
reduced to the theory of polynomial representations of $GL(n,\C)$.
This is accomplished as follows.
The determinant $\det\colon GL(n,\C)\to \C$
defines a $1$-dimensional (polynomial) representation of $GL(n,\C)$.
Any integral power
$$\matrix{
\hbox{$\det^k$}\colon &GL(n,\C) &\longrightarrow &\C \cr
&g &\longrightmapsto &\det(g)^k \cr
}$$
of the determinant also determines a $1$-dimensional representation
of $GL(n,\C)$.
{\bf
All irreducible rational representations $GL(n,\C)$ can be constructed
in the form
$$\hbox{$\det^k$}\otimes V^\lambda,$$
for some $k\in \Z$ and some irreducible polynomial representation
$V^\lambda$ of $GL(n,\C)$.  
}
\medskip
There exist representations of $GL(n,\C)$ which are not
rational representations, for example
$$g\mapsto \pmatrix{ 1 &\ln|\det(g)| \cr 0 &1\cr }.$$
There is no known classification of representations
of $GL(n,\C)$ which are not rational.

\bigskip\noindent
{\bf References}
\smallskip\noindent
See [Ste1] for a study of the combinatorics of the
rational representations of $GL(n,\C)$.

\appsection{Schur-Weyl duality and Young symmetrizers}

\medskip
Let $V$ be the usual $n$-dimensional representation
of $GL(n,\C)$ on column vectors of length $n$,
that is
$$V=\span\{b_1,\ldots, b_n\}
\quad\hbox{where}\quad b_i = (0, \ldots, 0, 1, 0, \ldots, 0)^t,$$
and the $1$ in $b_i$ appears in the $i$th entry.
Then
$$V^{\otimes k} = \span\{b_{i_1}\otimes \cdots \otimes b_{i_k}
\ |\ 1\le i_1,\ldots, i_k\le n \}
$$
is the span of the words of length $k$ in the letters $b_i$ 
(except that the letters are separated by tensor symbols).  
The general linear group $GL(n,\C)$ and the symmetric group $S_k$ act
on $V^{\otimes k}$ by 
$$g(v_1\otimes \cdots \otimes v_k)
=gv_1\otimes \cdots \otimes gv_k,
\quad\hbox{and}\quad
(v_1\otimes \cdots \otimes v_k)\sigma
= v_{\sigma(1)}\otimes \cdots \otimes v_{\sigma(k)},$$
where $g\in GL(n,\C)$, $\sigma\in S_k$, and 
$v_1,\ldots, v_k\in V$.
(We have chosen to make the $S_k$-action a right action
here, one could equally well choose the action of $S_k$
to be a left action but then the formula would be
$\sigma(v_1\otimes \cdots \otimes v_k)
= v_{\sigma^{-1}(1)}\otimes \cdots \otimes v_{\sigma^{-1}(k)}.$)

\smallskip
The following theorem is the amazing relationship between
the group $S_k$ and the group $GL(n,\C)$ which was discovered
by Schur [Sc1] and exploited with such success by Weyl [Wy1].

\appthm (Schur-Weyl duality)
\smallskip
\begin{enumerate}
\item[(a)]  The action of $S_k$ on $V^{\otimes k}$ generates
$\End_{Gl(n,\C)}(V^{\otimes k})$.
\item[(b)]  The action of $GL(n,\C)$ on $V^{\otimes k}$ generates
$\End_{S_k}(V^{\otimes k})$.
\end{enumerate}
\endthm

This theorem has the following important corollary, which provides
a intimate correspondence between the representation theory
of $S_k$ and {\it some} of the representations of $GL(n,\C)$
(the ones indexed by partitions of $k$).

\appcor As $GL(n,\C)\times S_k$ bimodules
$$V^{\otimes k} \cong \bigoplus_{\lambda\vdash k} V^\lambda\otimes S^\lambda,
$$
where $V^\lambda$ is the irreducible $GL(n,\C)$-module and 
$S^\lambda$ is the irreducible $S_k$-module indexed by $\lambda$.
\endcor

If $\lambda$ is a partition of $k$, then the irreducible 
$GL(n,\C)$-representation $V^\lambda$ is given by
$$V^\lambda\cong V^{\otimes k}P(T)N(T),$$
where $T$ is a tableau of shape $\lambda$ and 
$P(T)$ and $N(T)$ are as defined
in Section 2, Question {\bf C}.

\appsection{The Borel-Weil-Bott construction}

\medskip
Let $G=GL(n,\C)$ and let $B=B_n$ be the subgroup
of upper triangular matrices in $GL(n,\C)$.
A {\it line bundle on $G/B$}
is a pair $({\cal L}, p)$  where ${\cal L}$ is an algebraic variety
and $p$ is a
map (morphism of algebraic varieties)
$$p\;\colon {\cal L} \longrightarrow G/B,$$
such that the fibers of $p$ are lines and such that
${\cal L}$ is a locally trivial family of lines.
In this definition, {\it fibers} means
the sets $p^{-1}(x)$ for $x\in G/B$ and 
{\it lines} means one-dimensional
vector spaces. For the definition of {\it locally trivial
family of lines} see [Sh] Chapt. VI \S 1.2.
By abuse of language, a line bundle $({\cal L}, p)$ is simply denoted 
by ${\cal L}$.
Conceptually, a line bundle on $G/B$ means that we are
putting a one-dimensional vector space over each point in
$G/B$.
\smallskip\noindent
A {\it global section} of the line bundle ${\cal L}$ is a
map (morphism of algebraic varieties) 
$$s\;\colon G/B\to {\cal L}$$
such that $p\circ s$ is the identity map on $G/B$.
In other words a global section is any possible ``right inverse map''
to the line bundle.

Each partition $\lambda=(\lambda_1,\ldots,\lambda_n)$
determines a character (i.e. $1$-dimensional
representation) of the group $T_n$ of diagonal
matrices in $GL(n,\C)$ via
$$\lambda \left( 
\pmatrix{ t_1      &0      &\cdots &0      \cr
          0        &t_2    &       &\vdots \cr
          \vdots   &       &\ddots &0      \cr
          0        &\cdots &0      &t_n    \cr}
\right) 
=
t_1^{\lambda_1} t_2^{\lambda_2} \cdots t_n^{\lambda_n}.$$
Extend this character to be a character of $B=B_n$ by letting
$\lambda$ ignore the strictly upper triangular part
of the matrix, that is $\lambda(u)=1$,  for all $u\in U_n$.
Let ${\cal L}_{\lambda}$ be the fiber product $G\times_B \lambda$,
i.e. the set of
equivalence classes of pairs $(g,c)$, $g\in G$, $c\in \C^*$,
under the equivalence relation
$$(gb,c) \sim (g,\lambda(b^{-1})c),\qquad
\hbox{for all $b\in B$.}$$
Then ${\cal L}_{\lambda}= G\times_B \lambda$
with the map 
$$\matrix{
p\colon &{G\times_B \lambda} &\longrightarrow &G/B\cr
&(g,c) &\longrightmapsto &gB \cr}$$
is a line bundle on $G/B.$

The Borel-Weil-Bott theorem says that the irreducible
representation $V^{\lambda}$ of $GL_n(\C)$ is
$$V^\lambda \cong H^0(G/B, {\cal L}_\lambda),$$
where 
$H^0(G/B,{\cal L}_\lambda)$ is the space
of global sections of the line bundle ${\cal L}_\lambda$.

\bigskip\noindent
{\bf References}
\smallskip\noindent
See [FH] and G. Segal's article in [CMS] for further information and
references on this very important construction.

\appsection{Complex semisimple Lie algebras}

\medskip
A {\it finite dimensional complex semisimple Lie algebra} 
is a finite dimensional Lie algebra ${\goth g}$ over $\C$ 
such that $\rad({\goth g})=0$.
The following theorem classifies all finite dimensional
complex semisimple Lie algebras.

\appthm
\smallskip\noindent
\begin{enumerate}
\item[(a)]  Every finite dimensional complex semisimple Lie 
algebra ${\goth g}$ is a direct sum of complex simple Lie algebras.
\item[(b)]
There is one complex simple Lie algebra corresponding to
each of the following types
$$A_{n-1},\enspace B_n,\enspace C_n,\enspace D_n,\enspace 
E_6,\enspace E_7,\enspace E_8,\enspace F_4,\enspace G_2.$$ 
\end{enumerate}
\endthm

\noindent
The complex simple Lie algebras of types $A_n$, $B_n$, $C_n$ and $D_n$
are the ones of {\it classical type} and they are 
$$\matrix{
\hbox{Type $A_{n-1}$:} &\qquad &{\goth{sl}}(n,\C) = 
\{ A\in M_n(\C) \ |\ \Tr(A) = 0\}, \cr
\cr
\hbox{Type $B_n$:} &\qquad &{\goth{so}}(2n+1,\C) = 
\{ A\in M_{2n+1}(\C) \ |\ A+A^t = 0\}, \cr
\cr
\hbox{Type $C_n$:} &\qquad &{\goth{sp}}(2n,\C) = 
\{ A\in M_{2n}(\C) \ |\ AJ+JA^t = 0\}, \cr
\cr
\hbox{Type $D_n$:} &\qquad &{\goth{so}}(2n) = 
\{ A\in M_{2n}(\C) \ |\ A+A^t = 0\}, \cr
}$$
where $J$ is the matrix of a skew-symmetric form on a
$2n$-dimensional space.

\smallskip
Let ${\goth g}$ be a complex semisimple Lie algebra.  A 
{\it Cartan subalgebra} of ${\goth g}$ is a maximal
abelian subalgebra ${\goth h}$ of ${\goth g}$.
Fix a Cartan subalgebra ${\goth h}$ of ${\goth g}$.
If $V$ is a finite dimensional ${\goth g}$-module and
$\mu\colon {\goth h}\to \C$
is any linear function, define
$$V_\mu = \{ v\in V \ |\  hv=\mu(h)v, ~~\hbox{for all $h\in {\goth h}$}\}.$$
The space $V_\mu$ is the {\it $\mu$-weight space} of $V$.
It is a nontrivial theorem (see [Se2]) that
$$V=\bigoplus_{\mu\in P} V_\mu,$$
where $P$ is a $\Z$-lattice in ${\goth h}^*$ which can be identified
with the $\Z$-lattice $P$ which is defined below in Appendix A8.
The vector space ${\goth h}^*$ is the space of linear functions
from ${\goth h}$ to $\C$.

Let $\C[P]$ be the group algebra of $P$.  It can be given explicitly
as
$$\C[P]~=~\hbox{$\C$-span}\{ e^\mu \ |\ \mu\in P\},
\quad\hbox{with multiplication $e^\mu e^\nu=e^{\mu+\nu}$, for $\mu,\nu\in P$,}
$$
where the $e^\mu$ are formal variables indexed by the elements of $P$.
The {\it character} of a ${\goth g}$-module is 
$$\chr(V)=\sum_{\mu\in P} \dim(V_\mu)e^\mu.$$

\bigskip\noindent
{\bf References}
\smallskip\noindent
Theorem (A7.1) is due to the founders of the theory,
Cartan and Killing, from the late 1800's.  The beautiful 
text of Serre [Se2] gives a review of the definitions
and theory of complex semisimple Lie algebras.
See [Hu1] for further details.

\appsection{Roots, weights and paths}

\medskip
To each of the ``types'', $A_n$, $B_n$, etc.,
there is an associated hyperplane arrangement ${\cal A}$ in
$\R^n.$ 
\smallskip
$$
\displaylines{
\hyparr \cr
\hbox{Hyperplane arrangement for $A_2$} \cr
}
$$
\smallskip\noindent
The space $\R^n$ has the usual
Euclidean inner product $\langle\;,\;\rangle$.
For each hyperplane in
the arrangement ${\cal A}$ we choose two vectors
orthogonal to the hyperplane and pointing in opposite directions.
This set of chosen vectors is called the {\it root system}
$R$ associated to ${\cal A}$.
\smallskip
$$
\displaylines{
\root \cr
\hbox{Root system for $A_2$} \cr
}
$$
\smallskip\noindent
There is a convention for choosing the lengths of these vectors
but we shall not worry about that here.

Choose a chamber (connected component)
$C$ of $\R^n \backslash \bigcup_{H\in {\cal A}} H$.
\smallskip
$$
\displaylines{
\chamber \cr
\hbox{A chamber for $A_2$} \cr
}
$$
\smallskip\noindent
For each root $\alpha \in R$ we say that $\alpha$ is
{\it positive} if it points toward the same side of the hyperplane
as $C$ is and {\it negative} if points toward the opposite
side.  
It is standard notation to write
$$
\alpha >0, \qquad\hbox{if $\alpha$ is a positive root, \enspace and}\qquad
\alpha <0,  \qquad\hbox{if $\alpha$ is a negative root.}
$$
The positive roots which are associated to hyperplanes
which form the walls of $C$ are the {\it simple roots}
$\{\alpha_1,\ldots,\alpha_n\}$.
The {\it fundamental weights}
are the vectors $\{\omega_1,\ldots,\omega_n\}$
in $\R^n$ such that
$${\langle \omega_i, \alpha_j^\vee\rangle
= \delta_{ij}},
\quad
\hbox{where}\quad
\alpha_j^\vee
= {2\alpha_j \over \langle \alpha_j,\alpha_j\rangle}.$$
%The change from $alpha_j$ to $\alpha_j^\vee$ is
%just a renormalization.  Note that the $\{\omega_i\}$
%are the dual basis to the $\{\alpha_+j^\vee\}$.
Then $$P=\sum_{i=1}^r \Z\omega_i,
\qquad\hbox{and}\qquad
P^+=\sum_{i=1}^n \N\; \omega_i,
\qquad\hbox{where}\quad \N=\Z_{\ge 0},$$
are the lattice of {\it integral weights} and the cone of
{\it dominant integral weights}, respectively.
\smallskip
$$
\matrix{
\intweight &\qquad\qquad\qquad &\cone \cr 
\hbox{Lattice of integral weights} & 
&\hbox{Cone of dominant integral weights} \cr
}
$$
\smallskip\noindent
There is a one-to-one correspondence
between the irreducible representations of ${\goth g}$
and the elements of the cone $P^+$ in the lattice $P$.

Let $\lambda$ be a point in $P^+$.  Then the straight line
path from $0$ to $\lambda$ is the map
$$\matrix{
\pi_\lambda \colon &[0,1] &\longrightarrow &\R^n \cr
&t &\longrightmapsto &t\lambda. \cr
}$$
\smallskip
$$
\displaylines{
\path \cr
\hbox{Path from $0$ to $\lambda$} \cr
}
$$
\smallskip\noindent
The set ${\cal P}\pi_\lambda$ is given by
$${\cal P}\pi_\lambda =
\{ f_{i_1}\cdots f_{i_k}\pi_\lambda \ |\  1\le i_1,\ldots, i_k\le n \}$$
where $f_1,\ldots,f_n$ are the path operators introduced in [Li2].
These paths might look like
\smallskip
$$
\displaylines{
\piecewiselin \cr
\hbox{Path in ${\cal P}\pi_\lambda$} \cr
}
$$
\smallskip\noindent
They are always piecewise linear and end in a point in $P$.

\bigskip\noindent
{\bf References}
\smallskip\noindent
The basics of root systems can be found in [Hu1].  The
reference for the path model of Littelmann is [Li2].

\clearpage
\setcounter{section}{0}
\def\thesection{B\arabic{section}}

\part*{Appendix B}
\bppsection{Coxeter groups, groups generated by reflections, and Weyl groups}

\bigskip
A {\it Coxeter group} is a group $W$ presented by generators
$S=\{s_1,\ldots, s_n\}$ and relations
$$\let\{=.
\cases{
s_i^2=1,  &for $1\le i\le n$, \cr
(s_is_j)^{m_{ij}}=1,  &for $1\le i\ne j\le n$,\cr}
$$
where each $m_{ij}$ is either $\infty$ or a positive integer greater than 1.

\smallskip\noindent
A {\it reflection} is a linear transformation of $\R^n$
which is a reflection in some hyperplane.
\smallskip\noindent
A {\it finite group generated by reflections} is a finite subgroup
of $GL(n,\R)$ which is generated by reflections.

\bppthm The finite Coxeter groups are exactly the 
finite groups generated by reflections.
\endthm

\noindent
A finite Coxeter group is {\it irreducible} if it cannot be written
as a direct product of finite Coxeter groups.

\bppthm (Classification of finite Coxeter groups)
\smallskip\noindent
\begin{enumerate}
\item[(a)] Every finite Coxeter group can be written as a direct product
of irreducible finite Coxeter groups.
\item[(b)] 
There is one irreducible finite Coxeter group
corresponding to each of the following ``types''
$$A_{n-1},\enspace B_n,\enspace D_n,\enspace E_6,\enspace E_7,
\enspace E_8,\enspace F_4,\enspace H_3,\enspace H_4,\enspace I_2(m).$$
\end{enumerate}
\endthm

The irreducible finite Coxeter groups of {\it classical type}
are the ones of types $A_{n-1}, B_n$, and $D_n$ and the others
are the irreducible finite Coxeter groups of {\it exceptional
type}.  
\smallskip\noindent
\begin{enumerate}
\item[(a)]  The group of type $A_{n-1}$ is the symmetric group
$S_n$.
\item[(b)]  The group of type $B_n$ is the hyperoctahedral group
$(\Z/2\Z)\wr S_n$, the wreath product of the group of order
$2$ and the symmetric group $S_n$. 
It has order $2^n n!$. 
\item[(c)]  The group of type $D_n$ is 
a subgroup of index 2 in the Coxeter group of type $B_n$.
\item[(d)]  The group of type $I_2(m)$ is a dihedral group of
order $2m$.
\end{enumerate}

\medskip\noindent
A finite group $W$ generated by reflections in $\R^n$
is {\it crystallographic}
if there is a lattice in $\R^n$ which is stable under the action
of $W$.  The crystallographic finite Coxeter groups are also called
{\it Weyl groups}.  The irreducible Weyl groups are the
irreducible finite Coxeter groups of types
$$A_{n-1},\enspace B_n,\enspace D_n,\enspace E_6,\enspace E_7,
\enspace E_8,\enspace F_4,\enspace G_2=I_2(6).$$

\bigskip\noindent
{\bf References}
\smallskip\noindent
The most comprehensive reference for finite groups generated
by reflections is [Bou1].  See also the book of Humphreys [Hu2].

\bppsection{Complex reflection groups}

\bigskip
A {\it complex reflection} is an invertible linear transformation
of $\C^n$ of finite order which has exactly one eigenvalue
that is not $1$.  A {\it complex reflection group} is a group generated
by complex reflections in $\C^n$.  The finite complex reflection groups
have been classified by Shepard and Todd [ST].    
Each finite complex reflection group is either
\smallskip\noindent
\begin{enumerate}
     \item[(a)] $G(r,p,n)$ for some positive integers $r,p,n$ such 
that $p$ divides $r$, or
     \item[(b)] one of 34 other ``exceptional'' finite complex reflection groups.
\end{enumerate}

Let $r,p,d$ and $n$ be positive integers such that $pd=r$.  
The complex reflection
group $G(r,p,n)$ is the set of $n\times n$ matrices such that
\begin{enumerate}
     \item[(a)] The entries are either $0$ or $r$th roots of unity,
     \item[(b)] There is exactly on nonzero entry in each row and each column,
     \item[(c)] The $d$th power of the product for the nonzero entries is $1$.
\end{enumerate}
The group $G(r,p,n)$ is a normal subgroup of $G(r,1,n)$ of index $p$ and
$$|G(r,p,n)|=dr^{n-1}n!.$$
In addition
\smallskip\noindent
\begin{enumerate}
     \item[(a)] $G(1,1,n)\cong S_n$ the symmetric group or 
{\it Weyl group of type $A_{n-1}$},
     \item[(b)] $G(2,1,n)$ is the {\it hyperoctahedral group} or 
{\it Weyl group of type $B_n$},
     \item[(c)] $G(r,1,n)\cong (\Z/r\Z)\wr S_n$, the wreath product 
of the cyclic group of order $r$ with $S_n$,
     \item[(d)] $G(2,2,n)$ is the {\it Weyl group of type $D_n$}.
\end{enumerate}

\subsection*{Partial results for $G(r,1,n)$}

The following are answers to the main 
questions {\bf (Ia-c)} for the groups $G(r,1,n)\cong (\Z/r\Z)\wr S_n$.  
For the general $G(r,p,n)$ case see [HR2].

\bigskip\noindent
{\bf I.  What are the irreducible $G(r,1,n)$-modules?}
\medskip\noindent
\begin{enumerate}
\item[(a)] How do we index/count them?
\smallskip
\begin{enumerate}
{\bf 
     \item[]
There is a bijection
$$
\matrix{
\hbox{ $r$-tuples $\lambda=(\lambda^{(1)},\ldots,\lambda^{(r)})$ 
of partitions} \cr
&\ \mapleftright{1-1}\ 
&\hbox{Irreducible representations $C^\lambda$.} \cr
\hbox{such that $\sum_{i=1}^r |\lambda^{(i)}|=n$} \cr
}
$$
}
\end{enumerate}

\item[(b)] What are their dimensions?
\smallskip\noindent
\begin{enumerate}
     \item[]
The dimension of the irreducible representation
$C^\lambda$ is given by
$$\eqalign{
{\bf dim}(C^\lambda) 
&= \hbox{{\bf \# of standard tableaux of shape $\lambda$}} \cr
&= n! \prod_{i=1}^r 
~\prod_{x\in \lambda^{(i)} } {1\over h_x}, \cr}$$
where $h_x$ is the hook length at the box $x$.
A standard tableau of shape $\lambda=(\lambda^{(1)},\ldots,\lambda^{(r)})$
is any filling of the boxes of the $\lambda^{(i)}$ with
the numbers $1,2,\ldots, n$ such that the rows and the columns of
each $\lambda^{(i)}$ are increasing.
\end{enumerate}
\item[(c)] What are their characters?
\begin{enumerate}
     \item[]  
A Murnaghan-Nakayama type rule for the characters of the 
groups $G(r,1,n)$ was originally given by Specht [Spc].
See also [Osi] and [HR2].
\end{enumerate}
\end{enumerate}

\bigskip\noindent
{\bf References}
\smallskip\noindent
The original paper of Shepard and Todd [ST] remains a 
basic reference.  Further information about these
groups can be found in [HR2].  The articles [OS],
[Leh], [Ste3], [Mal] contain other recent work on the combinatorics
of these groups.

\bppsection{Hecke algebras and ``Hecke algebras'' of Coxeter groups}

\bigskip
Let $G$ be a finite group and let $B$ be a subgroup of
$G$.
The {\it Hecke algebra} of the pair $(G,B)$ is the subalgebra
$$
{\cal H}(G,B) = \left\{ \sum_{g\in G} a_g g \ \big|\ 
a_g\in \C, \hbox{ and $a_g = a_h$ if $BgB=BhB$}.\right\}
$$
of the group algebra of $G$.
The elements
$$T_w ={1\over |B|}\sum_{g\in BwB} g,$$
as $w$ runs over a set of representatives of the double cosets
$B\backslash G/B$, form a basis of ${\cal H}(G,B)$.

\smallskip
Let $G$ be a finite
Chevalley group over the field $\F_q$ with $q$ elements
and fix a Borel subgroup $B$ of $G$.  The pair $(G,B)$
determines a pair $(W,S)$ where $W$ is the Weyl group
of $G$ and $S$ is a set of simple reflections in $W$
(with respect to $B$).
The {\it Iwahori-Hecke algebra} corresponding to 
$G$ is the Hecke algebra ${\cal H}(G,B).$
In this case the basis elements $T_w$ are indexed by the 
elements $w$ of the Weyl group $W$ corresponding to the
pair $(G,B)$ and the multiplication is given by
$$T_sT_w
= \cases{
T_{sw}, &if $\ell(sw)>\ell(w)$,\cr
\cr
(q-1)T_{w}+qT_{sw}, &if $\ell(sw)<\ell(w)$,\cr}$$
if $s$ is a simple reflection in $W$.  In this formula
$\ell(w)$ is the {\it length} of $w$, i.e. the minimum
number of factors needed to write $w$ as a product of
simple reflections.

\smallskip
A particular example of the Iwahori-Hecke algebra
occurs when $G=GL(n,\F_q)$ and $B$ is the
subgroup of upper triangular matrices.  Then
the Weyl group 
$W$, is  the symmetric group $S_n$,  and the
simple reflections in the set $S$ are the transpositions
$s_i=(i,i+1)$, $1\le i\le n-1$.
In this case the algebra ${\cal H}(G,B)$ is the
{\it Iwahori-Hecke algebra of type $A_{n-1}$} and
(as we will see later) 
can be presented by generators $T_1,\ldots, T_{n-1}$
and relations
\medskip
\vbox{\openup4pt
\halign{\indent\indent #\hfil\quad&#\hfil\quad&#\hfil\cr
&$T_iT_j=T_jT_i$, & for $|i-j|>1$,\cr
&$T_iT_{i+1}T_i = T_{i+1}T_iT_{i+1}$, & for $1 \le i\le n-2$,\cr
&$T_i^2 = (q-1)T_i+q$, & for $2\le i\le n$.\cr
}}

\medskip\noindent
See Section B5 for more facts about the Iwahori-Hecke 
algebras of type $A$.
In particular, these Iwahori-Hecke algebras also appear as
tensor power centralizer algebras, see Theorem B5.3.
This is some kind of miracle: the Iwahori-Hecke algebras
of type $A$ are the only Iwahori-Hecke algebras which
arise naturally as tensor power centralizers.

\smallskip
In view of the multiplication rules for the Iwahori-Hecke
algebras of Weyl groups it is easy to define a
``Hecke algebra'' for all Coxeter groups $(W,S)$, 
just by defining
it to be the algebra with basis $T_w$, $w\in W$, 
and multiplication
$$T_sT_w
= \cases{
T_{sw}, &if $\ell(sw)>\ell(w)$,\cr
\cr
(q-1)T_{w}+qT_{sw}, &if $\ell(sw)<\ell(w)$,\cr}$$
if $s \in S$.  These algebras are not true Hecke algebras
except when $W$ is a Weyl group.

\bigskip\noindent
{\bf References}
\smallskip\noindent
For references on Hecke algebras see [CR2] (Vol I, Section 11).
For references on Iwahori-Hecke algebras see
[Bou1] Chpt. IV \S2 Ex. 23-25,
[CR2] Vol. II \S 67, and [Hu2] Chpt. 7.
The article [Cu] is also very informative.

\bppsection{``Hecke algebras'' of the groups $G(r,p,n)$}

\bigskip
Let $q$ and $u_0, u_1, \ldots, u_{r-1}$ be indeterminates. 
Let $H_{r,1,n}$ be
the algebra over the field $\C(u_0, u_1,
\ldots, u_{r-1},q)$ given by generators $T_1, T_2, \ldots, T_n$ and
relations
\medskip
\vbox{\openup4pt
\halign{\indent\indent #\hfil\quad&#\hfil\quad&#\hfil\cr
(1)& $T_iT_j=T_jT_i,$ &for $|i-j|>1,$ \cr
(2)& $T_iT_{i+1}T_i=T_{i+1}T_iT_{i+1},$ &for  $2\le i\le n-1,$ \cr
(3)& $T_1T_2T_1T_2=T_2T_1T_2T_1,$ & \cr
(4)& $(T_1 - u_0) (T_1 - u_1) \cdots (T_1 - u_{r-1}) =0,$\cr
(5)& $(T_i-q)(T_i+q^{-1})=0,$ &for $2\le i\le n.$\cr
}}
\medskip\noindent
Upon setting $q=1$ and $u_{i-1} = \xi^{i-1}$, where $\xi$ is a
primitive $r$th  root of unity, one obtains the group algebra
$\C G(r,1,n)$.  
In the special
case where $r = 1$ and $u_0 = 1$, we have $T_1 = 1$, and $H_{1,1,n}$
is isomorphic to an Iwahori-Hecke algebra of type $A_{n-1}$.
The case $H_{2,1,n}$ when $r=2$, $u_0=p$, and $u_1=p^{-1}$, is isomorphic
to an Iwahori-Hecke algebra of type $B_n$.

Now suppose that $p$ and $d$ are positive integers such that
$pd=r$.
Let $x_0^{1/p},\ldots, x_{d-1}^{1/p}$ be indeterminates,
let $\varepsilon=e^{2\pi i/p}$ be a primitive $p$th root of unity
and specialize the variables $u_0,\ldots,u_{r-1}$ according to the
relation 
$$u_{\ell d + kp +1}=\varepsilon^\ell x_k^{1/p},$$
where the subscripts on the $u_i$ are taken mod $r$.
The {\it ``Hecke algebra'' $H_{r,p,n}$ 
corresponding to the group $G(r,p,n)$} 
is the subalgebra of $H_{r,1,n}$ generated by the elements
$$a_0=T^p_1,\quad a_1=T_1^{-1}T_2T_1,
\quad\hbox{and}\quad
a_i=T_i,\quad 2\le i\le n.$$
Upon specializing $x_k^{1/p} = \xi^{kp}$, where $\xi$ is a
primitive $r$th root of unity,
$H_{r,p,n}$ becomes the group algebra
$\C G(r,p,n)$.  Thus $H_{r,p,n}$ is a ``$q$-analogue'' of the
group algebra of the group $G(r,p,n)$.

\bigskip\noindent
{\bf References}
\smallskip\noindent
The algebras $H_{r,1,n}$ were first constructed by Ariki
and Koike [AK], and they were classified as cyclotomic Hecke
algebras of type $B_n$ by Brou\'e and Malle [BM] and the
representation theory of $H_{r,p,n}$ was studied by Ariki [Ari].
See [HR2] for information about the characters of these 
algebras.

\bppsection{The Iwahori-Hecke algebras $H_k(q)$ of type $A$}

\bigskip
A {\it $k$-braid\/} is viewed  as two rows of $k$ vertices,
one above the other, 
and $k$ strands that connect top vertices to bottom 
vertices in such a way that
each vertex is incident to precisely one strand.  
Strands cross over and under
each other in three-space as they pass from one vertex to the next.
$$t_1 ~=~~ \bdone \;, \qquad t_2 ~=~ \bdtwo \;.$$
We multiply $k$-braids $t_1$ and $t_2$ using the concatenation product
given by identifying the vertices in the top row of 
$t_2$ with the corresponding
vertices in the bottom row of $t_1$ to obtain
the product $t_1 t_2$.  
$$t_1t_2 ~=~~ \prodbd$$

Given a permutation $w\in S_k$ we will make a $k$-braid $T_w$ by
tracing the edges in order from left to right across the top row.
Any time
an edge that we are tracing crosses an edge that has been already traced
we raise the pen briefly so that the edge being traced goes under 
the edge which is already there.  Applying this process 
to all of the permutations in $S_k$
produces a set of $k!$ braids.  

$$
w ~=~~ \perm \qquad\qquad\qquad T_w ~=~~ \Tw 
$$

Fix $q\in \C$.
The {\it Iwahori-Hecke algebra $H_k(q)$ of type $A_{k-1}$\/} 
is the span of the $k!$ braids produced by tracing permutations in
$S_k$ with multiplication determined by the braid multiplication
and the following identity.
$$ 
\poscross ~~=~ (q-1)~\idcross ~~+~ q~\negcross ~.
$$
This identity can be applied in any local portion of the braid.

\bppthm The algebra $H_k(q)$ is the associative algebra over $\C$
presented by generators $T_1,\ldots, T_{k-1}$ and relations
\smallskip
\medskip
\vbox{\openup4pt
\halign{\indent\indent #\hfil\quad&#\hfil\quad&#\hfil\cr
&  $T_iT_j=T_jT_i$, & for $|i-j|>1$,\cr
&  $T_iT_{i+1}T_i = T_{i+1}T_iT_{i+1}$, & for $1 \le i\le n-2$,\cr
&  $T_i^2 = (q-1)T_i+q$, & for $2\le i\le n$.\cr
}}
\endthm

The Iwahori-Hecke algebra of type $A$ is a 
$q$-analogue of the group algebra of the symmetric group.
If we allow ourselves to be imprecise (about the limit) we can write
$$\lim_{q\to 1} H_k(q) = \C S_k.$$

Let $q$ be a power of a prime $G=GL(n,\F_q)$ where $\F_q$ 
is the finite field with $q$ elements.
Let $B$ be the subgroup of upper triangular matrices in $G$  
and let $\One_B^G$ be the 
trivial representation of $B$ induced to $G$, i.e.
the $G$-module given by 
$$\One_B^G = \hbox{$\C$-span}\{ gB \ |\ g\in G\},$$
where $G$ acts on the cosets by left multiplication.  
Using the description, see \S B3, of $H_n(q)$ as a double coset
algebra one gets an action of $H_n(q)$ on $\One_B^G$,
by right multiplication. This action commutes with the $G$ action.

\bppthm
\begin{enumerate}
\item[(a)] The action of $H_n(q)$ on $\One_B^G$ generates
$\End_G(\One_B^G).$
\item[(b)] The action of $G$ on $\One_B^G$ generates
$\End_{H_n(q)}(\One_B^G).$
\end{enumerate}
\endthm
\noindent
This theorem gives a ``duality'' between $GL(n,\F_q)$ and $H_n(q)$
which is similar to a Schur-Weyl duality, but it differs in a crucial way:
the representation $\One_B^G$ is not a tensor power representation, and
thus this is not yet realizing $H_n(q)$ as a tensor power centralizer.

\smallskip
The following result gives a true analogue of the Schur-Weyl duality
for the Iwahori-Hecke algebra of type $A$, it realizes $H_k(q)$
as a tensor power centralizer.
Assume that $q\in \C$ is not $0$ and is not a root of unity.  Let
$U_q{\goth{sl}}_n$ be the Drinfel'd-Jimbo quantum group
of type $A_{n-1}$ and let $V$ be the $n$-dimensional irreducible representation
of $U_q{\goth{sl}}_n$ with highest weight $\omega_1$.
There is an action, see [CP], of $H_k(q^2)$ on $V^{\otimes k}$ which commutes
with the $U_q{\goth{sl}}_n$ action.

\bppthm
\begin{enumerate}
\item[(a)] The action of $H_k(q^2)$ on $V^{\otimes k}$ generates
$\End_{U_q{\goth{sl}}_n}(V^{\otimes k}).$
\item[(b)] The action of $U_q{\goth{sl}}_n$ on $V^{\otimes k}$ generates
$\End_{H_k(q^2)}(V^{\otimes k}).$
\end{enumerate}
\endthm

\bppthm The Iwahori-Hecke algebra of type $A_{k-1}$, $H_k(q)$, is semisimple
if and only if $q\ne 0$ and $q$ is not a $j$th root of unity for any
$2\le j\le n$.
\endthm

\vskip-.1in
\subsection*{Partial results for $H_k(q)$}

\noindent
The following results giving answers to the main questions
{\bf (Ia-c)} for the Iwahori-Hecke algebras of type
$A$ hold when $q$ is such that $H_k(q)$ is semisimple.

\bigskip\noindent
{\bf I.  What are the irreducible $H_k(q)$-modules?}
\medskip\noindent
\begin{enumerate}
\item[(a)] How do we index/count them?
\smallskip
\begin{enumerate}
{\bf 
     \item[]
There is a bijection
$$
\hbox{Partitions $\lambda$ of $n$}
\qquad\mapleftright{1-1}\qquad
\hbox{Irreducible representations $H^\lambda$ }
$$
}
\end{enumerate}

\item[(b)] What are their dimensions?
\smallskip\noindent
\begin{enumerate}
     \item[]
The dimension of the irreducible representation
$H^\lambda$ is given by
$$
\eqalign{
{\bf dim}(H^\lambda)
&= \hbox{{\bf \# of standard tableaux of shape $\lambda$}}\cr
&= {n!\over \prod_{x\in \lambda} h_x },\cr}
$$
where $h_x$ is the hook length at the box $x$ in $\lambda$.
\end{enumerate}

\item[(c)] What are their characters?
\smallskip\noindent
\begin{enumerate}
     \item[]  For each partition 
$\mu=(\mu_1,\mu_2,\ldots,\mu_\ell)$ of $k$
let $\chi^\lambda(\mu)$ be the character of the irreducible
representation $H^\lambda$ evaluated at the element $T_{\gamma_\mu}$
where $\gamma_\mu$ is the permutation 
$$\gamma_\mu = \underbrace{\gammud}_{\mu_1}\quad
\underbrace{\gammut}_{\mu_2}\quad \cdots\quad 
\underbrace{\gammuq}_{\mu_\ell}\quad .$$ 
Then {\bf the character $\chi^\lambda(\mu)$ is given by
$$
\chi^\lambda(\mu) = \sum_{T} \wt^\mu(T), $$
where the sum is over all standard tableaux $T$ of
shape $\lambda$ and
$$
\wt^\mu(T) = \prod_{i=1}^n f(i,T),$$
where
$$
f(i,T) = \cases{
-1, &if $i\not\in B(\mu)$ and $i+1$ is sw of $i$,\cr
0, &if $i,i+1\not\in B(\mu)$,
$i+1$ is ne of $i$, and
$i+2$ is sw of $i+1$,\cr
q, &otherwise,\cr}
$$
and
$B(\mu)= \{ \mu_1+\mu_2+\cdots+\mu_k | 1\le k\le \ell \}
$.
}
In the formula for $f(i,T)$, {\bf sw} means strictly south and weakly west
and {\bf ne} means strictly north and weakly east.
\end{enumerate}
\end{enumerate}

\bigskip\noindent
{\bf References}
\smallskip\noindent
The book [CP] contains a treatment of the Schur-Weyl duality type
theorem given above.  See also the references there.  Several
basic results on the Iwahori-Hecke algebra are given in the 
book [GHJ].  The theorem giving the explicit values of $q$ such 
that $H_k(q)$ is semisimple is due to Gyoja and Uno [GU].
The character formula given above is due to Roichman [Ro].
See [Ra3] for an elementary proof.

\bppsection{The Brauer algebras $B_k(x)$}

\bigskip
Fix $x\in \C.$
A {\it Brauer diagram on $k$ dots\/}
is a graph on two rows of $k$-vertices, one above
the other, and $k$ edges such that each vertex is incident to precisely
one edge.
The product of two $k$-diagrams
$d_1$ and $d_2$ is obtained by placing $d_1$ above $d_2$ and
identifying the
vertices in the bottom row of $d_1$ with the corresponding vertices
in the top row of $d_2$.   The resulting graph contains
$k$ paths and some number
$c$ of closed loops.   If $d$ is the $k$-diagram
with the edges that are the
paths in this graph but with the closed loops removed, then the product $d_1d_2$
is given by $d_1d_2=\eta^c d$.  For example, if
$$d_1=\ \ \donediag
\qquad\hbox{and}\qquad
d_2=\ \ \dtwodiag~~,$$
then
$$d_1d_2 = \dcomp = x^2\ \ \dprod\ \ .$$

The {\it Brauer algebra} $B_k(x)$ is the span
of the $k$-diagrams with multiplication given by
the linear extension of the diagram multiplication.
The dimension of the Brauer algebra is
$$\dim(B_k(x)) = (2k)!! = (2k-1)(2k-3)\cdots 3 \cdot 1,$$
since the number of $k$-diagrams is $(2k)!!$.

The diagrams in $B_k(x)$ which have all  their edges
connecting  top vertices to
bottom vertices form a symmetric group $S_k$.
The elements
$$s_i ~=~ \sidiag \quad\hbox{and}\quad e_i ~=~ \eidiag\ ,$$
$1\le i\le k-1$,
generate the Brauer algebra $B_k(x)$.

\bppthm  The Brauer algebra $B_k(x)$  has a
presentation as an algebra by generators
$s_1,s_2,\ldots,s_{k-1}$, $e_1,e_2,\ldots, e_{k-1}$ and relations
$$s_i^2=1, \quad e_i^2=x e_i, \quad e_i s_i=s_i e_i=e_i,
\qquad 1\le i\le k-1,$$
$$s_i s_j=s_j s_i,\quad s_i e_j=e_j s_i,\quad e_i e_j=e_j e_i,
\qquad |i-j|>1,$$
$$s_i s_{i+1}s_i = s_{i+1}s_i s_{i+1}, \quad
e_i e_{i+1} e_i = e_i, \quad
e_{i+1}e_i e_{i+1} = e_{i+1},\qquad 1\le i\le k-2,$$
$$s_ie_{i+1}e_i = s_{i+1}e_i,\quad
e_{i+1}e_i s_{i+1} = e_{i+1}s_i,\qquad 1\le i\le k-2.$$
\endthm

There are two different Brauer algebra analogues of the 
Schur Weyl duality theorem, Theorem A5.1.
In the first one the orthogonal group $O(n,\C)$ plays
the same role that $GL(n,\C)$ played in the $S_k$-case,
and in the second, the symplectic group $Sp(2n,\C)$
takes the $GL(n,\C)$ role.
\smallskip\noindent
Let $O(n,\C) = \{ A\in M_n(\C) \ |\ AA^t=I \}$ be the orthogonal
group and let $V$ be the usual $n$-dimensional representation
of the group $O(n,\C)$.  There is an action of the Brauer algebra
$B_k(n)$ on $V^{\otimes k}$ which commutes with the action
of $O(n,\C)$ on $V^{\otimes k}$.

\bppthm
\begin{enumerate}
\item[(a)] The action of $B_k(n)$ on $V^{\otimes k}$ generates
$\End_{O(n)}(V^{\otimes k}).$
\item[(b)] The action of $O(n,\C)$ on $V^{\otimes k}$ generates
$\End_{B_k(n)}(V^{\otimes k}).$
\end{enumerate}
\endthm

\noindent
Let $Sp(2n,\C)$ be the 
symplectic group and let $V$ be the usual $2n$-dimensional representation
of the group $Sp(2n,\C)$.  There is an action of the 
Brauer algebra
$B_k(-2n)$ on $V^{\otimes k}$ which commutes with the action
of $Sp(2n,\C)$ on $V^{\otimes k}$.

\bppthm
\begin{enumerate}
\item[(a)] The action of $B_k(-2n)$ on $V^{\otimes k}$ generates
$\End_{Sp(2n,\C)}(V^{\otimes k}).$
\item[(b)] The action of $Sp(2n,\C)$ on $V^{\otimes k}$ generates
$\End_{B_k(-2n)}(V^{\otimes k}).$
\end{enumerate}
\endthm

\bppthm The Brauer algebra $B_k(x)$ is semisimple if 
$x \not\in  \{-2k+3, -2k+2, \ldots, k-2\}$.
\endthm

\vskip-.1in
\subsection*{Partial results for $B_k(x)$}

\noindent
The following results giving answers to the main questions
{\bf (Ia-c)} for the Brauer algebras 
hold when $x$ is such that $B_k(x)$ is semisimple.

\bigskip\noindent
{\bf I.  What are the irreducible $B_k(x)$-modules?}
\medskip\noindent
\begin{enumerate}
\item[(a)] How do we index/count  them?
\smallskip
\begin{enumerate}
{\bf 
     \item[]
There is a bijection
$$
\hbox{Partitions of $k-2h$, $h=0,1,\ldots, \lfloor k/2\rfloor$}
\qquad
\mapleftright{1-1}
\qquad
\hbox{Irreducible representations $B^\lambda$.}
$$
}
\end{enumerate}

\item[(b)] What are their dimensions?
\smallskip\noindent
\begin{enumerate}
     \item[]
The dimension of the irreducible representation
$B^\lambda$ is given by
$$\eqalign{
{\bf dim}(B^\lambda) &=
\hbox{{\bf \# of up-down tableaux of shape $\lambda$ and length $k$}} \cr
&= 
{k\choose 2h} (2h-1)!! {(k-2h)!\over \prod_{x\in \lambda} h_x},
\cr}$$
where $h_x$ is the hook length at the box $x$ in $\lambda$.
An up-down tableau of shape $\lambda$ and length $k$
is a sequence $(\emptyset =\lambda^{(0)}, \lambda^{(1)},
\cdots \lambda^{(k)}=\lambda)$ of partitions,
such that each partition in the sequence differs
from the previous one by either addingor removing
a box.
\end{enumerate}
\item[(c)] What are their characters?
\smallskip\noindent
\begin{enumerate}
     \item[]  
A Murnaghan-Nakayama type rule for the characters of the Brauer
algebras was given in [Ra1].
\end{enumerate}
\end{enumerate}

\medskip\noindent
{\bf References}
\smallskip\noindent
{\bf (1)} The Brauer algebra was defined originally by R. Brauer [Br]
in 1937.  H. Weyl treats it in his book [Wy].
\smallskip\noindent
{\bf (2)} The Schur-Weyl duality type theorems are due to
Brauer [Br], from his original paper.
See also [Ra1] for a detailed description of these Brauer 
algebra actions.
\smallskip\noindent
{\bf (3)} 
The theorem giving values of $x$ for which
the Brauer algebra is semisimple is due to Wenzl, see [Wz2].

\bppsection{The Birman-Murakami-Wenzl algebras $BMW_k(r,q)$}

\bigskip
A {\it $k$-tangle\/} is viewed  as two rows of $k$ vertices,
one above the other, and $k$ strands that connect vertices in such a way that
each vertex is incident to precisely one strand.  Strands cross over and under
each other in three-space as they pass from one vertex to the next.
For example, the following are $7$-tangles:
$$
t_1 ~=~ \stanglesa \qquad t_2 ~=~ \stanglesb\  .
$$
We multiply $k$-tangles $t_1$ and $t_2$ using the concatenation product
given by identifying the vertices in the top row of 
$t_2$ with the corresponding
vertices in the bottom row of $t_1$ to obtain
the product tangle $t_1 t_2$.  Then we allow the following ``moves.''
\bigskip
{\it Reidemeister moves II and III}:
\medskip
\begin{enumerate}
     \item[(R2)]  \qquad$\displaystyle{
\rtwoa ~~\longleftrightarrow~~ \rtwob ~~\longleftrightarrow~~ \rtwoc}$
     \item[(R3)] \qquad$\displaystyle{
\rthreea ~~\longleftrightarrow~~ \rthreeb}$
\end{enumerate}

\bigskip
Given a Brauer diagram $d$ we will make a tangle $T_d$
tracing the edges in order from left to right across the top row and
then from left to right across the bottom row.  Any time
an edge that we are tracing crosses and edge that has been already traced
we raise the pen briefly so that the edge being traced goes under 
the edge which is already there.  Applying this process 
to all of the Brauer diagrams on $k$ dots produces a set 
of $(2k)!!$ tangles.  

$$
d~=~~ \brauertang\qquad\qquad T_d ~=~~\stanglesa ~
$$

Fix $r,q\in \C$.  The {\it Birman-Murakami-Wenzl algebra\/}
$BMW_k(r,q)$
is the span of the $(2k)!!$ tangles produced by tracing the Brauer
diagrams with multiplication determined
by the tangle multiplication and the Reidemeister moves and
the following tangle identities.
\bigskip\noindent
$$ 
\tangleida  ~~=~  (q-q^{-1}) \left( ~~\tangleidb~~ \right). 
%\leqno(T1)
$$
\smallskip
$$
\Tda ~~=~  r^{-1} ~~\Tdb~ , \qquad\qquad
\Tdaa ~~=~ r ~~\Tdb\ . 
$$
\smallskip
$$
\Tt ~~=~  ~x,\qquad \hbox{where}\qquad
x = {r-r^{-1}\over q-q^{-1}} + 1.
%\leqno(T3)
$$
\smallskip
The Reidemeister moves and the tangle identities
can be applied in any appropriate local portion of the tangle.

\bppthm  Fix $r,q\in \C$.
The {\it Birman-Murakami-Wenzl algebra} $BMW_k(r,q)$ is the 
algebra generated over $\C$ by $1, g_1, g_2, \ldots, g_{k-1}$,
which are assumed to be invertible,
subject to the relations
\medskip
\begin{enumerate}
     \item[] $\quad g_i g_{i+1} g_i = g_{i+1} g_i g_{i+1},$
     \item[] $\quad g_i g_j = g_j g_i \qquad\hbox{if  } |i-j| \ge 2,$
     \item[]$ \quad (g_i - r^{-1})(g_i + q^{-1})(g_i - q) = 0,$
     \item[]$ \quad E_i g_{i-1}^{\pm 1} E_i ~=~ r^{\pm 1} E_i$ \quad and
\quad $E_i g_{i+1}^{\pm 1} E_i ~=~ r^{\pm 1} E_i,$ 
\end{enumerate}
where $E_i$ is defined by the equation
$$
(q-q^{-1})(1 - E_i) = g_i - g_i^{-1}.
$$
\endthm

The BMW-algebra is a $q$-analogue of the Brauer algebra in the
same sense that the Iwahori-Hecke algebra of type $A$ is a 
$q$-analogue of the group algebra of the symmetric group.
If we allow ourselves to be imprecise (about the limit) we can write
$$\lim_{q\to 1} BMW_k(q^{n+1},q) = B_k(n).$$
It would be interesting to sharpen the following theorem to
make it an if and only if statement.

\bppthm [Wz3] The Birman-Murakami-Wenzl algebra is semisimple if
$q$ is not a root of unity and $r \ne q^{n+1}$ for any $n \in \Z$. 
\endthm

\subsection*{Partial results for $BMW_k(r,q)$}

The following results hold when $r$ and $q$ are such that 
$BMW_k(r,q)$ is semisimple.

\bigskip\noindent
{\bf I.  What are the irreducible $BMW_k(r,q)$-modules?}
\medskip\noindent
\begin{enumerate}
\item[(a)] How do we index/count  them?
\smallskip
\begin{enumerate}
{\bf
     \item[]
There is a bijection
$$
\hbox{Partitions of $k-2h$, $h=0,1,\ldots, \lfloor k/2\rfloor$}
\qquad
\mapleftright{1-1}
\qquad
\hbox{Irreducible representations $W^\lambda$.}
$$
}
\end{enumerate}
\item[(b)] What are their dimensions?
\smallskip\noindent
\begin{enumerate}
     \item[]
The dimension of the irreducible representation
$W^\lambda$ is given by
$$\eqalign{
{\bf dim}(W^\lambda) &=
\hbox{{\bf \# of up-down tableaux of shape $\lambda$ and length $k$}} \cr
&= 
{k\choose 2h} (2h-1)!! {(k-2h)!\over \prod_{x\in \lambda} h_x},
\cr}$$
where $h_x$ is the hook length at the box $x$ in $\lambda$,
and up-down tableaux is as in the case of the Brauer algebra,
see Section B6 {\bf (Ib)}.
\end{enumerate}

\item[(c)] What are their characters?
\smallskip\noindent
\begin{enumerate}
     \item[]  
A Murnaghan-Nakayama rule for the irreducible characters of
the BMW-algebras was given in [HR1].
\end{enumerate}
\end{enumerate}

\bigskip\noindent
{\bf References}
\smallskip\noindent
{\bf (1)}  The Birman-Murakami-Wenzl algebra was defined independently
by Birman and Wenzl in [BW] and by Murakami in [Mu1].
See [CP] for references to the analogue of Schur-Weyl duality for
the BMW-algebras.  The articles [HR1], [LR], [Mu2], [Re], and [Wz3]  
contain further important information about the BMW-algebras.
\smallskip\noindent
{\bf (2)} Although the tangle description of the BMW algebra
was always in everybody's minds it was Kaufmann that really made
it precise see [Ka2].

\bppsection{The Temperley-Lieb algebras $TL_k(x)$}

\bigskip
A {\it $TL_k$-diagram} is a Brauer diagram on $k$ dots
which can be drawn with no crossings of edges.
$$
\TLkdiagram
$$
The {\it Temperley-Lieb algebra} $TL_k(x)$ is the
subalgebra of the Brauer algebra $B_k(x)$ which is
the span of the $TL_k$-diagrams.

\bppthm 
The Temperley-Lieb algebra $TL_k(x)$ 
is the algebra over $\C$ given by 
generators $E_1,E_2,\ldots,E_{k-1}$ and relations
\medskip
\begin{enumerate}
     \item[] \quad $E_iE_j ~=~ E_jE_i,$\quad if $|i-j|>1$,
     \item[] \quad $E_iE_{i\pm1}E_i ~=~ E_i,$\quad and
     \item[] \quad $E_i^2~=~\displaystyle{x}E_i$ .
\end{enumerate}
\endthm

\bppthm Let $q\in \C^*$ be such that $q+q^{-1}+2 = 1/x^2$
and let $H_k(q)$ be the Iwahori-Hecke algebra of type $A_{k-1}$.
Then the map
$$\matrix{
H_k(q) &\longrightarrow &TL_k(x) \cr
\cr
T_i &\longrightmapsto &\displaystyle{{q+1\over x}E_i-1} \cr
}$$
is a surjective homomorphism and the kernel of this 
homomorphism is the ideal generated by the elements
$$
T_iT_{i+1}T_i+T_iT_{i+1}+T_{i+1}T_i+ T_i+T_{i+1}+1 , 
\quad\hbox{for $1\le i\le n-2$.}
$$  
\endthm

The Schur Weyl duality theorem for $S_n$ has the following
analogue for the Temperley-Lieb algebras.
Let $U_q{\goth{sl}}_2$ be the Drinfeld-Jimbo quantum group
corresponding to the Lie algebra ${\goth{sl}}_2$
and let $V$ be the $2$-dimensional representation
of $U_q{\goth{sl}}_2$.  There is an action, see [CP], of the 
Temperley-Lieb algebra
$TL_k(q+q^{-1})$ on $V^{\otimes k}$ which commutes with the action
of $U_q{\goth{sl}}_2$ on $V^{\otimes k}$.

\bppthm
\begin{enumerate}
\item[(a)] The action of $TL_k(q+q^{-1})$ on $V^{\otimes k}$ generates
$\End_{U_q{\goth{sl}}_2}(V^{\otimes k}).$
\item[(b)] The action of $U_q{\goth{sl}}_2$ on $V^{\otimes k}$ generates
$\End_{TL_k(q+q^{-1})}(V^{\otimes k}).$
\end{enumerate}
\endthm

\bppthm  The Temperley-Lieb algebra is semisimple if and only if
$1/x^2 \ne 4\cos^2(\pi/\ell)$, for any $2\le \ell\le k$.
\endthm

\vskip-.1in
\subsection*{Partial results for $TL_k(x)$}

\noindent
The following results giving answers to the main questions
{\bf (Ia-c)} for the Temperley-Lieb algebras 
hold when $x$ is such that $TL_k(x)$ is semisimple.

\bigskip\noindent
{\bf I.  What are the irreducible $TL_k(x)$-modules?}
\medskip\noindent
\begin{enumerate}
\item[(a)] How do we index/count  them?
\smallskip
\begin{enumerate}
{\bf 
     \item[]
There is a bijection
$$
\hbox{ Partitions of $k$ with at most two rows}
\qquad
\mapleftright{1-1}
\qquad
\hbox{Irreducible representations $T^\lambda$.}
$$
}
\end{enumerate}

\item[(b)] What are their dimensions?
\smallskip\noindent
\begin{enumerate}
     \item[]
The dimension of the irreducible representation
$T^{(k-\ell,\ell)}$ is given by
$$\eqalign{
{\bf dim}(T^{(k-\ell,\ell)}) &=
\hbox{{\bf \# of standard tableaux of shape $(k-\ell,\ell)$}} \cr
&= {k \choose \ell} - {k \choose \ell-1}.\cr
}$$  
\end{enumerate}

\item[(c)] What are their characters?
\smallskip\noindent
\begin{enumerate}
     \item[]  
The {\bf character of the irreducible
representation $T^{(k-\ell,\ell)}$ evaluated at the element 
$$d_{2h} ~=~ 
\underbrace{\idkdhAa\quad\cdots\quad \idkdhAb}_{k-2h}
\quad\underbrace{\idkdhBa\quad\cdots\quad \idkdhBb}_{2h}$$
is
$$
\chi^{(k - \ell,\ell)}(d_{2h})
~=~ \cases{
\displaystyle{k - 2h \choose \ell - h} - 
\displaystyle{k - 2h \choose \ell - h - 1}, 
& if $ \ell \ge h,$ \cr
\cr
0, & if $\ell < h.$ \cr}
$$
}
There is an algorithm for writing the character
$\chi^{(k - \ell,\ell)}(a)$ of a general element
$a\in TL_k(x)$ as a linear combination of the characters 
$\chi^{(k - \ell,\ell)}(d_{2h})$.
\end{enumerate}
\end{enumerate}

\bigskip\noindent
{\bf References}
\smallskip\noindent
The book [GHJ] contains a comprehensive treatment of the 
basic results on the Temperley-Lieb algebra.  
The Schur-Weyl duality theorem is treated in the book [CP],
see also the references there.
The character formula given above
is derived in [HR1].

\bppsection{Complex semisimple Lie groups}

\bigskip
We shall not define Lie groups and Lie algebras let us only
recall that a complex Lie group is a differential $\C$-manifold
and a real Lie group is a differential $\R$-manifold and
that every Lie group has an associated Lie algebra,
see [CMS].

\smallskip
If $G$ is a complex Lie group then the word {\it representation}
is usually used to refer to a {\it holomorphic representation}, i.e.
the homomorphism
$$\rho\colon G \to GL(V)$$
determined by the module $V$ should be a morphism of (complex) analytic
manifolds.
Strictly speaking there are representations which are not holomorphic
but there is a good theory only for holomorphic representations, so
one usually abuses language and assumes that representation means
holomorphic representation.  The terms holomorphic
representation and {\it complex analytic representation}
are used interchangeably.
Similarly, if $G$ is a real Lie group then representation usually
means {\it real analytic representation}.  See [Va] p. 102 
for further details.
Every holomorphic representation of $GL(n,\C)$ is also
rational representation, see [FH].

\smallskip
A {\it complex semisimple Lie group} is a connected
complex Lie group
$G$ such that its Lie algebra ${\goth g}$ is a complex semisimple
Lie algebra.


\begin{thebibliography}{}

\bibitem[Al]{Al}  E. Allen, {\it New bases for the decompositions of
the regular representations of wreath products of $S_n$},
preprint 1997.

\bibitem[Ari]{Ari}  S. Ariki, 
{\it Representation theory of a Hecke algebra of $G(r,p,n)$},
J. of Algebra {\bf 177} (1995), 164-185.

\bibitem[AK]{AK}  S. Ariki and K. Koike, {\it  A Hecke algebra of $(\Z/r\Z) \wr S_n$
and construction of its irreducible representations},
Adv.~in Math. {\bf  106} (1994), 216-243.

\bibitem[Ba]{Ba}  H. Barcelo, {\it  Young straightening in a quotient $S_n$-module},
J. Alg. Comb. {\bf  2} (1993), 5-23.

\bibitem[Ba2]{Ba2} H. Barcelo, 
{\it On the action of the symmetric group on the free Lie algebra and
the partition lattice}, J. Combin. Theory Ser. A {\bf 55} 
(1990), 93--129. 

\bibitem[BBL]{BBL} G. Benkart, D. Britten, F. Lemire,
{\it Stability in modules for classical Lie algebras--a constructive
approach}, Mem. Amer. Math. Soc. {\bf 85} (1990).

\bibitem[BC]{BC} G. Benkart, M. Chakrabarti, T. Halverson, R. Leduc,
C. Lee, J. Stroomer, {\it Tensor product representations of general
linear groups and their connections with Brauer algebras},
J. Algebra {\bf 166} (1994), 529-567.

\bibitem[Be1]{Be1} A. Berele, {\it A Schensted type correspondence for the 
symplectic group}, J. Combin. Theory Ser. A {\bf 43} (1986), 320-328.

\bibitem[Be2]{Be2} A. Berele, {\it Construction of $Sp$ modules by tableaux},
Linear and Multilinear algebra {\bf 19} (1986), 299-307.

\bibitem[BGG]{BGG} Bernstein, Gelfan'd, and S. Gelfan'd,
{\it  Schubert cells and the cohomology of the spaces $G/P$},
Russian Math. Surveys {\bf 28:3} (1973), 1-26.

\bibitem[BM]{BM}   M. Brou\'e and G. Malle,
{\it Zyklotomische Heckealgebren},   Ast\'erisque
{\bf 212} (1993),  119-189.

\bibitem[Bo]{Bo} A. Borel, {\it Sur la cohomologie des
espaces fibr\'es principaux et des espaces homog\`enes de 
groupes de Lie compacts}, Ann. Math. (2) {\bf 57} (1953),
115-207.

\bibitem[Bou1]{Bou1} N. Bourbaki, 
{\it Groupes et Alg\`ebre de Lie, Chapitre IV, V, VI},
El\'ements de Math\'ematique, Hermann, Paris, (1968)

\bibitem[Bou2]{Bou2} N. Bourbaki, {\it Alg\`ebre, Chapt. VIII},
Hermann, Paris 1958.

\bibitem[BR]{BR} A. Berele and A. Regev, {\it Hook Young diagrams with
applications to combinatorics and to representations of Lie superalgebras},
Adv. in Math. {\bf 64} (1987), 118-175.

\bibitem[Br]{Br}   R. Brauer,
{\it On algebras which are connected with the semisimple continous groups},  
Ann. Math. {\bf 38} (1937),  854-872.

\bibitem[Bre]{Bre} F. Brenti, 
{\it A combinatorial formula for Kazhdan-Lusztig polynomials},
Invent. Math. {\bf 118} (1994), 371--394.

\bibitem[BSR]{BSR} G. Benkart, C. Lee Shader, and A. Ram,
{\it Tensor representations of orthosymplectic Lie color algebras},
to appear in J. Pure and Applied Algebra.

\bibitem[BtD]{BtD} T. Br\"ocker and T. tom Dieck, {\it Representations
of compact Lie groups}, Graduate Texts in Mathematics {\bf 98},
Springer-Verlag, New York-Berlin, 1985.

\bibitem[BW]{BW} J. Birman and H. Wenzl,
{\it Braids, link polynomials and a new algebra},
Trans. Amer. Math. Soc. {\bf 313} (1989), 249-273.

\bibitem[Ca]{Ca}   R. W. Carter, 
{\it Finite Groups of Lie Type: Conjugacy Classes and Complex Characters},
Wiley, 1985.

\bibitem[CG]{CG} N. Chriss and V. Ginzburg, {\it Representation
theory and complex geometry}, Birkhauser, Boston 1997.

\bibitem[Cl]{Cl}   A. H. Clifford, 
{\it Representations induced in an invariant subgroup}, 
Ann. of Math. {\bf 38} (1937), 533-550.

\bibitem[CMS]{CMS} R.W. Carter, I.G. Macdonald, G. Segal,
{\it Lectures on Lie groups and Lie algebras}, 
London Mathematical Society Student Texts {\bf 32},
Cambridge University Press, Cambridge, 1995.

\bibitem[CP]{CP} V. Chari and A. Pressley, {\it A guide to quantum
groups}, Cambridge University Press, Cambridge, 1994.

\bibitem[CR1]{CR1}   C. W. Curtis, and I. Reiner,
{\it Representation Theory of Finite groups and Associative Algebras},
Wiley, 1988.

\bibitem[CR2]{CR2}   C. W. Curtis, and I. Reiner,
{\it Methods of Representation Theory with Applications to Finite Groups
and Orders}, Vol. I, Vol II, Wiley, 1987.

\bibitem[Cu]{Cu} C.W. Curtis,
{\it Representations of Hecke algebras}, Orbites unipotentes et
représentations, I. Astérisque {\bf 168} (1988), 13--60.

\bibitem[DJ1]{DJ1}   R. Dipper and G. James,
{\it Representations of Hecke algebras of general linear groups},
Proc. London Math. Soc. (3) {\bf 52} (1986), no. 1, 20-52.

\bibitem[DJ2]{DJ2}   R. Dipper and G. James,
{\it Representations of Hecke algebras of type $B_n$}
J. Alg. {\bf 146} (1992), no. 2, 454-481.

\bibitem[DJM]{DJM} R. Dipper, G. James, and E. Murphy,
{\it Hecke algebras of type $B_n$ at roots of unity},
Proc. London Math. Soc. (3) {\bf 70} (1995), 505-528.

\bibitem[DK]{DK} G. Duchamp, D. Krob, A. Lascoux, B. Leclerc, 
T. Scharf, and J.-Y. Thibon, {\it Euler-Poincar\'e 
characteristic and polynomial representations of 
Iwahori-Hecke algebras}, Publ. Res. Inst. Math. Sci.
{\bf 31} (1995), 179-201.

\bibitem[DL]{DL} 
P. Deligne and G. Lusztig, 
{\it Representations of reductive groups over finite fields},
Ann. of Math. (2) {\bf 103} (1976), 103--161.

\bibitem[DM]{DM} F. Digne and J. Michel,  
{\it Representations of finite groups of Lie type}, 
London Mathematical Society Student Texts {\bf 21}, 
Cambridge University Press, Cambridge, 1991.

\bibitem[DP]{DP} C. DeConcini and C. Procesi,
{\it Symmetric functions, conjugacy classes and the flag variety},
Invent. Math. {\bf 64} (1981), 203-219.

\bibitem[FF]{FF} 
B.L. Feigin and  D.B. Fuchs, 
{\it Representations of the Virasoro algebra},
Representation of Lie groups and related topics,  
Adv. Stud. Contemp. Math. {\bf 7}, Gordon and Breach, New York, 1990,
465-554.

\bibitem[FG]{FG} S. Fomin and C. Greene, {\it Noncommutative Schur functions
and their applications}, to appear in Discrete Math. J. (Special
volume in Honour of A. Garsia). 

\bibitem[FH]{FH} W. Fulton and J. Harris,
{\it Representation theory, A first course},
Graduate Texts in Mathematics {\bf 129}, Springer-Verlag, New York, 1991.

\bibitem[Fr]{Fr}   F.G. Frobenius, 
{\it \"Uber die Charaktere der symmetrischen Gruppe},
Sitzungsberichte der K\"oniglich Preussischen Akademie der
Wissenschaften zu Berlin (1900), 516-534 (Reprinted in Gesammelte 
Abhandlungen, {\bf 3}, 
148-166).

\bibitem[FRT]{FRT}   G.S. Frame, G. de B. Robinson, and R.M. Thrall,
{\it The hook graphs of $S_n$}, Canadian J. Math. {\bf 6} (1954), 316-324.

\bibitem[Gar]{Gar} A. Garsia, {\it 
Combinatorics of the free Lie algebra and the symmetric
group}, Analysis, et cetera, Academic Press, Boston, MA, 1990,
309-382.

\bibitem[GHJ]{GHJ}  F. Goodman, P. de la Harpe, and V.F.R. Jones,
{\it Coxeter graphs and towers of algebras}, 
Mathematical Sciences Research Institute Publications
{\bf 14}, Springer-Verlag, New York-Berlin, 1989.

\bibitem[Gi]{Gi} V. Ginzburg, {\it ``Lagrangian'' construction for
representations of Hecke algebras},
Adv. in Math. {\bf 63} (1987), 100-112.

\bibitem[GL]{GL} J. Graham and G. Lehrer, {\it Cellular algebras},
Invent. Math. {\bf 123} (1996), 1-34.

\bibitem[GM]{GM} A. M. Garsia and T. J. Maclarnan, {\it Relations
between Young's natural and the Kazhdan-Lusztig representations of $S_n$},
Adv. in Math. {\bf 69} (1988), no. 1, 32-92.

\bibitem[GR]{GR} A. Garsia and J. Remmel, {\it $q$-counting rook
configurations and a formula of Frobenius}, J. Combin. Theory
Ser. A {\bf 41} (1986), 246-275.

\bibitem[Gr]{Gr} J. A. Green, {\it  The characters of the finite general
linear groups}, Trans. Amer. Math. Soc. {\bf 80} (1955), 402-447.

\bibitem[GT1]{GT1} I.M. Gelfand and M.L. Tsetlin,
{\it Finite-dimensional representations of the group of
unimodular matrices}, Doklady Akad. Nauk SSSR (N.S.) {\bf 71} (1950), 
825--828.

\bibitem[GT2]{GT2} I.M. Gelfand and M.L. Tsetlin,
{\it Finite-dimensional representations of groups
of orthogonal matrices}, Doklady Akad. Nauk SSSR 
(N.S.) {\bf 71} (1950), 1017--1020.

\bibitem[GU]{GU} A. Gyoja and K. Uno,
{\it On the semisimplicity of Hecke algebras}, J. Math.
Soc. Japan {\bf 41} (1989), 75-79.

\bibitem[GW]{GW} A. Garsia and M. Wachs, 
{\it Combinatorial aspects of skew representations of
the symmetric group}, J. Combin. Theory Ser. A {\bf 50} (1989), 47--81.

\bibitem[Gy]{Gy} A. Gyoja, 
{\it A $q$-analogue of Young symmetrizer},
Osaka J. Math. {\bf 23} (1986), no. 4, 841-852.

\bibitem[Ha1]{Ha1} T. Halverson, {\it Characters of the centralizer algebras
of mixed tensor representations of $GL(r,\C)$ and the
quantum group $U_q{\goth{gl}}(r,\C)$},
Pacific J. Math. {\bf 174} (1996), 359-410.

\bibitem[Ha2]{Ha2} T. Halverson, {\it A $q$-rational Murnaghan-Nakayama
rule}, J. Combin. Theory Ser. A {\bf 71} (1995), 1-18.

\bibitem[Hfs]{Hfs}   P. N. Hoefsmit,
Representations of Hecke algebras of finite groups with BN-pairs of
classical type, Thesis, University of British Columbia,  1974.

\bibitem[Hn]{Hn} P. Hanlon, {\it The fixed point partition lattices},
Pacific J. Math. {\bf 96} (1981), 319-341.

\bibitem[Ho]{Ho} R. Howe, {\it Hecke algebras and $p$-adic ${\rm GL}\sb n$},
Representation theory and analysis on homogeneous spaces 
(New Brunswick, NJ, 1993),  Contemp.  Math. {\bf 177}, 
Amer. Math. Soc., Providence, RI, 1994, 65-100. 

\bibitem[HR1]{HR1}  T. Halverson and A. Ram,
{\it Characters of algebras containing a Jones basic construction:
The Temperley-Lieb, Okada, Brauer, and Birman-Wenzl algebras},
Adv. Math. {\bf 116} (1995), 263--321.

\bibitem[HR2]{HR2}   T. Halverson and A. Ram,
{\it Murnaghan-Nakayama rules for the characters of Iwahori-Hecke
algebras of the complex reflection groups $G(r,p,n)$},
to appear in Canadian ~J. ~Math.

\bibitem[HS]{HS} R. Hotta and T. Springer,
{\it A specialization theorem for certain Weyl group
representations and an application to Green polynomials of
unitary groups}, Invent. Math. {\bf 41} (1977), 113-127.

\bibitem[Hu1]{Hu1} J.E. Humphreys,
{\it Introduction to Lie algebras and representation theory},
Graduate Texts in Mathematics {\bf 9}, Springer-Verlag,
New York-Berlin,  1978.

\bibitem[Hu2]{Hu2}   J. E. Humphreys,
{\it Reflection Groups and Coxeter Groups}, Cambridge Studies in
Advanced Mathematics No. 29, Cambridge University Press, 1990.

\bibitem[HW]{HW} P. Hanlon and D. Wales,
{\it On the decomposition of Brauer's centralizer algebras},
J. Algebra {\bf 121} (1989), 409-445.

\bibitem[Iw]{Iw} N. Iwahori, {\it On the structure of a Hecke ring
of a Chevalley group over a finite field}, J. Fac. Sci. Univ. Tokyo
Sect. 1, {\bf 10} (1964), 215-236.

\bibitem[Ja]{Ja} J. Jantzen, {\it Lectures on quantum groups},
Graduate Studies in Mathematics {\bf 6}, Amer. Math. Soc.,
Providence, RI 1996.

\bibitem[Ji]{Ji}   M. Jimbo,
{\it A $q$-analog of $U(gl(N+1))$, Hecke algebra, and the Yang-Baxter
equation}, Lett. Math. Phys. {\bf II} (1986), 247-252.

\bibitem[Jo1]{Jo1} V. Jones, {\it The Potts model and the symmetric group},
in {\sl Subfactors} (Kyuzeso, 1993), World Sci. Publishing, River
Edge, NJ (1994) 259-267.

\bibitem[Jo2]{Jo2} V. Jones, {\it Index for subfactors},
Invent. Math. {\bf 72} (1983), 1-25.

\bibitem[Jy]{Jy} A. Joyal, {\it Foncteurs analytiques et esp\`eces de
structures}, Lect. Notes in Math. {\bf 1234}, Springer-Verlag (1986),
126-160.

\bibitem[Ka1]{Ka1} L. Kauffman, {\it State models and the Jones polynomial},
Topology {\bf 26} (1987), 395-407.

\bibitem[Ka2]{Ka2} L. Kauffman, {\it An invariant of regular isotopy},
Trans. Amer. Math. Soc. {\bf 318} (1990), 417-471.

\bibitem[Kat]{Kat} S. Kato, 
{\it A realization of irreducible representations of affine Weyl groups},
Nederl. Akad. Wetensch. Indag. Math. {\bf 45} (1983), 193--201. 

\bibitem[Kc]{Kc} V. Kac, {\it Infinite dimensional Lie algebras},
Third edition, Cambridge University Press, Cambridge 1990.

\bibitem[Ke]{Ke} S.V. Kerov, {\it Characters of Hecke and Birman-Wenzl
algebras}, in Quantum groups (Leningrad 1990), Lecture Notes in Math.
{\bf 1510}, Springer, Berlin (1992), 179-187.

\bibitem[KK]{KK} B. Kostant and S. Kumar, {\it The nil Hecke ring and
cohomology of $G/P$ for a Kac-Moody group $G$}, Adv. in Math. {\bf 62}
(1986), 187-237.

\bibitem[Kk]{Kk} C. Kostka, {\it \"Uber den Zusammenhang zwischen einigen
Formen von symmetrischen Funktionen}, Crelle's Journal {\bf 93} (1882),
89-123.

\bibitem[KL1]{KL1} D. Kazhdan and G. Lusztig, {\it Representations of
Coxeter groups and Hecke algebras}, Invent. Math. {\bf 53} (1979),
165-184.

\bibitem[KL2]{KL2} D. Kazhdan and G. Lusztig, {\it Proof of the 
Deligne-Langlands conjecture for Hecke algebras}, Invent. Math.
{\bf 87} (1987), 153-215.

\bibitem[Kl]{Kl} A.A. Klyachko, {\it Lie elements in the tensor algebra},
Siberian Math. J. {\bf 15} (1974), 1296-1304.

\bibitem[KM1]{KM1} M. Kosuda and J. Murakami, {\it Centralizer algebras of the 
mixed tensor representations of quantum group $U_q{\goth{gl}}(n,\C)$},
Osaka J. Math. {\bf 30} (1993), 475-507.

\bibitem[KM2]{KM2} M. Kosuda and J. Murakami, {\it The centralizer algebras
of mixed tensor representations of $U_q{\goth{gl}}_n$ and the 
HOMFLY polynomial of links}, Proc. Japan Acad. Ser. A Math. Sci.
{\bf 68} (1992), 148-151.

\bibitem[Ko]{Ko} K. Koike, {\it On the decomposition of tensor products
of the representations of the classical groups: by means of universal
characters}, Adv. Math. {\bf 74} (1989), 57-86.

\bibitem[Kr]{Kr} H. Kraft, {\it Conjugacy classes and Weyl group
representations}, Proc. 1980 Torun Conf. Poland, Asterisque
{\bf 87-88} (1981), 195-205.

\bibitem[Ksh]{Ksh} M. Kashiwara, 
{\it Crystalizing the $q$-analogue of universal enveloping
algebras}, Comm. Math. Phys. {\bf 133} (1990), 249--260.

\bibitem[Ku1]{Ku1} G. Kuperberg,
{\it Symmetries of plane partitions and the permanent-determinant
method}, J. Combin. Theory Ser. A {\bf 68} (1994), 115--151.

\bibitem[Ku2]{Ku2} G. Kuperberg,
{\it Self-complementary plane partitions by Proctor's minuscule
method}, European J. Combin. {\bf 15} (1994), 545--553.

\bibitem[Ku3]{Ku3} G. Kuperberg,
{\it Four symmetry classes of plane partitions under one roof},
J.  Combin. Theory Ser. A {\bf 75} (1996), 295--315.

\bibitem[Ku4]{Ku4} G. Kuperberg,
{\it Spiders for rank $2$ Lie algebras}, Comm. Math. Phys. {\bf 180}
(1996), 109--151.

\bibitem[Ku5]{Ku5} G. Kuperberg,
{\it The quantum $G\sb 2$ link invariant}, Internat. J. Math. 
{\bf 5} (1994), 61--85. 

\bibitem[KV]{KV}  A.W. Knapp and D. Vogan, 
{\it Cohomological induction and unitary
representations}, Princeton Mathematical Series {\bf 45}, 
Princeton University Press, Princeton, NJ, 1995.

\bibitem[KW]{KW} R.C. King and B.G. Wybourne,
{\it Representations and traces of the Hecke algebras $H_n(q)$
of type $A_{n-1}$}, J. Math. Phys. {\bf 33} (1992), 4-14.

\bibitem[KWe]{KWe} R.C. King and T.A. Welsh, 
{\it Construction of orthogonal group modules using tableaux},
Linear and Multilinear Algebra {\bf 33} (1993), 251--283.

\bibitem[L]{L} A. Lascoux,
{\it Cyclic permutations on words, tableaux and harmonic
polynomials}, Proceedings of the Hyderabad Conference on 
Algebraic Groups (Hyderabad, 1989), Manoj Prakashan, Madras, 1991, 323-347. 

\bibitem[Le]{Le} R. Leduc, Thesis, University of Wisconsin--Madison,
1994.

\bibitem[Leh]{Leh} G. Lehrer, 
{\it Poincar\'e polynomials for unitary reflection groups}, 
Invent. Math. {\bf 120} (1995), 411--425.

\bibitem[Li1]{Li1} P. Littelmann, {\it A Littlewood-Richardson rule for
symmetrizable Kac-Moody Lie algebras}, Invent. Math.
{\bf 116} (1994), 329-346.

\bibitem[Li2]{Li2} P. Littelmann, {\it Paths and root
operators in representation theory},
Ann. Math. (2) {\bf 142} (1995), 499-525.

\bibitem[LR]{LR} R. Leduc and A. Ram, {\it A ribbon Hopf algebra approach to
the irreducible representations of centralizer algebras:
The Brauer, Birman-Wenzl and type A Iwahori-Hecke algebras},
Adv. in Math. {\bf 125} (1997), 1-94.

\bibitem[LS]{LS} V. Lakshmibai and C.S. Seshadri,
{\it Standard monomial theory}, 
Proceedings of the Hyderabad Conference on Algebraic Groups 
(Hyderabad, 1989), Manoj Prakashan, Madras, (1991), 279-322. 

\bibitem[Lu1]{Lu1} G. Lusztig, {\it Cells in affine Weyl groups II},
J. Algebra {\bf 109} (1987), 536-548.

\bibitem[Lu2]{Lu2} G. Lusztig, {\it Cells in affine Weyl groups III},
J. Fac. Sci. Univ. Tokyo Sect. 1A Math. {\bf 34} (1987), 223-243.

\bibitem[Lu3]{Lu3} G. Lusztig, {\it Affine Hecke algebras and their graded
version}, J. Amer. Math. Soc. {\bf 2} (1989), 599-635.

\bibitem[Lu4]{Lu4} G. Lusztig, {\it Cuspidal local systems and graded
Hecke algebra}, Inst. Hautes Etudes Sci. Publ. Math. {\bf 67}
(1988), 145-202.

\bibitem[Lu5]{Lu5} G. Lusztig, {\it Representations of affine Hecke algebras,
Orbites unipotentes et repr\'esentations, II}, Ast\'erisque {\bf 171-172}
(1989), 73-84.

\bibitem[Lu6]{Lu6} G. Lusztig, {\it Cuspidal local systems and graded Hecke
algebras II}, CMS Conf. Proc. {\bf 16} Representations of groups
(Banff, AB, 1994), Amer. Math. Soc. Providence, RI (1995), 217-275.

\bibitem[Lu7]{Lu7} G. Lusztig, {\it Classification of unipotent representations
of simple $p$-adic groups}, Internat. Math. Res. Notices (1995),
517-589.

\bibitem[Lu8]{Lu8}
G. Lusztig, 
{\it Representations of finite Chevalley groups}, Expository lectures
from the CBMS Regional Conference held at Madison, Wis., August 8--12, 1977. 
CBMS Regional Conference Series in Mathematics {\bf 39},
American Mathematical Society, Providence, R.I., 1978.

\bibitem[Lu9]{Lu9}
G. Lusztig, {\it 
Classification des représentations irréductibles des groupes
classiques finis}, C. R. Acad. Sci. Paris Sér. A-B {\bf 284} 
(1977), A473--A475.

\bibitem[Lu10]{Lu10}
G. Lusztig, {\it Characters of reductive groups over a finite field}, 
Annals of Mathematics Studies {\bf 107}, 
Princeton University Press, Princeton, N.J. 1984.

\bibitem[Lu11]{Lu11}
G. Lusztig, {\it 
The discrete series of $GL\sb{n}$ over a finite field},
Annals of Mathematics Studies {\bf 81}, 
Princeton University Press, Princeton, N.J., 1974.

\bibitem[Lu12]{Lu12}
G. Lusztig and N. Spaltenstein,
{\it On the generalized Springer correspondence for 
calssical groups}, Algebraic groups and related topics
(kyoto/Nagoya, 1983), Adv. Stud. Pure Math. {\bf 6},
North-Holland, Amsterdam-New York, 1995.

\bibitem[Lw]{Lw} D.E. Littlewood,
{\it The Theory of Group Characters and Matrix Representations of
Groups}, Oxford University Press, New York 1940.

\bibitem[M1]{M1} G.E. Murphy, {\it On the representation theory of
the symmetric groups and associated Hecke algebras}, J. Algebra
{\bf 152} (1992), 492-513.

\bibitem[M2]{M2} G.E. Murphy, {\it The representations of Hecke algebras
of type $A_n$}, J. Algebra {\bf 173} (1995), 97-121.

\bibitem[M3]{M3} G.E. Murphy, {\it Representations of Hecke algebras of
$C_r\wr S_n$ at roots of unity}, in preparation (1996).

\bibitem[Ma]{Ma} P. Martin, {\it The structure of the partition algebras},
J. Algebra {\bf 183} (1996), 319-358.

\bibitem[Mac]{Mac}   I. G. Macdonald,
Symmetric Functions and Hall Polynomials (2nd Edition)
Oxford Univ. Press,  1995.

\bibitem[Mal]{Mal} G. Malle, 
{\it Unipotente Grade imprimitiver komplexer Spiegelungsgruppen},
J. Algebra {\bf 177} (1995), 768--826.

\bibitem[Mth]{Mth} A. Mathas, {\it Canonical bases and the decomposition
matrices of Ariki-Koike algebras}, to appear in J. Algebra.

\bibitem[Mu1]{Mu1} J. Murakami, {\it The Kauffman polynomial of links
and representation theory}, Osaka J. Math. {\bf 24} (1987), 
745-758.

\bibitem[Mu2]{Mu2} J. Murakami, {\it The representations of the
$q$-analogue of Brauer's centralizer algebras and the 
Kauffman polynomial of links}, Publ. Res. Inst. Math. Sci.
{\bf 26} (1990), 935-945.

\bibitem[Mur]{Mur}   F.D.~Murnaghan,
{\it The characters of the symmetric group}, Amer. J. Math. {\bf 59}
(1937),  739-753.

\bibitem[Nak]{Nak}   T. Nakayama, {\it
On some modular properties of irreducible representations
of a symmetric group I and II}, Jap. J. Math. {\bf 17} (1940),
165-184, 411-423.

\bibitem[Nz]{Nz} M. Nazarov, {\it Young's orthogonal form for Brauer's
centralizer algebra}, J. Algebra {\bf 182} (1996), 664-693.

\bibitem[OS]{OS} P. Orlik and L. Solomon,
{\it Unitary reflection groups and cohomology}, 
Invent.  Math. {\bf 59} (1980), 77--94.

\bibitem[Osi]{Osi}   M. Osima,
{\it On the representations of the generalized symmetric group},
Math. J. Okayama Univ. {\bf 4} (1954), 39-56.

\bibitem[P1]{P1} R. Penrose, {\it Angular momentum: an approach to combinatorial
space-time},  Quantum theory and beyond (T.A. Bastin, eds.),
Cambridge Univ. Press, 1969.

\bibitem[P2]{P2} R. Penrose, {\it Applications of negative dimensional tensors},
Combinatorial mathematics and its applications (D.J.A. Welsh, eds.),
Academic Press, 1971.

\bibitem[Pa]{Pa} C. Pallikaros, {\it Representations of Hecke algebras of
type $D_n$}, J. Algebra {\bf 169} (1994), 20--48.

\bibitem[Pr]{Pr} C. Procesi, {\it The invariant theory of $n\times n$
matrices}, Adv. in Math. {\bf 19} (1976), 330-354.

\bibitem[Prc]{Prc} R. Proctor, {\it
Bruhat lattices, plane partition generating functions, 
and minuscule representations},
European J. Combin. {\bf 5} (1984), 331--350. 

\bibitem[Ra1]{Ra1} A. Ram, {\it Characters of Brauer's centralizer algebras},
Pacific J. Math. {\bf 169} (1995), 173-200.

\bibitem[Ra2]{Ra2} A. Ram, {\it A Forbenius formula for the characters
of the Hecke algebras}, Invent. Math. {\bf 106} (1991), 461-488.

\bibitem[Ra3]{Ra3} A. Ram, {\it Robinson-Schensted-Knuth insertion and
the characters of symmetric groups and Iwahori-Hecke algebras of
type A}, to appear in the special {\it Festschrift in honor of
Gian-Carlo Rota}, B. Sagan ed., Birkhauser, Boston, 1997.

\bibitem[Rem]{Rem} J. Remmel, {\it Formula for the expansion of the Kronecker
products $S_{(m,n)}\otimes S_{(1^{p-r}, r)}$ and $S_{(1^k, 2^l)} \otimes
S_{(1^{p-r},r)}$}, Discrete Math. {\bf 99} (1992), no 1-3, 265-287.

\bibitem[Res]{Res} N. Reshetikhin, {\it Quantized universal enveloping algebras,
the Yang-Baxter equation and invariants of links I}, LOMI Preprint no.
E-4-87 (1987).

\bibitem[Ro]{Ro}  Y. Roichman, 
{\it A recursive rule for Kazhdan-Lusztig characters},
to appear in Adv. in Math.

\bibitem[RTW]{RTW} G. Rumer, E. Teller, H. Weyl, {\it Eine fur die Valenztheorie
geeignete Basis ber bin\"aren Vektor-invarianten}, Nachr. Ges. Wiss. G\"ottingen
Math.-Phys. Kl., 1932, pp. 499-504.

\bibitem[Ru]{Ru} D.E. Rutherford,
{\it Substitutional Analysis}, Edinburgh, at the University
Press, 1948.

\bibitem[Sc1]{Sc1}   I. Schur, 
{\it \"Uber eine Klasse von Matrizen, die sich einer gegeben Matrix
zuordnen lassen}, Dissertation, 1901; (Reprinted in Gesammelte Abhandlungen
{\bf 1}, 1-72).

\bibitem[Sc2]{Sc2}   I. Schur, 
{\it \"Uber die rationalen Darstellungen der allgemeinen linearen
Gruppe}, Sitzungsberichte der K\"oniglich Preussisschen Akademie der
Wissenschaften zu Berlin (1927), 58-75; (Reprinted in Gesammelte
Abhandlungen {\bf 3}, 68-85).

\bibitem[Sch]{Sch}  M. Scheunert, {\it The theory of Lie superalgebras:
An introduction}, Lecture Notes in Mathematics {\bf 716},
Springer-Berlin, 1979.

\bibitem[Se1]{Se1}   J.-P. Serre, 
{\it Linear Representations of Finite Groups},
Springer-Verlag, New York, Berlin, 1977.

\bibitem[Se2]{Se2}   J.-P. Serre,
{\it Complex semisimple Lie algebras}, 
Translated from the French by G. A. Jones, 
Springer-Verlag, New York-Berlin, 1987.

\bibitem[Sh]{Sh}  I.R. Shafarevich, {\it Basic algebraic geometry 2:
Schemes and complex manifolds, Second edition},
Springer-Verlag, Berlin, 1996.

\bibitem[Shj]{Shj} T. Shoji, 
{\it On the Springer representations of the Weyl groups of classical
algebraic groups}, Comm. Algebra {\bf 7} (1979), 1713--1745.

\bibitem[Snv]{Snv} V. Serganova, 
{\it Kazhdan-Lusztig polynomials and character formula for the 
Lie superalgebra ${\goth{gl}}(m|n)$}, preprint 1997.


\bibitem[So1]{So1}  L. Solomon, {\it The Bruhat decomposition,
Tits system and Iwahori ring for the monoid of matrices
over a finite field}, Geom. Dedicata {\bf 36} (1990), 15-49.

\bibitem[So2]{So2} L. Solomon, Abstract No. 900-16-169,
Abstracts presented to the American Math. Soc., Vol. 16, No. 2,
Spring 1995.

\bibitem[Spc]{Spc} W. Specht, {\it Eine Verallgemeinerung der 
symmetrischen Gruppe}, Schriften Math. Seminar (Berlin) {\bf 1}
(1932), 1-32.

\bibitem[Spl]{Spl} N. Spaltenstein, {\it On the fixed point set of
a unipotent transformation on the flag manifold}, Proc. Kon.
Nederl. Akad. Wetensch. {\bf 79} (1976), 452-456.

\bibitem[Spr]{Spr}   T. Springer, 
{\it A Construction of Representations of Weyl Groups},
Invent. Math. {\bf 44} (1978), 279-293.

\bibitem[Srg]{Srg}  A.N. Sergeev, {\it Representations of the Lie superalgebras
${\goth{gl}}(n,m)$ and $Q(n)$ in a space of tensors},
Funct. Anal. App. {\bf 18} (1984), 80-81.

\bibitem[ST]{ST}   G.~C.~Shephard and J.~A.~Todd
{\it Finite unitary reflection groups},  Canad. J. Math.
{\bf 6} (1954),  274-304.

\bibitem[Sta1]{Sta1} R.P. Stanley,
{\it Theory and application of plane partitions. I, II},
Studies in Appl. Math. {\bf 50} (1971), 167-188, 259-279.

\bibitem[Sta2]{Sta2} R.P. Stanley,
{\it Some aspect of groups acting on finite posets},
J. Combin. Theory Ser. A {\bf 32} (1982), 132-161.

\bibitem[Ste1]{Ste1} J. Stembridge, 
{\it A combinatorial theory for rational actions of $GL_n$},
Contemp. Math., {\bf 88} (1989) 163-176.

\bibitem[Ste2]{Ste2} J. Stembridge,
{\it Rational tableaux and the tensor algebra of
${\goth{gl}}_n$}, J. Combin. Theory Ser. A {\bf 46}
(1987) 79-120.

\bibitem[Ste3]{Ste3}   J. Stembridge, 
{\it On the eigenvalues of representations of reflection groups and
wreath products},  Pacific J. Math. {\bf 140} (1989), 353-396.

\bibitem[Ste4]{Ste4} J. Stembridge, 
{\it On miniscule representations, plane partitions and involutions
in complex Lie groups}, Duke Math. J. {\bf 73} (1994), 469-490.

\bibitem[Ste5]{Ste5} J. Stembridge,
{\it  The enumeration of totally symmetric plane partitions}, Adv.
Math. {\bf 111} (1995), 227--243. 

\bibitem[Ste6]{Ste6} J. Stembridge,
{\it Some hidden relations involving the ten symmetry classes of
plane partitions}, J. Combin. Theory Ser. A {\bf 68} (1994), 372--409.

\bibitem[SU]{SU} Y. Shibukawa and K. Ueno, {\it Character table of
Hecke algebra of type $A_{N-1}$ and representations of the quantum
group $U_q({\goth{gl}_{n+1})}$}, in {\sl Infinite Analysis} (Kyoto,1991),
Adv. Ser. Math. Phys. {\bf 16}, World Sci. Publishing,
River Edge, NJ, (1992), 977-984.

\bibitem[Su1]{Su1} S. Sundaram, 
{\it Tableaux in the representation theory of the classical Lie groups},
Invariant theory and tableaux'', IMA Vol. Math. Appl. {\bf 19}
Springer, New York (1990), 191-225.

\bibitem[Su2]{Su2} S. Sundaram, 
{\it Orthogonal tableaux and an insertion algorithm for
$SO(2n+1)$}, J. Combin. Theory Ser. A {\bf 53} (1990), 239-256.

\bibitem[Su3]{Su3} S. Sundaram,
{\it The Cauchy identity for $Sp(2n)$}, J. Combin.
Theory Ser. A {\bf 53} (1990), 209-238.

\bibitem[Ta]{Ta} T. Tanisaki, {\it Defining ideals of the closures os
conjugacy classes and representations of Weyl groups},
Tohoku Math. J. {\bf 34} (1982), 575-585.

\bibitem[TL]{TL}   N.Temperley, E. Lieb,
{\it Relations between the ``percolation'' and ``colouring'' problem
and other graph theoretical problems associated with regular plane lattices:
Some exact results for the ``percolation'' problem}, Proc. Roy. Soc. London
Ser. A, {\bf 322} (1971), 547-597.

\bibitem[Va]{Va} V.S. Varadarajan, {\it Lie groups, Lie algebras,
and Their Representations}, Graduate Texts in Mathematics {\bf 102},
Springer-Verlag, New York-Berlin, 1984.

\bibitem[vdJ]{vdJ} J. van der Jeugt, {\it An algorithm for 
characters of Hecke algebras $H\sb n(q)$ of type $A\sb {n-1}$},
J. Phys. A {\bf 24} (1991), 3719--3725. 

\bibitem[Vg1]{Vg1}
D. Vogan, {\it Unitary representations of reductive Lie groups}, 
Annals of Mathematics Studies {\bf 118}, 
Princeton University Press, Princeton, NJ, 1987.

\bibitem[Vg2]{Vg2} D. Vogan, 
{\it Representations of real reductive Lie groups}, Progress in
Mathematics {\bf 15}, Birkhäuser, Boston, Mass., 1981.

\bibitem[VK]{VK} A. Vershik and S. Kerov,
{\it Characters and realizations of representations of the
infinite-dimensional Hecke algebra, and knot invariants},
Soviet Math. Dokl. {\bf 38}, (1989), 134--137.

\bibitem[Wa]{Wa} N. Wallach, {\it Real Reductive Groups I, II},
Pure and Applied Mathematics, {\bf 132},
Academic Press, Inc., Boston, MA, 1988, 1992.

\bibitem[Wm]{Wm} J. Weyman,
{\it The equations of conjugacy classes of nilpotent matrices},
Invent.  Math. {\bf 98} (1989), 229--245.

\bibitem[Wy1]{Wy1} H. Weyl, {\it The classical groups, their
invariants and their representations}, Princeton University Press, 
Princeton, 1946.

\bibitem[Wy2]{Wy2} H. Weyl, {\it Theorie der Darstellung kontinuerlicher
halb-einfacher Gruppen durch lineare Transformationen I, II, III},
Mathematische Zeitschrift {\bf 23, 24}, (1925-26), 271-309, 328-376,
377-395.

\bibitem[Wz1]{Wz1}   H.  Wenzl,
{\it Hecke algebras of type $A_n$ and subfactors},
Invent. Math. {\bf 92} (1988), no. 2, 349-383. 

\bibitem[Wz2]{Wz2} H. Wenzl,
{\it On the structure of Brauer's centralizer algebras},
Ann. Math. {\bf 128} (1988), 173-193.

\bibitem[Wz3]{Wz3}   H.  Wenzl,
{\it Quantum groups and subfactors of type $B$, $C$, and $D$},
Commun. Math. Phys. {\bf 133} (1990), 383-432.

\bibitem[Y]{Y}   A.~Young,
{\it On Quantitative substitutional analysis I-IX},
Proc. London Math. Soc. (1901-1952), reprinted in
{\it The Collected Papers of Alfred Young 1873-1940},
Mathematical Exposition No. {\bf 21} U. of Toronto Press, (1977).

\bibitem[Y1]{Y1}   A.~Young,
{\it On Quantitative Substitutional Analysis} (First paper)
Proc. London Math. Soc. {\bf 33} (1900), 97-146.

\bibitem[Y2]{Y2}   A.~Young,
{\it On Quantitative Substitutional Analysis} (Second paper)
Proc. London Math. Soc. {\bf 34} (1902), 361-397.

\bibitem[Y3]{Y3}   A.~Young,
{\it On Quantitative Substitutional Analysis} (Fifth paper)
Proc. London Math. Soc. (2) {\bf 31} (1930), 273-288.

\bibitem[Y4]{Y4}   A.~Young,
{\it On Quantitative Substitutional Analysis} (Sixth paper)
Proc. London Math. Soc. (2) {\bf 34} (1931), 196-230.

\bibitem[Y5]{Y5}   A.~Young,
{\it On Quantitative Substitutional Analysis} (Eighth paper)
Proc. London Math. Soc. (2) {\bf 37} (1934), 441-495.

\bibitem[Zh1]{Zh1} D.P. Zhelobenko, 
{\it An analogue of the Gelfand-Tsetlin basis for symplectic Lie algebras},
Russian Math. Surveys {\bf 42} (1987), 247--248.

\bibitem[Zh2]{Zh2} D.P. Zhelobenko, 
{\it Compact Lie groups and their representations},
Translations of Mathematical Monographs v. 40,
Amer. Math. Soc. 1973.






\end{thebibliography}
\end{document}